\makeatletter

\newdimen\paperwidth
\newdimen\paperheight

\def\papersize#1#2{\let\p@persize\relax\paperwidth#1\paperheight#2}

\def\Afour{\papersize{210truemm}{297truemm}}

\let\p@persize\Afour

\let\onesidestyle\@twosidefalse
\let\twosidestyle\@twosidetrue

\def\margins{\@ifnextchar[{\@margins}{\@margins[\z@]}}

\def\@margins[#1]#2#3{
  \p@persize\dimen0 #3\dimen0 .5\dimen0\normalsize%
  \oddsidemargin-1truein\advance\oddsidemargin#2%
  \evensidemargin-1truein\advance\evensidemargin#2%
  \topmargin-1truein\advance\topmargin\dimen0\headsep\dimen0\footskip\dimen0%
  \textwidth\paperwidth\advance\textwidth-#2\advance\textwidth-#2%
  \textheight\paperheight\advance\textheight-#3\advance\textheight-#3%
  \headheight\baselineskip\advance\topmargin-.5\baselineskip%
  \advance\headsep-.5\baselineskip%
  \footheight\baselineskip
  \advance\textwidth-#1\advance\oddsidemargin#1
  \if@twoside\def\@themargin%
    {\ifodd\count\z@\oddsidemargin\else\evensidemargin\fi}\fi}

\def\headlinesep#1{\advance\topmargin\headsep\advance\topmargin -#1
  \advance\topmargin.5\baselineskip\headsep #1\advance\headsep-.5\baselineskip}

\def\headline{\if@twoside\let\n@xt\h@dlin@\else\let\n@xt\h@@dlin@\fi\n@xt}
  
\def\h@dlin@#1#2{%
  \def\@oddhead{%
    {{\leftskip\z@\rightskip\z@\noindent\normalsize#1}}}
  \def\@evenhead{%
    {{\leftskip\z@\rightskip\z@\noindent\normalsize#2}}}}

\def\h@@dlin@#1{%
  \def\@oddhead{{{\leftskip\z@\rightskip\z@\noindent\normalsize#1}}}}

\def\footline{\if@twoside\let\n@xt\f@tlin@\else\let\n@xt\f@@tlin@\fi\n@xt}

\def\f@tlin@#1#2{%
  \def\@oddfoot{%
    {{\leftskip\z@\rightskip\z@\noindent\normalsize#1}}}
  \def\@evenfoot{%
    {{\leftskip\z@\rightskip\z@\noindent\normalsize#2}}}}

\def\f@@tlin@#1{%
  \def\@oddfoot{{{\leftskip\z@\rightskip\z@\noindent\normalsize#1}}}}

\def\normalpage{\global\@specialpagefalse}

\makeatother

\makeatletter

\def\ft{\@ifnextchar[{\ft@s}{\ft@}}
\def\ft@{\ft@@@s[\f@size]}
\def\ft@s[{\@ifnextchar{a}{\ft@sz[}{\ft@@s[}}
\def\ft@@s[{\@ifnextchar{s}{\ft@sz[}{\ft@@@s[}}
\def\ft@@@s[#1]{\ft@sz[at #1pt]}
\def\ft@sz[#1]#2{\font\fonttemp=#2 #1\fonttemp\ignorespaces}

\makeatother

\makeatletter

\input{epsf.sty}

\def\showfig#1#2{\epsfbox{#2}}
\def\fig@#1#2{\leavevmode{\framebox{\figstyl@\strut{ #1 }}}}

\def\figstyle#1{\def\figstyl@{#1}}

\figstyle{\shape{n}\series{m}\selectfont\normalsize}

\def\showfigurestrue{\let\fig\showfig}
\def\showfiguresfalse{\let\fig\fig@}

\makeatother

\showfiguresfalse

\def\smallcircc{\mathop{\mkern3.5mu\hbox{\raise.58ex\hbox{\ft{lcircle10}a}}}}
\def\varemptyset{{\text{\raise.21ex\hbox{$\not$}}\mkern.15mu\mathrm{O}\mkern.15mu}}

  \let\epsilon\varepsilon
      \let\theta\vartheta
          \let\phi\varphi

\documentstyle[12pt]{article}

\makeatletter

\let\Larg@\Large
\let\hug@\huge

\def\usepackage#1{\input{#1.sty}}

\input{geom.sty}

\def\r@adlabel#1#2{\global\@namedef{#1@\the\@key}{#2}}

\let\Large\Larg@
\let\huge\hug@

\def\smallskip{\vskip\smallskipamount}
\def\medskip{\vskip\medskipamount}
\def\bigskip{\vskip\bigskipamount}

\def\mytrivlist{\parsep\parskip\@nmbrlistfalse
  \my@trivlist \labelwidth\z@ \leftmargin\z@
  \itemindent\z@ \def\makelabel##1{##1}}

\def\my@trivlist{\global\@newlisttrue \@outerparskip\parskip}

\def\end#1{\csname end#1\endcsname\@checkend{#1}%
  \expandafter\endgroup\if@endpe\@doendpe\fi
  \if@ignore \global\@ignorefalse \ignorespaces\fi}

\def\put{\@ifnextchar[{\@put}{\@@rput[\z@,\z@][r]}}
\def\@put[#1]{\@ifnextchar[{\@@put[#1]}{\@@@@@put[#1]}}
\def\@@put[#1][{\@ifnextchar{l}{\@@lput[#1][}{\@@@put[#1][}}
\def\@@@put[#1][{\@ifnextchar{c}{\@@cput[#1][}{\@@@@put[#1][}}
\def\@@@@put[#1][{\@ifnextchar{r}{\@@rput[#1][}{\relax}}
\def\@@@@@put[{\@ifnextchar{l}{\@@lput[\z@,\z@][}{\@@@@@@put[}}
\def\@@@@@@put[{\@ifnextchar{c}{\@@cput[\z@,\z@][}{\@@@@@@@put[}}
\def\@@@@@@@put[{\@ifnextchar{r}{\@@rput[\z@,\z@][}{\@@@@@@@@put[}}
\def\@@@@@@@@put[#1]{\@@rput[#1][r]}

\let\hm@d@\leavevmode

\long\def\@@lput[#1,#2][l]#3{\setbox0\hbox{#3}\hm@d@\raise#2\hbox to\z@{\dimen0 #1%
  \advance\dimen0-\wd0\kern\dimen0\dp0\z@\ht0\z@\wd0\z@\box0\hss}\ignorespaces}
\long\def\@@cput[#1,#2][c]#3{\setbox0\hbox{#3}\hm@d@\raise#2\hbox to\z@{\dimen0 #1%
  \advance\dimen0-.5\wd0\kern\dimen0\dp0\z@\ht0\z@\wd0\z@\box0\hss}\ignorespaces}
\long\def\@@rput[#1,#2][r]#3{\setbox0\hbox{\kern#1\raise#2\hbox{#3}}%
  \dp0\z@\ht0\z@\wd0\z@\hm@d@\box0\ignorespaces}

\def\flbox{\@ifnextchar[{\@flbox}{\@@rflbox[\z@,\z@][r]}}
\def\@flbox[#1]{\@ifnextchar[{\@@flbox[#1]}{\@@@@@flbox[#1]}}
\def\@@flbox[#1][{\@ifnextchar{l}{\@@lflbox[#1][}{\@@@flbox[#1][}}
\def\@@@flbox[#1][{\@ifnextchar{c}{\@@cflbox[#1][}{\@@@@flbox[#1][}}
\def\@@@@flbox[#1][{\@ifnextchar{r}{\@@rflbox[#1][}{\relax}}
\def\@@@@@flbox[{\@ifnextchar{l}{\@@lflbox[\z@,\z@][}{\@@@@@@flbox[}}
\def\@@@@@@flbox[{\@ifnextchar{c}{\@@cflbox[\z@,\z@][}{\@@@@@@@flbox[}}
\def\@@@@@@@flbox[{\@ifnextchar{r}{\@@rflbox[\z@,\z@][}{\@@@@@@@@flbox[}}
\def\@@@@@@@@flbox[#1]{\@@rflbox[#1][r]}
\long\def\@@lflbox[#1,#2][l]#3{\@@lput[#1,#2][l]{%
  \vtop{\leftskip\z@\parindent\z@\raggedleft\hm@d@#3}}}
\long\def\@@cflbox[#1,#2][c]#3{\@@cput[#1,#2][c]{%
  \vtop{\leftskip\z@\parindent\z@\raggedcenter\hm@d@#3}}}
\long\def\@@rflbox[#1,#2][r]#3{\@@rput[#1,#2][r]{%
  \vtop{\leftskip\z@\parindent\z@\raggedright\hm@d@#3}}}

\def\maketitle{\par
 \begingroup
 \def\thefootnote{\fnsymbol{footnote}}
 \def\@makefnmark{\hbox to 0pt{$^{\@thefnmark}$\hss}} 
 \if@twocolumn 
 \twocolumn[\@maketitle] 
 \else 
 \global\@topnum\z@ \@maketitle \fi\thispagestyle{plain}\@thanks
 \endgroup
 \setcounter{footnote}{0}
 \let\maketitle\relax
 \let\@maketitle\relax
 \gdef\@thanks{}\gdef\@author{}\gdef\@title{}\let\thanks\relax}

\def\@maketitle{ 
 \null
 \vskip 2em \begin{center}
 {\LARGE \@title \par} \vskip 1.5em {\large \lineskip .5em
\begin{tabular}[t]{c}\@author 
 \end{tabular}\par} 
 \vskip 1em {\large \@date} \end{center}
 \par
 \vskip 1.5em}
 
\def\partbeforeskip#1{\def\p@rtbeforeskip{#1}}
\def\partstyle#1{\def\p@rtstyl@{#1}}
\def\partdot#1{\def\partd@t{#1}}
\def\partafterskip#1{\def\p@rtafterskip{#1}}
\def\partintrostyle#1{\def\partintr@styl@{#1}}
\def\partintrodot#1{\def\partintr@dot{#1}}
\long\def\partintrosep#1{\long\def\partintr@sep{#1}}
\def\partnewpagetrue{\def\p@rtnewp@ge{\newpage}}
\def\partnewpagefalse{\long\def\p@rtnewp@ge{\par}}

\partbeforeskip{4ex}
\partstyle{\centering\Large\bf}
\partdot{}
\partafterskip{3ex}
\partintrostyle{\large}
\partintrodot{}
\partintrosep{\par}
\partnewpagefalse

\def\partname{Part}
\def\part{\p@rtnewp@ge\addvspace\p@rtbeforeskip\@afterindentfalse\secdef\@part\@spart}

\def\@part[#1]#2{\ifnum \c@secnumdepth >-1\relax  
        \refstepcounter{part}                     
        \def\@tempa{\addcontentsline{toc}{part}}  %
        \expandafter\@tempa\expandafter{\thepart  
          \hspace{1em}#1}\else                    
        \addcontentsline{toc}{part}{#1}\fi        
   {\p@rtstyl@                       
    \ifnum \c@secnumdepth >-1\relax        
      {\partintr@styl@\partname\ \thepart  
       \partintr@dot}\partintr@sep\nobreak 
    \fi                                    
    #2\partd@t\markboth{}{}\par}
    \nobreak                       
    \vskip\p@rtafterskip           
   \@afterheading                  
    }                              

\def\@spart#1{{\p@rtcentering\p@rtstyl@                      
    #1\partd@t\par}                 
    \nobreak                        
    \vskip\p@rtafterskip            
    \@afterheading                  
  }                                 

\newif\ifsection@ftind
\newif\ifsection@ftpar

\def\sectionbeforeskip#1{\def\s@ctbeforeskip{#1}}
\def\sectionstyle#1{\def\s@ctstyl@{#1}}
\def\sectiondot#1{\def\sectiond@t{#1}}
\def\sectionafterskip#1{\def\s@ctafterskip{#1}}
\def\sectionintrostyle#1{\def\sectionintr@styl@{#1}}
\def\sectionintro#1{\def\sectionintr@{#1}}
\def\sectionintrodot#1{\def\sectionintr@dot{#1}}
\def\sectionintrosep#1{\def\sectionintr@sep{#1}}
\def\sectionindenttrue{\def\s@ctind{\parindent}}
\def\sectionindentfalse{\def\s@ctind{\z@}}
\def\sectionafterindenttrue{\section@ftindtrue}
\def\sectionafterindentfalse{\section@ftindfalse}
\def\sectionafternewlinetrue{\section@ftpartrue}
\def\sectionafternewlinefalse{\section@ftparfalse}

\newif\ifsubsection@ftind
\newif\ifsubsection@ftpar

\def\subsectionbeforeskip#1{\def\ss@ctbeforeskip{#1}}
\def\subsectionstyle#1{\def\ss@ctstyl@{#1}}
\def\subsectiondot#1{\def\subsectiond@t{#1}}
\def\subsectionafterskip#1{\def\ss@ctafterskip{#1}}
\def\subsectionintrostyle#1{\def\subsectionintr@styl@{#1}}
\def\subsectionintro#1{\def\subsectionintr@{#1}}
\def\subsectionintrodot#1{\def\subsectionintr@dot{#1}}
\def\subsectionintrosep#1{\def\subsectionintr@sep{#1}}
\def\subsectionindenttrue{\def\ss@ctind{\parindent}}
\def\subsectionindentfalse{\def\ss@ctind{\z@}}
\def\subsectionafterindenttrue{\subsection@ftindtrue}
\def\subsectionafterindentfalse{\subsection@ftindfalse}
\def\subsectionafternewlinetrue{\subsection@ftpartrue}
\def\subsectionafternewlinefalse{\subsection@ftparfalse}

\newif\ifsubsubsection@ftind
\newif\ifsubsubsection@ftpar

\def\subsubsectionbeforeskip#1{\def\sss@ctbeforeskip{#1}}
\def\subsubsectionstyle#1{\def\sss@ctstyl@{#1}}
\def\subsubsectiondot#1{\def\subsubsectiond@t{#1}}
\def\subsubsectionafterskip#1{\def\sss@ctafterskip{#1}}
\def\subsubsectionintrostyle#1{\def\subsubsectionintr@styl@{#1}}
\def\subsubsectionintro#1{\def\subsubsectionintr@{#1}}
\def\subsubsectionintrodot#1{\def\subsubsectionintr@dot{#1}}
\def\subsubsectionintrosep#1{\def\subsubsectionintr@sep{#1}}
\def\subsubsectionindenttrue{\def\sss@ctind{\parindent}}
\def\subsubsectionindentfalse{\def\sss@ctind{\z@}}
\def\subsubsectionafterindenttrue{\subsubsection@ftindtrue}
\def\subsubsectionafterindentfalse{\subsubsection@ftindfalse}
\def\subsubsectionafternewlinetrue{\subsubsection@ftpartrue}
\def\subsubsectionafternewlinefalse{\subsubsection@ftparfalse}

\newif\ifparagraph@ftind
\newif\ifparagraph@ftpar

\def\paragraphbeforeskip#1{\def\p@rbeforeskip{#1}}
\def\paragraphstyle#1{\def\p@rstyl@{#1}}
\def\paragraphdot#1{\def\paragraphd@t{#1}}
\def\paragraphafterskip#1{\def\p@rafterskip{#1}}
\def\paragraphintrostyle#1{\def\paragraphintr@styl@{#1}}
\def\paragraphintro#1{\def\paragraphintr@{#1}}
\def\paragraphintrodot#1{\def\paragraphintr@dot{#1}}
\def\paragraphintrosep#1{\def\paragraphintr@sep{#1}}
\def\paragraphindenttrue{\def\p@rind{\parindent}}
\def\paragraphindentfalse{\def\p@rind{\z@}}
\def\paragraphafterindenttrue{\paragraph@ftindtrue}
\def\paragraphafterindentfalse{\paragraph@ftindfalse}
\def\paragraphafternewlinetrue{\paragraph@ftpartrue}
\def\paragraphafternewlinefalse{\paragraph@ftparfalse}

\newif\ifsubparagraph@ftind
\newif\ifsubparagraph@ftpar

\def\subparagraphbeforeskip#1{\def\sp@rbeforeskip{#1}}
\def\subparagraphstyle#1{\def\sp@rstyl@{#1}}
\def\subparagraphdot#1{\def\subparagraphd@t{#1}}
\def\subparagraphafterskip#1{\def\sp@rafterskip{#1}}
\def\subparagraphintrostyle#1{\def\subparagraphintr@styl@{#1}}
\def\subparagraphintro#1{\def\subparagraphintr@{#1}}
\def\subparagraphintrodot#1{\def\subparagraphintr@dot{#1}}
\def\subparagraphintrosep#1{\def\subparagraphintr@sep{#1}}
\def\subparagraphindenttrue{\def\sp@rind{\parindent}}
\def\subparagraphindentfalse{\def\sp@rind{\z@}}
\def\subparagraphafterindenttrue{\subparagraph@ftindtrue}
\def\subparagraphafterindentfalse{\subparagraph@ftindfalse}
\def\subparagraphafternewlinetrue{\subparagraph@ftpartrue}
\def\subparagraphafternewlinefalse{\subparagraph@ftparfalse}

\sectionbeforeskip{\bigskipamount}
\sectionstyle{\large\bf}
\sectiondot{}
\sectionafterskip{.5\bigskipamount}
\sectionintrostyle{}
\sectionintro{}
\sectionintrodot{.}
\sectionintrosep{1.25ex}
\sectionindentfalse
\sectionafterindenttrue
\sectionafternewlinetrue

\subsectionbeforeskip{.8\bigskipamount}
\subsectionstyle{\normalsize\bf}
\subsectiondot{}
\subsectionafterskip{.4\bigskipamount}
\subsectionintrostyle{}
\subsectionintro{}
\subsectionintrodot{.}
\subsectionintrosep{1.25ex}
\subsectionindentfalse
\subsectionafterindenttrue
\subsectionafternewlinetrue

\subsubsectionbeforeskip{.6\bigskipamount}
\subsubsectionstyle{\normalsize\bf}
\subsubsectiondot{}
\subsubsectionafterskip{.3\bigskipamount}
\subsubsectionintrostyle{}
\subsubsectionintro{}
\subsubsectionintrodot{.}
\subsubsectionintrosep{1.25ex}
\subsubsectionindentfalse
\subsubsectionafterindenttrue
\subsubsectionafternewlinetrue

\paragraphbeforeskip{.5\bigskipamount}
\paragraphstyle{\normalsize\bf}
\paragraphdot{.}
\paragraphafterskip{1.25ex}
\paragraphintrostyle{}
\paragraphintro{}
\paragraphintrodot{.}
\paragraphintrosep{1.25ex}
\paragraphindentfalse
\paragraphafterindenttrue
\paragraphafternewlinefalse

\subparagraphbeforeskip{.5\bigskipamount}
\subparagraphstyle{\normalsize\bf}
\subparagraphdot{.}
\subparagraphafterskip{1.25ex}
\subparagraphintrostyle{}
\subparagraphintro{}
\subparagraphintrodot{.}
\subparagraphintrosep{1.25ex}
\subparagraphindenttrue
\subparagraphafterindenttrue
\subparagraphafternewlinefalse

\let\@partoken\par
\long\def\@@gobble#1{}
\def\ignorepar{\@ifnextchar\@partoken{\expandafter\ignorepar\@@gobble}{\ignorespaces}}

\def\@startsection#1#2#3#4#5#6{
   \@tempskipa #4\relax
   \csname if#1@ftind\endcsname\@afterindenttrue\else\@afterindentfalse\fi
   \advance\@tempskipa by\presection
   \if@nobreak \everypar{}\else
     \addpenalty{\@secpenalty}\addvspace{\@tempskipa}%
     \allowbreak\vskip -\presection \fi \@ifstar
     {\@ssect{#1}{#2}{#3}{#4}{#5}{#6}}{\@dblarg{\@sect{#1}{#2}{#3}{#4}{#5}{#6}}}}

\def\@sect#1#2#3#4#5#6[#7]#8{\def\object@type{#1}%
   \ifnum #2>\c@secnumdepth\def\@svsec{}\def\@tempb{}%
      \else\refstepcounter{#1}\def\@svsec{{\csname #1intr@styl@\endcsname%
        {\csname #1intr@\endcsname}\csname the#1\endcsname%
        \csname #1intr@dot\endcsname\kern\csname #1intr@sep\endcsname}}%
        \edef\@tempb{\noexpand\numberline{\csname the#1\endcsname}}\fi%
   \def\@tempa{\addcontentsline{toc}{#1}}%
   \csname if#1@ftpar\endcsname%
      \begingroup #6\relax%
        \@hangfrom{\hskip #3\relax\@svsec}{\interlinepenalty \@M{#8}%
        \csname #1d@t\endcsname\par}%
      \endgroup%
      \csname #1mark\endcsname{#7}%
      \expandafter\@tempa\expandafter{\@tempb #7}%
      \ifautolabel\label*{#8}\fi%
   \else%
      \def\@svsechd{#6\hskip #3\relax%
         \@svsec{#8}\csname #1mark\endcsname{#7}%
         \expandafter\@tempa\expandafter{\@tempb #7}%
         \ifautolabel\label*{#8}\fi}\fi%
   \@xsect{#1}{#5}\ignorepar}

\def\@ssect#1#2#3#4#5#6#7{%
   \ifnum #2>\c@secnumdepth\def\@tempb{}\else \def\@tempb{\numberline{}}\fi%
     \def\@tempa{\addcontentsline{toc}{s#1}}%
     \csname if#1@ftpar\endcsname
        \begingroup #6\relax
           \@hangfrom{\hskip #3}{\interlinepenalty \@M{#7}%
           \csname #1d@t\endcsname\par}%
        \endgroup
        \csname s#1mark\endcsname{#7}%
        \ifstarredcontents\expandafter\@tempa\expandafter{\@tempb #7}\fi%
        \ifautolabel\label*{#7}\fi%
     \else%
        \def\@svsechd{#6\hskip #3\relax{#7}\csname s#1mark\endcsname{#7}%
        \ifautolabel\label*{#7}\fi}\fi
   \@xsect{#1}{#5}\ignorepar}

\def\@xsect#1#2{
   \csname if#1@ftpar\endcsname 
       \par \nobreak \vskip #2\relax \@afterheading
    \else \global\@nobreakfalse \global\@noskipsectrue
       \everypar{\if@noskipsec \global\@noskipsecfalse
                   \clubpenalty\@M \hskip -\parindent
                   \begingroup \@svsechd \endgroup \unskip
                   \hskip #2\relax  
                  \else \clubpenalty \@clubpenalty
                    \everypar{}\fi}\fi\ignorespaces}

\def\section{\@startsection{section}{1}{\s@ctind}
  {\s@ctbeforeskip}{\s@ctafterskip}{\s@ctstyl@}}
\def\subsection{\@startsection{subsection}{2}{\ss@ctind}
  {\ss@ctbeforeskip}{\ss@ctafterskip}{\ss@ctstyl@}}
\def\subsubsection{\@startsection{subsubsection}{3}{\sss@ctind}
  {\sss@ctbeforeskip}{\sss@ctafterskip}{\sss@ctstyl@}}
\def\paragraph{\@startsection{paragraph}{4}{\p@rind}
  {\p@rbeforeskip}{\p@rafterskip}{\p@rstyl@}}
\def\subparagraph{\@startsection{subparagraph}{4}{\sp@rind}
  {\sp@rbeforeskip}{\sp@rafterskip}{\sp@rstyl@}}

\def\statementabove#1{\def\th@bove{#1}}
\def\statementstyle#1{\def\thstyl@{#1}}
\def\statementbelow#1{\def\thb@low{#1}}
\def\statementindentfalse{\let\thind@nt\relax}
\def\statementindenttrue{\let\thind@nt\indent}

\def\statementintrostyle#1{\def\thintr@style{#1}}
\def\statementintrodot#1{\def\thintr@dot{#1}}
\def\statementintrosep#1{\def\thintr@sep{#1}}
\def\statementintrobrackets#1#2{\def\thintr@left{#1}\def\thintr@right{#2}}

\statementabove{\medskip}
\statementstyle{\sl}
\statementbelow{\medskip}
\statementindenttrue

\statementintrostyle{\normalshape\bf}
\statementintrodot{.}
\statementintrosep{\kern1.25ex}
\statementintrobrackets{(}{)}

\def\@thskip{\dimen0\lastskip\vskip-\dimen0%
  \th@bove\dimen1\lastskip\vskip-\dimen1%
  \ifdim\dimen0>\dimen1\else\dimen0\dimen1\fi\vskip\dimen0}

\long\def\@@newtheorem#1#2#3{%
  \newenvironment{#3}%
    {\def\object@type{#3}\par\@thskip%
     \@ifnextchar[{\@enva{#3}{\thstyl@#1{#2}}}{\@envb{#3}{\thstyl@#1{#2}}}}%
    {\end{#3@}}%
  \@ifnextchar[{\@nothm{#3}}{\@nnthm{#3}}}

\def\@nothm#1[#2]#3{%
  \@ifundefined{c@#2}{\@latexerr{No theorem environment `#2' defined}\@eha}%
  {\expandafter\@ifdefinable\csname #1@\endcsname
  {\global\@namedef{the#1}{\@nameuse{the#2}}%
   \global\@namedef{c@#1}{\@nameuse{c@#2}}
   \global\@namedef{p@#1}{\@nameuse{p@#2}}
   \global\@namedef{#1@}{\@nnnthm{#2}{#3}}%
   \global\@namedef{end#1@}{\@endtheorem}}}}

\def\@nnnthm#1#2{\refstepcounter
    {#1}\@ifnextchar[{\@ynnnthm{#1}{#2}}{\@xnnnthm{#1}{#2}}}
 
\def\@xnnnthm#1#2{\@begintheorem{#2}{\csname the#1\endcsname}\ignorespaces}
\def\@ynnnthm#1#2[#3]{\@opargbegintheorem{#2}{\csname
       the#1@\endcsname}{#3}\ignorespaces}

\def\renewtheorem{\@ifnextchar[{\@renewtheorem}{\@renewtheorem[{}{}]}}

\long\def\@renewtheorem[#1]{\@@renewtheorem#1}

\long\def\@@renewtheorem#1#2#3{%
  \expandafter\let\csname#3@\endcsname\undefined
  \renewenvironment{#3}%
    {\def\object@type{#3}\par\@thskip%
     \@ifnextchar[{\@enva{#3}{\thstyl@#1{#2}}}{\@envb{#3}{\thstyl@#1{#2}}}}%
    {\end{#3@}}%
  \@ifnextchar[{\@nothm{#3}}{\@nnthm{#3}}}

\def\@begintheorem#1#2{\@opargbegintheorem{#1}{#2}{}}

\def\@opargbegintheorem#1#2#3{%
        \def\@tempx{#1}%
        \expandafter\let\expandafter\@tempy#2
        \def\@tempz{#3}%
        \mytrivlist\item[\thind@nt\hskip\labelsep%
        {\thintr@style%
          #1\if\@tempx\@empty\else\if\@tempy\relax\else\kern1ex\fi\fi#2%
          \ifx\@tempz\@empty%
            \if\@tempx\@empty\if\@tempy\relax%
            \else\thintr@dot\thintr@sep\fi\else\thintr@dot\thintr@sep\fi%
            \else%
            \if\@tempx\@empty\if\@tempy\relax\else\kern1ex\fi\else\kern1ex\fi%
           \thintr@left{#3}\thintr@right\thintr@dot\thintr@sep\fi}%
            \hskip-\labelsep]%
        \ifautolabel\label*{#3}\fi}

\def\@endtheorem{\strut\endtrivlist\thb@low}

\def\proofabove#1{\def\pf@bove{#1}}
\def\proofstyle#1{\def\pfstyl@{#1}}
\def\proofbelow#1{\def\pfb@low{#1}}
\def\proofindentfalse{\let\pfind@nt\relax}
\def\proofindenttrue{\let\pfind@nt\indent}

\def\proofintrostyle#1{\def\pfintr@style{#1}}
\def\proofintrodot#1{\def\pfintr@dot{#1}}
\def\proofintrosep#1{\def\pfintr@sep{#1}}
\def\proofintrobrackets#1#2{\def\pfintr@left{#1}\def\pfintr@right{#2}}

\proofabove{\medskip}
\proofstyle{}
\proofbelow{\medskip}
\proofindenttrue

\proofintrostyle{\sl}
\proofintrodot{.}
\proofintrosep{\kern1.25ex}
\proofintrobrackets{of\kern1ex}{}

\def\@pfskip{\dimen0\lastskip\vskip-\dimen0%
  \pf@bove\dimen1\lastskip\vskip-\dimen1%
  \ifdim\dimen0>\dimen1\else\dimen0\dimen1\fi\vskip\dimen0}

\renewenvironment{proof}%
  {\@pfskip\mytrivlist\item[\pfind@nt]\@ifnextchar[{\pro@f}{\pro@f[\prooftag]}}
  {\ifvoid\provedbox\else\hproved\fi\endtrivlist\pfb@low}

\def\pro@f[#1]{\setbox\provedbox\hbox{\provedboxcontents{#1}}\proofintro{#1}}

\def\proofintro#1{\expandafter\def\expandafter\@tempa\expandafter{#1}%
  {\pfintr@style{Proof\ifx\@tempa\empty\else\kern1ex\pfintr@left{#1}%
  \pfintr@right\fi}\pfintr@dot\pfintr@sep}\pfstyl@\ignorespaces}

\def\provedmark#1{\def\prm@rk{#1}}
\def\provedsep#1{\def\prs@p{#1}}

\provedmark{$\square$}
\provedsep{\kern1.25ex}

\def\provedtexttrue{\def\prb@x##1{\fbox{\small##1}}}
\def\provedtextfalse{\def\prb@x##1{\prm@rk}}
\def\provedmarkrighttrue{\let\prhf@l\hfill}
\def\provedmarkrightfalse{\let\prhf@l\relax}

\provedtextfalse
\provedmarkrightfalse

\def\provedboxcontents#1{\expandafter\def\expandafter\@tempa\expandafter{#1}%
  \ifx\@tempa\empty\prm@rk\else\prb@x{#1}\fi}

\def\proved{\ifmmode\eqno{\box\provedbox}\else\hproved\fi}

\def\hproved{\unskip\nobreak\prhf@l\penalty50\prs@p\hbox{}\nobreak\prhf@l
  \box\provedbox{\finalhyphendemerits=0\par}}

\def\captionstyle#1{\def\c@ptstyl@{#1}}
\def\captionintrostyle#1{\def\c@pintr@style{#1}}
\def\captionintrodot#1{\def\c@pintr@dot{#1}}
\def\captionintrosep#1{\def\c@pintr@sep{#1}}

\captionstyle{\small\sf}
\captionintrostyle{\bf}
\captionintrodot{.}
\captionintrosep{\hskip1.25ex}

\long\def\@makecaption#1#2{%
    \vskip\captionskip
    \setbox\@tempboxa\hbox{%
      \ifproofing\@ifundefined{the@label}{}
        {\hbox to 0pt{\vbox to 0pt{\vss\hbox{\tiny\the@label}\bigskip}\hss}}\fi
      \c@ptstyl@{\c@pintr@style #1\c@pintr@dot}\ignorespaces #2}%
    \@captionwidth=\hsize \advance\@captionwidth-2\@captionmargin
    \ifdim \wd\@tempboxa >\@captionwidth {%
        \rightskip=\@captionmargin\leftskip=\@captionmargin
        \unhbox\@tempboxa\par}%
      \else
        \hbox to\hsize{\hfil\box\@tempboxa\hfil}%
    \fi}

\def\end@Float#1{%
  \expandafter\caption\expandafter[\the@title]{%
   {\c@pintr@style%
   \ifx\the@caption\empty\ifx\the@title\empty
   \else\c@pintr@sep\fi\else\c@pintr@sep\fi
    \the@title\ifx\the@caption\empty%
     \expandafter\label\expandafter*\expandafter{\the@label}%
    \else\ifx\the@title\empty%
     \expandafter\label\expandafter*\expandafter{\the@label}%
    \else\c@pintr@dot\c@pintr@sep%
     \expandafter\label\expandafter*\expandafter{\the@label}\fi\fi}%
   \ignorespaces\the@caption}%
  \end{#1}}

\@definecounter{bibenumi}

\def\thebibliography#1{%
 \section*{\refname}\vskip-\lastskip%
 \list{[\arabic{bibenumi}]}{\topsep0pt\settowidth\labelwidth{[#1]}%
 \leftmargin\labelwidth\advance\leftmargin\labelsep\usecounter{bibenumi}}%
 \def\newblock{\hskip .11em plus .33em minus .07em}%
 \sloppy\clubpenalty4000\widowpenalty4000\sfcode`\.=1000\relax}

\makeatother

\proofingfalse
\autolabelfalse

\showfiguresfalse

\newtheorem{stat}{\statname}  \unnumbered{stat}
\newenvironment{statement}[1]{\def\statname{#1}\begin{stat}}{\end{stat}}

\newtheorem{nstat}{\nstatname}[section]

\newtheorem{lemma}[nstat]{Lemma}
\newtheorem{proposition}[nstat]{Proposition}
\newtheorem{theorem}[nstat]{Theorem}

\newtheorem[{\ns}{}]{remark}[nstat]{Remark}

\papersize{215.9truemm}{279.4truemm}
\margins{3.295cm}{2.12cm}

\headline{\hfill}
\footline{\small\hfill--\kern1ex\thepage\kern1ex--\hfill}

\def\ns{\normalshape}

\showfigurestrue

\sectionbeforeskip{1.5\bigskipamount}
\sectionstyle{\centering\normalsize\bf}
\sectionafterskip{\bigskipamount}
\sectionafterindenttrue

\subsectionbeforeskip{\bigskipamount}
\subsectionstyle{\normalsize\bf}
\subsectionafterskip{.5\bigskipamount}
\subsectionafterindenttrue

\paragraphbeforeskip{\bigskipamount}
\paragraphstyle{\centering\normalsize\sl}
\paragraphafterskip{.5\bigskipamount}
\paragraphafternewlinetrue
\paragraphafterindenttrue
\paragraphdot{}

\statementintrostyle{\sc}
\renewtheorem[{\ns}{}]{remark}[nstat]{Remark}

\captionstyle{\small}
\captionintrostyle{\sc}

\newcommand{\Cl}{\mathop{\mathrm{Cl}}\nolimits} 
\newcommand{\Int}{\mathop{\mathrm{Int}}\nolimits} 
\newcommand{\Bd}{\mathop{\mathrm{Bd}}\nolimits}
\newcommand{\Lk}{\mathop{\mathrm{Lk}}\nolimits} 
\newcommand{\Wr}{\mathop{\mathrm{Wr}}\nolimits}

\def\(#1){$(${\sl #1}\/$)$}

\def\bA{\smash{\mkern5mu\overline{\mkern-5mu\phantom{A}}\mkern-13.5mu}A}
\def\bH{\smash{\mkern4.5mu\overline{\mkern-5.5mu\phantom{H}}\mkern-15.425mu}H}
\def\bL{\smash{\mkern3mu\overline{\mkern-4mu\phantom{L}}\mkern-11.25mu}L}

\def\bm{\smash{\mkern2.5mu\overline{\mkern-3.5mu\phantom{m}}\mkern-14.804mu}m}
\def\bn{\smash{\mkern2.5mu\overline{\mkern-3.5mu\phantom{n}}\mkern-9.804mu}n}

\def\be{\smash{\mkern2mu\overline{\mkern-2mu\phantom{\eta}}\mkern-9.583mu}\eta}

\font\ftt cmtt10 at 11pt
\font\fsc cmcsc10 at 12pt
\font\fsl cmsl12 at 12pt

\def\mypagebreak{{\parfillskip0pt\relax\par}\break\noindent}

\begin{document}

\title{\large\bf COVERING MOVES AND KIRBY CALCULUS%
       \label{Version 0.53 / \today}}
\author{\fsc I. Bobtcheva\\
\fsl Dipartimento di Scienze Matematiche\\[-3pt]
\fsl Universit\`a Politecnica delle Marche -- Italia\\
\ftt bobtchev@dipmat.unian.it
\and 
\fsc R. Piergallini\\
\fsl Dipartimento di Matematica e Informatica\\[-3pt]
\fsl Universit\`a di Camerino -- Italia\\
\ftt riccardo.piergallini@unicam.it}
\date{}

\maketitle

\begin{abstract}
\baselineskip13.5pt
\smallskip

\noindent
We show that simple coverings of $B^4$ branched over ribbon surfaces up to certain
local ribbon moves bijectively represent orientable 4-dimensional 2-handlebodies up
to handle sliding and addition/deletion of  cancelling handles. As a consequence, we
obtain an equivalence theorem for simple coverings of $S^3$ branched over links, in
terms of local moves. This result generalizes  to coverings of any degree results by
the second author and Apostolakis, concerning respectively the case of degree $3$ and
$4$. We also provide an extension of our equivalence theorem to possibly non-simple
coverings of $S^3$ branched over embedded graphs.

This work represents the first part of our study of 4-dimensional 2-handle\-bodies. In
the second part \cite{BP05}, we factor such bijective correspondence between simple
coverings of $B^4$ branched over ribbon surfaces and orientable 4-dimensional
2-handlebodies through a map onto the closed morphisms in a universal braided category
freely generated by a Hopf algebra object.



\end{abstract}

\renewtheorem{theorem}{Theorem}

\section*{Introduction}

In the early 70's Hilden \cite{Hi74,Hi76}, Hirsch \cite{Hs74} and Montesinos
\cite{Mo74,Mo76} independently proved that every closed connected oriented
3-manifold can be represented as a 3-fold simple covering of $S^3$ branched over a
link. Successively, Montesinos \cite{Mo78} obtained an analogous representation of
any connected oriented 4-manifold admitting a finite handlebody decomposition with
handles of indices $\leq 2$ as a simple 3-fold covering of $B^4$ branched over a
possibly non-orientable ribbon surface. Actually, the branching surface can always
be made orientable as we remark at the end of Section \ref{todiagram/sec} (cf.
\cite{LP01,PZ03} for other constructions giving directly orientable
ribbon surfaces).

\break

The problem of finding moves relating any two such covering representations of the
same manifold was first considered by Montesinos.
For the 3-dimensional case, in \cite{Mo85} he proposed the two local moves $M_1$
and $M_2$ of Figure \ref{movesM/fig}, where $i$, $j$, $k$ and $l$ are all distinct,
in terms of branching links and monodromy. Here, as well as in all the following
pictures of moves, we draw only the part of the labelled branching set inside the
relevant cell, assuming it to be fixed outside this cell.

\begin{Figure}[htb]{movesM/fig}{}{}
\centerline{\fig{}{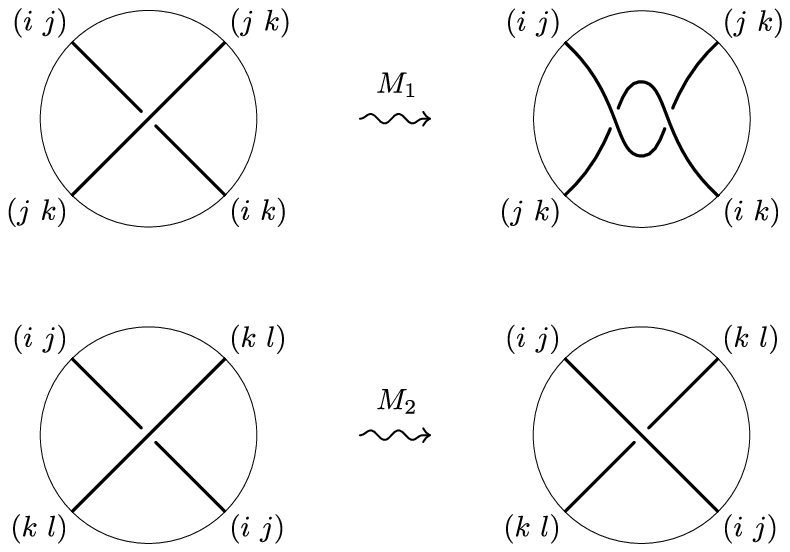}}
\end{Figure}

\begin{Figure}[b]{movesR/fig}{}{}
\centerline{\fig{}{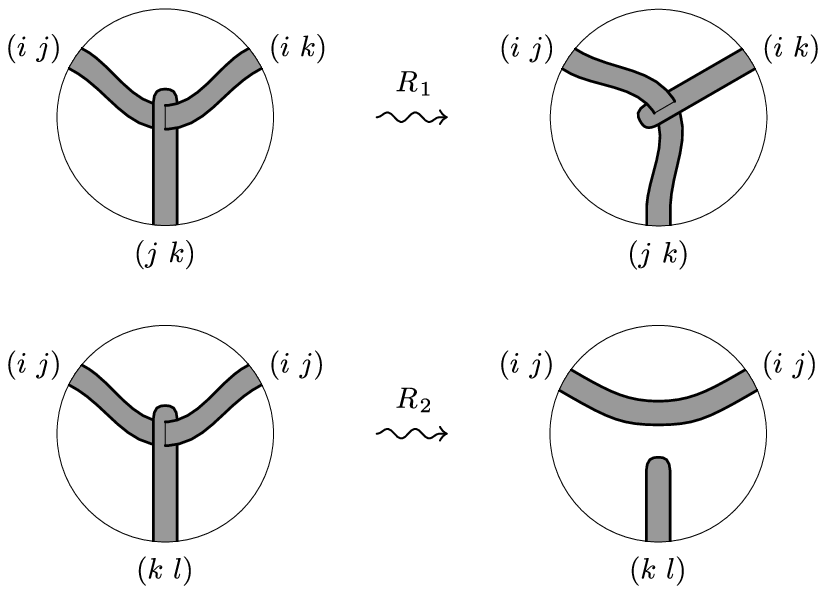}}
\end{Figure}

It is worth observing that the inverse move $M_1^{-1}$ can be realized, up to
labelled isotopy, by a composition of two moves $M_1$. We leave this easy exercise
to the reader, referring to Figure 11 of \cite{Pi95} for the solution. On the other
hand the inverse move $M_2^{-1}$ coincides up to isotopy with the move $M_2$,
becoming distinct from it only after an orientation is fixed on the branching link.

A complete set of moves for 3-fold simple coverings of $S^3$ branched over a link
was given in \cite{Pi91} by the second author. Such moves are non-local, but in
\cite{Pi95} the local moves $M_1$ and $M_2$ are shown to suffice after stabilization
with a fourth trivial sheet.  In \cite{Pi95} the question was also posed, whether
these local moves together with stabilization suffice for covering representations
of arbitrary degree. Recently, Apostolakis \cite{Ap03} answered this question
positively for coverings of degree $4$. 

In this paper, we derive the solution of the moves problem for arbitrary degree 
simple coverings of $S^3$ branched over links (cf. Theorem \ref{equiv3/thm}),
from an equivalence theorem for simple coverings of $B^4$ branched over ribbon
surfaces, that relates the local ribbon moves $R_1$ and $R_2$ of Figure
\ref{movesR/fig}, where $i$, $j$, $k$ and $l$ are all distinct, with the
4-dimensional Kirby calculus (cf. Theorem \ref{equiv4/thm}). In particular, our
result does not depend on the partial ones of \cite{Pi91}, \cite{Pi95} and
\cite{Ap03}.

Analogously to the Montesinos moves, also these ribbon moves generate their inverses
up to labelled isotopy (cf. Proposition \ref{movesInv/thm}). This is obvious for the
move $R_1$, if we think of it as rotation of $120^\circ$ (followed by relabelling),
being $R_1^{-1} = R_1^2$. We leave to the reader to verify that $R_2^{-1}$
coincide with $R_2$ up to labelled isotopy (they become distinct once the branching
ribbon surface is oriented).

Given a connected simple covering $p:M \to B^4$ branched over a ribbon surface $F
\subset B^4$, we have that any 2-dimensional 1-handlebody structure on $F$ induces a
4-dimensional 2-handlebody structure on $M$ (see Section \ref{prelim/sec} for the
definition of $m$-dimensional $n$-handlebody). In fact, the simple covering of $B^4$
branched over the disjoint union of trivial disks $F_0$, representing the 0-handles
of $F$, can be easily seen to be a 4-dimensional 1-handlebody $M_1$. Moreover,
following \cite{Mo78} (cf. also \cite{IP02}), any 1-handle of $F$ attached to $F_0$
corresponds to a 2-handle of $M$ attached to $M_1$. 

In Section \ref{todiagram/sec} we show that handle sliding and handle cancellation in
$F$ give raise to analogous modifications in $M$. Therefore, the 2-handlebody
structure of $M$ turns out to be uniquely determined by the labelled ribbon surface
$F$ up to 2-equivalence, that is up to handle sliding and addition/deletion of
cancelling pairs of handles of indices $\leq 2$ (cf. Section \ref{prelim/sec}).
In other words, any simple covering of $B^4$ branched over a ribbon surface
represents a well defined 2-equivalence class of 4-dimensional 2-handlebodies.

The main result of Montesinos \cite{Mo78} is that any connected oriented
4-dimensional 2-handle\-body $M$ has a 3-fold branched covering representation as
above. The corresponding labelled ribbon surface $F$, with the right 2-dimensional
1-handlebody structure, is obtained from a Kirby diagram of $M$, after it has been
suitably symmetrized with respect to a standard 3-fold simple covering
representation of $M_1$.

In Section \ref{toribbon/sec} (see also Remark \ref{trivialstate/rem}) we give a
different construction of the labelled ribbon surface $F$, similar to that one of
labelled links given in \cite{Mo80} for 3-manifolds (cf. Remark
\ref{smalldisks/rem}). Our construction is simpler and more effective than the
Montesinos one, is canonical up to ribbon moves and better preserves the structure
of the starting Kirby diagram, allowing us to interpret the Kirby calculus in terms
of ribbon moves.

At this point, we are ready to state our first theorem. In substance, it asserts that
simple coverings of $B^4$ branched over ribbon surfaces up certain local isotopy
moves, stabilization and ribbon moves $R_1$ and $R_2$ bijectively represent
4-dimensional 2-handlebodies up to 2-equivalence. For the sake of simplicity, we
consider only the connected case. Nevertheless, as we remark at the end of this
introduction, the statement essentially holds in the general case too, provided the
lower bound for the stabilization degree is replaced by the appropriate one (cf.
Proposition \ref{bijectivity/thm}).

\begin{theorem}\label{equiv4/thm}
Two connected simple coverings of $B^4$ branched over ribbon surfaces represent
2-equivalent 4-dimensional 2-handlebodies if and only if after stabilization to the
same degree $\geq 4$ their labelled branching surfaces can be related by labelled
1-isotopy and a finite sequence of moves $R_1$ and $R_2$.
\end{theorem}

The definition of 1-isotopy is given in Section \ref{prelim/sec}. Here, we limit
ourselves to say that it is essentially generated by the local isotopy moves
shown in Figure \ref{isotopy/fig} (cf. Proposition \ref{1-isotopy/thm}). We do not
know whether 1-isotopy coincides with isotopy of ribbon surfaces (see discussion in
Sections \ref{prelim/sec} and \ref{remarks/sec}). Anyway, we have 1-isotopy instead
of isotopy in the statement of Theorem \ref{equiv4/thm}, due to Lemma
\ref{1-isotopy/2-equiv/thm}. The proof of the theorem\break is achieved in Section
\ref{equivalence/sec}, as a consequence of the above mentioned covering
representation of Kirby calculus. Other main ingredients are Propositions
\ref{to3fold/thm} and \ref{tospecial/thm}.

\begin{Figure}[htb]{isotopy/fig}{}{}
\vskip4pt\centerline{\fig{}{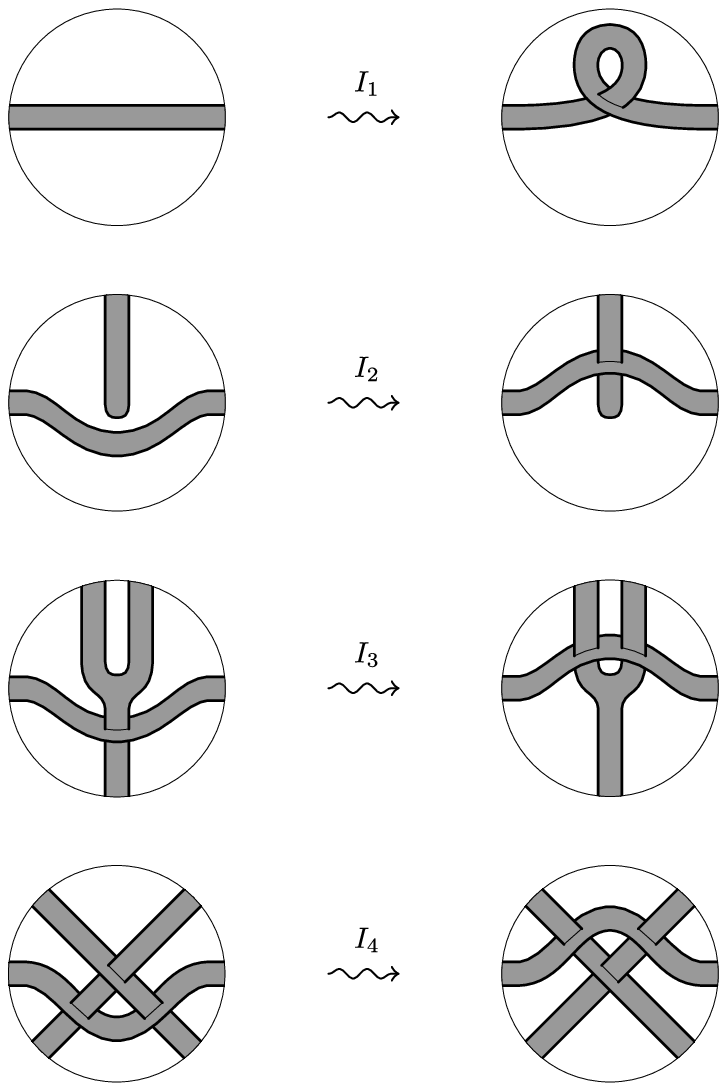}}
\end{Figure}

Now, in order to deal with 3-manifolds, we need to introduce the further moves
depicted in Figures \ref{movesP/fig} and \ref{moveT/fig}. As we see in Section
\ref{equivalence/sec}, these moves allow us to realize respectively
positive/negative blow up and handle trading. In particular, they are not covering
moves in the sense defined in Section \ref{prelim/sec}, since they change the
covering 4-manifold. On the other hand, they do not change the restriction of the
covering over $S^3$, leaving the boundary of the branching surface fixed up to
isotopy.

The next theorem, whose proof is given in Section \ref{equivalence/sec}, tells us
that these last moves together with their inverses and the previous ribbon moves
suffice to completely represent the Kirby calculus for 3-manifolds. Notice that
here, differently from the statement of Theorem \ref{equiv4/thm}, labelled isotopy
can be equivalently used instead of labelled 1-isotopy, since it preserves the
covering manifold up to diffeomorphism. 

\begin{theorem}\label{equiv4b/thm}
Two connected simple coverings of $B^4$ branched over ribbon surfaces represent
4-manifolds with diffeomorphic oriented boundaries if and only if after
stabilization to the same degree $\geq 4$ their labelled branching surfaces can be
related by labelled isotopy and a finite sequence of moves $R_1$, $R_2$,
$P_\pm^{\pm1}$ and
$T^{\pm1}$.
\end{theorem}

\begin{Figure}[htb]{movesP/fig}{}{}
\centerline{\fig{}{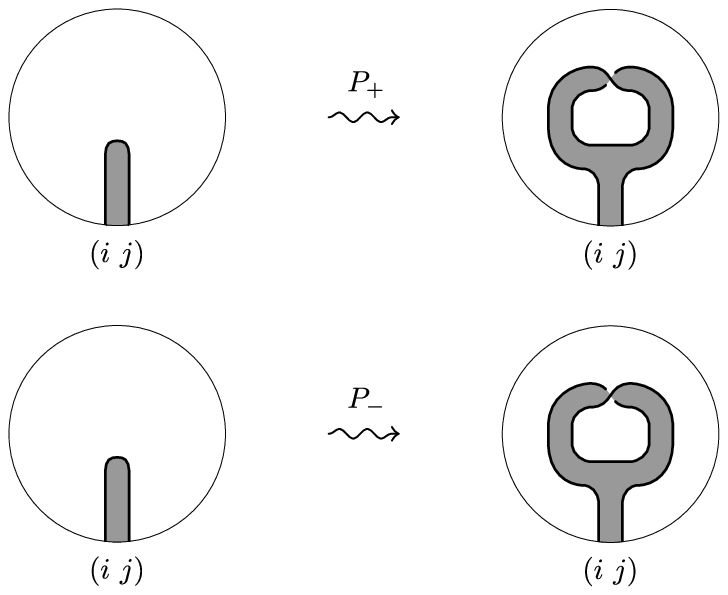}}
\end{Figure}

\begin{Figure}[htb]{moveT/fig}{}{}
\vskip-2pt\centerline{\fig{}{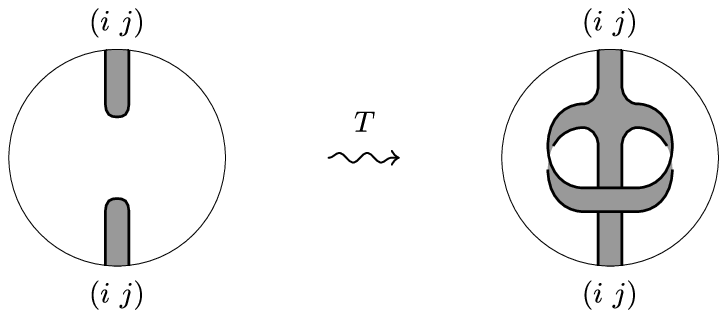}}
\end{Figure}

\medskip

By focusing on the boundary, we observe that the restrictions of the ribbon moves
$R_1$ and $R_2$ to $S^3$ can be realized respectively by Montesinos moves $M_1$ and
$M_2$. This is shown in Figure \ref{cusps/fig} for move $R_1$, while it is trivial
for move $R_2$. In both cases we can apply two Montesinos moves inverse to each
other (with respect to any local orientation of the link as boundary of the surface,
for move $M_2$).

\begin{Figure}[htb]{cusps/fig}{}{}
\vskip6pt\centerline{\fig{}{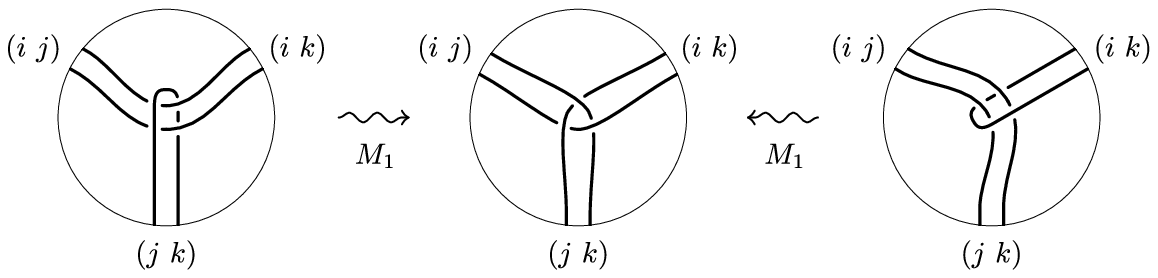}}
\end{Figure}

This observation allows us to derive from Theorem \ref{equiv4b/thm} the following
theorem for simple covering of $S^3$ branched over links. 

\begin{theorem}\label{equiv3/thm}
Two connected simple coverings of $S^3$ branched over links represent diffeomorphic
oriented 3-manifolds if and only if after stabilization to the same degree $\geq 4$
their labelled branching links can be related by labelled isotopy and a finite
sequence of moves $M_1$ and $M_2$.
\end{theorem}

We prove Theorem \ref{equiv3/thm} in Section \ref{equivalence/sec}, as a consequence
of Proposition \ref{link/ribbon/thm}. This says that any labelled link representing
a simple branched covering of $S^3$ can be transformed through Montesinos moves into
the boundary of a labelled ribbon surface representing a simple branched covering of
$B^4$.

Finally, we want to extend Theorem \ref{equiv3/thm} to arbitrary branched coverings
of $S^3$. To do that, we introduce the moves $S_1$ and $S_2$ depicted
in Figure \ref{movesS/fig}. Here, the branching set is allowed to be singular and the
monodromy is not necessarily simple. In fact, $\sigma_1$ and $\sigma_2$ are any
permutations, coherent in the sense defined Section \ref{prelim/sec}, and $\sigma =
\sigma_1 \sigma_2$. This is the reason why we need to specify orientations for the
arcs or equivalently positive meridians to which refer the monodromies.

\begin{Figure}[htb]{movesS/fig}{}{}
\centerline{\fig{}{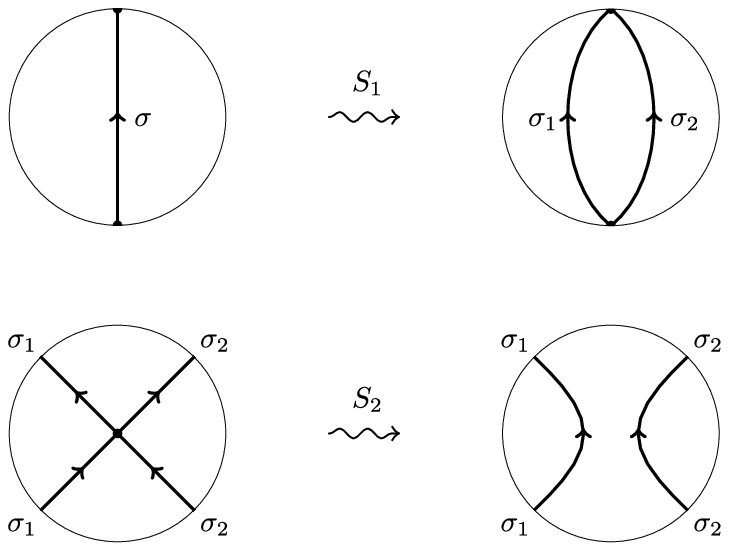}}
\end{Figure}

Our last theorem is the wanted extension of Theorem \ref{equiv3/thm}. Its proof,
given in Section \ref{equivalence/sec}, is based on the fact that moves $S_1$ and
$S_2$ suffice to make simple any branched covering of $S^3$ and to remove all the
singularities from its branching set.

\begin{theorem}\label{equiv3g/thm}
Two connected coverings of $S^3$ branched over a graph represent diffeomorphic
oriented 3-manifolds if and only if after stabilization to the same degree $\geq 4$
their branching graphs can be related by labelled isotopy and a finite sequence of
moves $M_1$, $M_2$, $S_1^{\pm1}$, $S_2^{\pm1}$.
\end{theorem}

We notice that all the above theorems could be easily reformulated to deal with non
connected branched coverings too. Since everything can be done componentwise,
possibly after labelling conjugation, it obviously suffice to stabilize the
coverings to have the same number of sheets $\geq 4$ for corresponding components.
Moreover, as will be clear at the end of Section \ref{equivalence/sec}, the total
degree can be lowered to $3c+1$, where $c$ is the maximum number of components of the
two coverings, if we allow stabilization/destabilization at intermediate stages (cf.
Proposition \ref{bijectivity/thm}).

In conclusion, it is also worth remarking that our results, beyond establishing a
strong relation between branched covering presentations and Kirby diagrams of 3- and
4-manifolds, also provide an effective way to pass from one to the other. We discuss
this aspect in Section \ref{remarks/sec}.


\renewtheorem{theorem}[nstat]{Theorem}

\section{Preliminaries\label{prelim/sec}}

Before going into details, we fix some general notations and conventions about
handlebodies, that will be used in various contexts in the following. We refer to
\cite{GS99} or \cite{Ki89} for all the definitions and basic results not explicitly
mentioned here.

We recall that an {\sl $i$-handle} of dimension $m$ is a copy $H^i$ of $B^i \times
B^{m-i}$ attached to the boundary of an $m$-manifold $M$ by an embedding $\phi:
S^{i-1} \times B^{m-i} \to \Bd M$. The two balls $B^i \times \{0\}$ and $\{0\}
\times B^{m-i}$ in $M' = M \cup_\phi H^i$ are called respectively the {\sl core} and
the {\sl cocore} of $H^i$, while their boundaries $S^{i-1} \times \{0\}$ and $\{0\}
\times S^{{m-i-1}}$ are called the {\sl attaching sphere} and the {\sl belt sphere}
of $H^i$. Inside $H^i$, {\sl longitudinal} means parallel to the core and {\sl
transversal} means parallel to the cocore. Up to isotopy, the attaching map $\phi$
is completely determined by the attaching sphere together with its {\sl framing} in
$\Bd M$, given by $S^{i-1} \times \{\ast\}$ for any $\ast \in B^{m-i} - \{0\}$.

Then, an {\sl $n$-handlebody} of dimension $m$ is defined by induction on $n$ to be
obtained by simultaneously smoothly attaching a finite number of $n$-handles to an
$(n-1)$-handlebody of the same dimension $m$, starting with a disjoint union of
0-handles for $n = 0$.

By a well known result of Cerf \cite{Ce70} (cf. \cite{GS99} or \cite{Ki89}), two
handlebodies of the same dimension are diffeomorphic (forgetting their handle
structure), if and only if they can be related by a finite sequence of the following
modifications: 1) isotoping the attaching map of $i$-handles; 2) adding/deleting a
pair of {\sl cancelling handles}, that is a $i$-handle $H^i$ and a $(i+1)$-handle
$H^{i+1}$, such that the attaching sphere of $H^{i+1}$ intersects the belt sphere of
$H^i$ transversally in a single point; 3) {\sl handle sliding} of one $i$-handle
$H^i_1$ over another one $H^i_2$, that means pushing the attaching sphere of $H^i_1$
through the belt sphere of $H^i_2$.

We call {\sl $k$-deformation} any finite sequence of the above modifications such
that at each stage we have an $n$-handlebody with $n \leq k$, that is we start from
a $n$-handlebody with $n \leq k$ and never add any cancelling $i$-handle with $i >
k$. Furthermore, we call {\sl $k$-equivalent} two handlebodies related by a
$k$-deformation.

In particular, any compact surface with non-empty boundary has a 1-handlebody
structure and any two such structures are easily seen to be 1-equivalent (cf.
proof of Proposition \ref{1-handles/thm}).

The other relevant case for our work is that one of orientable 4-manifolds (with
non-empty boundary) admitting a 4-dimensional 2-handlebody structure.
Any two such structures are 3-equivalent, but whether they are 2-equivalent is a
much more subtle open question, which is expected to have negative answer 
(cf. Section I.6 of \cite{Ki89} and Section 5.1 of \cite{GS99}). 
This question seems to be strongly related to the problem of finding isotopy
moves for ribbon surfaces in $B^4$. In fact, as we will see, 4-dimensional
2/3-deformations correspond by means of branched coverings to regularly embedded
2-dimensional 1/2-deformations of branching surfaces in $B^4$ (cf. Proposition
\ref{1-equiv/2-equiv/thm} and the discussion in Section \ref{remarks/sec}).

\paragraph{Links}

As usual, we represent a link $L \subset R^3 \subset R^3 \cup \infty \cong
S^3$ by a planar {\sl diagram} $D \subset R^2$, consisting of the orthogonal
projection of $L$ into $R^2$, that can be assumed self-transversal after a suitable
horizontal (height preserving) isotopy of $L$, with a {\sl crossing state} for each
double point, telling which arc passes over the other one.
\mypagebreak
Such a diagram $D$
uniquely determines $L$ up to vertical isotopy. On the other hand, link isotopy can
be represented in terms of diagrams by crossing preserving isotopy in $R^2$ and
Reidemeister moves.

A link $L$ is called {\sl trivial} if it bounds a disjoint union of disks in $R^3$.
It is well known that any link diagram $D$ can be transformed into a diagram $D'$ of
a trivial link by suitable crossing changes, that is by inverting the state of some
of its crossings. We say that $D'$ is a {\sl trivial state} of $D$. Actually, any
link diagram $D$ has many trivial states, but it is not clear at all how they are
related to each other. For this reason, we are lead to introduce the more
restrictive notions of vertically trivial link and vertically trivial state of a
link diagram.

We say that a link $L$ is {\sl vertically trivial} if it meets any horizontal plane
(parallel to $R^2$) in at most two points belonging to the same component. In this
case, the height function separates the components of $L$ (that is the height
intervals of different components are disjoint), so that we can vertically order the
components of $L$ according to their height. Moreover, each component can be
split into two arcs on which the height function is monotone, assuming the only
unique minimum and maximum values at the common endpoints. Then, all the (possibly
degenerate) horizontal segments spanned by $L$ in $R^3$ form a disjoint union of
disks bounded by $L$. This proves that $L$ is a trivial link.

By a {\sl vertically trivial state} of a link diagram $D$ we mean any trivial state
of $D$ which is the diagram of a vertically trivial link. A vertically trivial
state $D'$ of $D$ can be constructed by the usual naive unlinking procedure: 
1) number the components of the link $L$ represented by $D$ and fix on each component
an orientation and a starting point away from crossings; 2) order the points of $L$
lexicographically according to the numbering of the components and then to the
starting point and the orientation of each component; 3) resolve each double point
of $D$ into a crossings of $D'$ by letting the arc which comes first in the order
pass under the other one. The link $L'$ represented by $D'$ can be clearly assumed
to be vertically trivial, considering on it a height function which preserves the
order induced by the vertical bijection with $L$ except for a small arc at the end of
each component. Figure \ref{height/fig} \(a) shows how the height function of a
component looks like with respect to a parametrization having the starting point and
the orientation fixed above. Keeping the parametrization fixed but changing the
starting point or the orientation we get different height functions as in Figures
\ref{height/fig} \(b) and \(c) respectively.

\smallskip
\begin{Figure}[htb]{height/fig}{}{}
\centerline{\fig{}{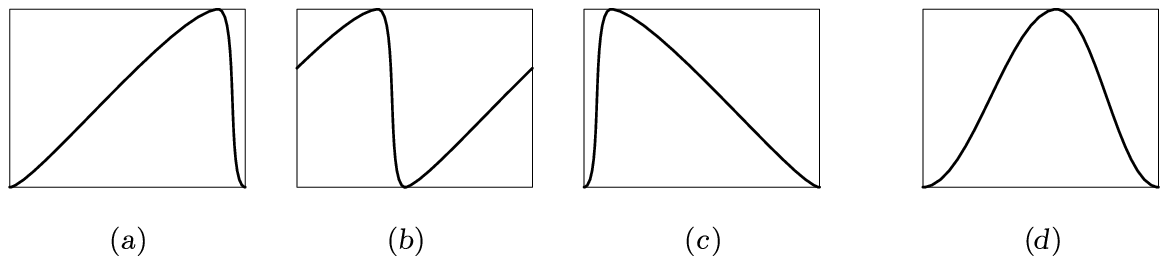}}
\end{Figure}

Notice that the above unlinking procedure gives us only very special vertically
trivial states. While it is clear how to pass from \(a) to \(b), by moving the
starting point along the component, going from \(a) to \(c) turns out to be
quite mysterious without considering generic vertically trivial states.
The height function of a component for such a state, with respect to a
parametrization starting from the unique minimum point, looks like in Figure
\ref{height/fig} \(e), that is apparently an intermediate state between \(a)
and \(c). The following proposition settles the problem of relating different
vertically trivial states of the same link diagram.

\begin{proposition} \label{vertstat/thm}
Any two vertically trivial states $D'$ and $D''$ of a link diagram $D$ are related by
a sequence $D_0, D_1, \dots, D_n$ of vertically trivial states of $D$, such that $D_0
= D'$, $D_n = D''$ and, for each $i = 1, \dots, n$, $D_i$ is obtained from $D_{i-1}$
by changing a single self-crossing of one component or by changing all the crossings
between two vertically adjacent components.
\end{proposition}

\begin{proof} 
Since the effect of changing all the crossings between two vertically adjacent
components is the transposition of these components in the vertical order, by
iterating this kind of modification we can permute as we want the vertical order of
all the components. Hence, we only need to address the case of a knot diagram.

Given a knot diagram $D \subset R^2$ with double points $x_1, \dots, x_n \in R^2$,
we consider a parametrization $h:S^1 \to D$ and denote by $t'_i,t''_i \in S^1$ the
two values of the parameter such that $h(t'_i) = h(t''_i) = x_i$, for any $i = 1,
\dots, n$.

For any smooth knot $K \subset R^3$ which projects to a vertically trivial state of
$D$, let $h_K: S^1 \to R^3$ be the parametrization of $K$ obtained by lifting $h$
and $f_K:S^1 \to R$ be the composition of $h_K$ with the height function. Then, $f_K$
is a smooth function with the following properties: 1)~$f_K$ has only one minimum
and one maximum; 2)~$f_K(t'_i) \neq f_K(t''_i)$, for any $i = 1, \dots, n$.
In this way, the space of all smooth knots which project to vertically trivial
states of $D$ can be identified with the space of all smooth functions $f:S^1 \to R$
satisfying properties 1 and 2.

Now, the space $\cal S$ of all smooth functions $f:S^1 \to R$ satisfying property 1
is clearly pathwise connected, while the complement $\cal C \subset \cal S$ of
property 2 is a closed codimension 1 stratified subspace. Therefore, if $K'$ and
$K''$ are knots projecting to the vertically trivial states $D'$ and $D''$, then we
can join $f_{K'}$ and $f_{K''}$ by a path in $\cal S$ transversal with respect to
$\cal C$. This, path gives rise to a finite sequence of self-crossing changes as in
the statement, one for each transversal intersection with $\cal C$.
\end{proof}

We remark that the singular link between two consecutive vertically trivial\break
states, obtained from each other by a single self-crossing change, is trivial.
Namely, the unique singular component spans a 1-point union of two disks, disjoint
from all the other components. This fact, which  will play a crucial role in the
proof of Proposition \ref{diag/ribbon/thm}, follows from \cite{Sc93} but can also
be easily proved directly by inspection.

\paragraph{Ribbon surfaces}

A smooth compact surface $F \subset B^4$ with $\Bd F \subset S^3$ is
called a {\sl ribbon surface} if the Euclidean norm restricts to a Morse function on
$F$ with no local maxima in $\Int F$. Assuming $F \subset R^4_- \subset R^4_- \cup
\{\infty\} \cong B^4$, this property is topologically\break equivalent to the fact
that the fourth Cartesian coordinate restricts to a Morse height function on $F$
with no local maxima in $\Int F$. Such a surface $F \subset R^4_-$ can be
horizontally (preserving the height function) isotoped to make its orthogonal
projection into $R^3$ a self-transversal immersed surface, whose double points form
disjoint arcs as in Figure \ref{ribbon/fig}.

\begin{Figure}[htb]{ribbon/fig}{}{}
\centerline{\fig{}{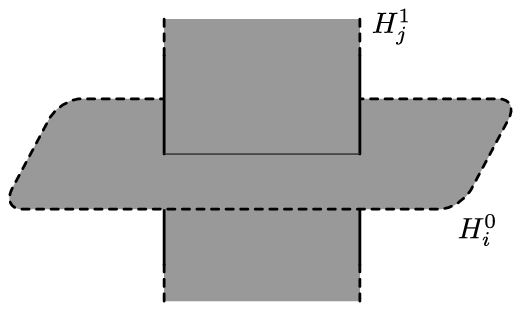}}
\end{Figure}

We will refer to such a projection as a 3-dimensional {\sl diagram} of $F$. Actually,
any immersed compact surface $F \subset R^3$ with no closed components and all
self-intersections of which are as above, is the diagram of a ribbon surface uniquely
determined up to vertical isotopy. This can be obtained by pushing $\Int F$ down
inside $\Int R^4_-$ in such a way that all self-intersections disappear. 

\medskip

In the following, {\sl ribbon surfaces will be always represented by diagrams and
considered up to vertical isotopy}. Moreover, {\sl we will use the same notations for
a ribbon surface and for its diagram in $R^3$}, disregarding the projection. By the
above observation there will be no danger of confusion, provided that the ambient
space will be clear.

\medskip

Since a ribbon surface $F$ has no closed components, it admits a handlebody
decomposition $F = H^0_1 \cup \dots \cup H^0_m \cup H^1_1 \cup \dots \cup H^1_n$
with only 0- and 1-handles. Such a 1-handlebody decomposition is called {\sl
adapted}, if each ribbon self-intersection involves an arc contained in the interior
of a 0-handle and a proper transversal arc inside a 1-handle. (cf. \cite{Ru85}). 

By an {\sl embedded 2-dimensional 1-handlebody} we mean a ribbon surface endowed with
an adapted 1-handlebody decomposition as above. Looking at the diagram, we have that
the $H^0_i$'s are disjoint non-singular disks, while the $H^1_j$'s are non-singular
bands attached to the $H^0_i$'s and possibly passing across them as shown in Figure
\ref{ribbon/fig}. Moreover, we can think of $F$ as a smooth perturbation of the
boundary of\break $((H^0_1 \cup \dots \cup H^0_m) \times [0,-1]) \cup ((H^1_1 \cup
\dots \cup H^1_n) \times [0,-1/2])$, in such a way that the handlebody decomposition
turns out to be induced by the height function. 

\medskip

We say that two embedded 2-dimensional 1-handlebodies are equivalent up to
{\sl embedded 1-deformation}, or briefly that they are {\sl 1-equivalent}, if they
are related by a finite sequence of the following modifications:
\begin{itemize}\itemsep3pt

\item[\(a)]\vskip-\lastskip\smallskip
{\sl adapted isotopy}, that is isotopy of 1-handlebodies in $R^4$, all adapted except
for a finite number of intermediate critical stages, at which one of the
modifications described in Figure \ref{ribbon1/fig} takes place (between any two
such critical stages, we have isotopy of diagrams in $R^3$, preserving ribbon
intersections);
\begin{Figure}[htb]{ribbon1/fig}{}{}
\centerline{\fig{}{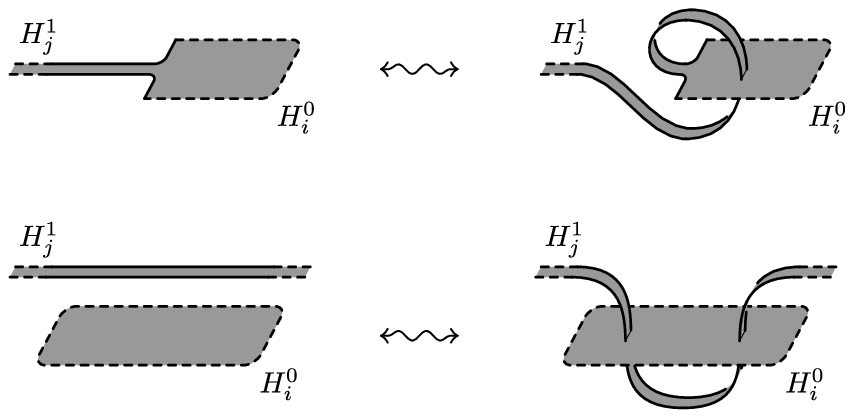}}
\end{Figure}

\item[\(b)] 
{\sl ribbon intersection sliding}, allowing a ribbon intersection to run along a
1-handle from one 0-handle to another one, as shown in Figure \ref{ribbon2/fig};%
\begin{Figure}[htb]{ribbon2/fig}{}{}
\centerline{\fig{}{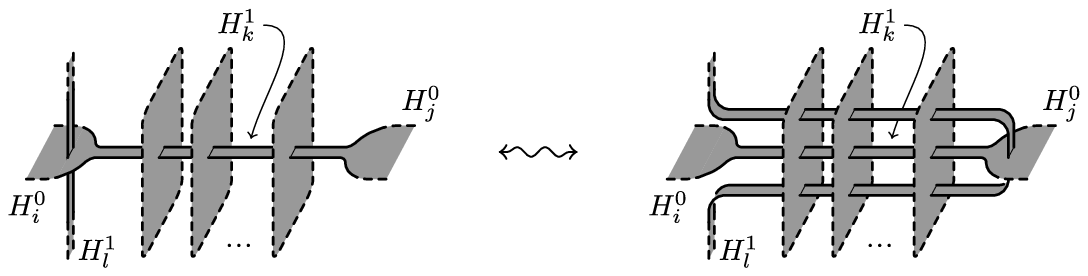}}
\end{Figure}

\item[\(c)] 
{\sl embedded 0/1-handles operations}, that is addition/delection of cancelling
pairs of 0/1-handles and embedded 1-handle slidings (see Figure \ref{ribbon3/fig}).
\begin{Figure}[htb]{ribbon3/fig}{}{}
\centerline{\fig{}{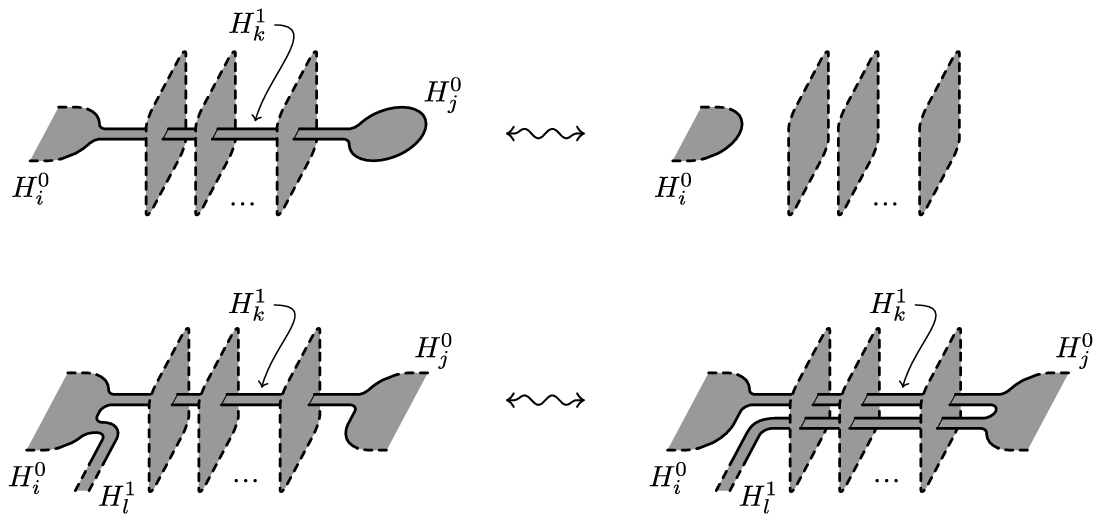}}
\end{Figure}

\end{itemize}\vskip-\lastskip\smallskip

\medskip

We observe that the second modification of Figure \ref{ribbon1/fig} is actually
redundant in presence of the handle operations of Figure \ref{ribbon3/fig} (cf.
proof of Proposition \ref{1-isotopy/thm}).\break
It is also worth noticing that no twist appears in the 1-handle $H^1_k$ of Figures
\ref{ribbon2/fig} and \ref{ribbon3/fig}, since $H^0_i$ and $H^0_j$ can be assumed to
be distinct in all the cases, up to addition/delection of cancelling pairs of
0/1-handles, where they are always distinct.

\begin{proposition} \label{1-handles/thm}
All the adapted 1-handlebody decompositions of a given ribbon surface are
1-equivalent as embedded 2-dimensional 1-handlebodies. More precisely, they are
related to each other by the special cases without vertical disks of the moves of
Figures \ref{ribbon2/fig} and \ref{ribbon3/fig}, realized (up to isotopy of diagrams)
in such a way that the surface is kept fixed.
\end{proposition}

\begin{proof}
First of all, we observe that the moves specified in the statement allow us to
realize the following two modifications: 1) split a 0-handle along any regular arc
avoiding ribbon intersections in the diagram, into two 0-handles joined by a new
1-handle;  2) split a 1-handle at any transversal arc avoiding ribbon intersections
in the diagram, into two 1-handles, by inserting a new 0-handle along it. We leave
the straightforward verification of this to the reader.

Let $F = H^0_1 \cup \dots \cup H^0_m \cup H^1_1 \cup \dots \cup H^1_n = \bH^0_1
\cup \dots \cup \bH^0_{\bm} \cup \bH^1_1 \cup \dots \cup \bH^1_{\bn}$\break be any
two 1-handlebody decompositions of a ribbon surface $F$, which we denote
respectively by $H$ and $\bH$. After having suitably split the 1-handles, we can
assume that any 1-handle contains at most one ribbon self-intersection of $F$ and
that this coincides with its cocore. Up to isotopy, we can also assume that the
1-handles of $H$ and $\bH$ whose cocore is the same self-intersection arc coincide.
Let $H_1 = \bH_1, \dots,\break H_k = \bH_k$ be these 1-handles. Then, it suffices to
see how to make the remaining 1-handles $H^1_{k+1}, \dots, H^1_n$ into $\bH^1_{k+1},
\dots, \bH^1_{\bn}$, without changing $H^1_1, \dots, H^1_k$.

Calling $\eta_i$ (resp. $\be_j$) the cocore of $H^1_i$ (resp. $\bH^1_j$), we
have $\eta_1 = \be_1, \dots, \eta_k = \be_k$, while the arcs $\eta_{k+1},
\dots, \eta_n$ can be assumed to be transversal with respect to the arcs $\be_{k+1},
\dots, \be_{\bn}$. Up to isotopy, we can think of each 1-handle as a tiny regular
neighborhood of its cocore, so that the intersection between $H^1_{k+1} \cup \dots
\cup H^1_n$ and $\bH^1_{k+1} \cup \dots \cup \bH^1_{\bn}$ consists only of a certain
number $h$ of small four-sided regions. 

We eliminate all these intersection regions in turn, by pushing them outside $F$
along the $\bH^1_j$'s. This is done by performing on $H$ moves of the types
specified in the statement, as suggested by the following Figure \ref{equiv1/fig},
which concerns the $l$-th elimination. Namely, in \(a) we assume that the
intersection is the first one along $\be_j$ starting from $\Bd F$, then we generate
the new 1-handle $H^1_{n+l}$ by 0-handle splitting to get \(b), finally \(c) is
obtained by handle sliding.

\begin{Figure}[htb]{equiv1/fig}{}{}
\centerline{\fig{}{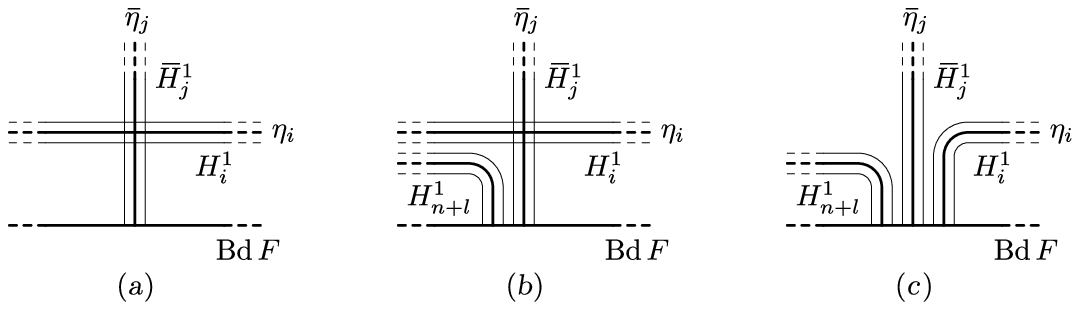}}
\end{Figure}

After that, $H$ has been changed into a new handlebody decomposition $H'$ with
1-handles $H^1_1, \dots, H^1_{n + h}$, such that $H^1_i$ is the same as above for $i
\leq k$, while it is disjoint from the $\bH^1_j$'s for $i > k$. Hence, $H^1_1,
\dots, H^1_k,H^1_{k+1}, \dots, H^1_{n+h}, \bH^1_{k+1}, \dots, \bH^1_{\bn}$ can be
considered as the 1-handles of a handlebody decomposition of $F$ which can be
obtained from both $H'$ and $\bH$ by 0-handle splitting.
\end{proof}

Now, forgetting the 1-handlebody structure, 1-equivalence of embedded 2-dimensional
1-handlebodies induces an equivalence relation between ribbon surfaces, that we call
{\sl 1-isotopy}. More precisely, two ribbon surfaces are 1-isotopic if and only if
they admit 1-equivalent 1-handlebody decompositions. By the above proposition,
this implies that actually all their 1-handlebody decompositions are 1-equivalent.

Of course 1-isotopy implies isotopy, but the converse is not known. In fact, the
problem of finding a complete set of moves representing isotopy of ribbon surfaces
is still open. We will come back to this delicate aspect later. 

As we anticipated in the Introduction, the next proposition says that 1-isotopy
is generated by the local isotopy moves of Figure \ref{isotopy/fig}, up to diagram
isotopy in $R^3$, that means isotopy preserving ribbon intersections.

\begin{proposition} \label{1-isotopy/thm}
Two ribbon surfaces are 1-isotopic if and only if they can be related by a finite
sequence of diagram isotopies and moves $I_1,\dots,I_4$ and their inverses.
\end{proposition}

\begin{proof}
On one hand, we have to realize the modifications of Figures \ref{ribbon1/fig},
\ref{ribbon2/fig} and \ref{ribbon3/fig}, disregarding the handlebody structure, by
moves $I_1, \dots, I_4$ and their inverses. Of course, it is enough to do that in
one direction, say from left to right. Proceeding in the order: one move $I_1$
suffices for the upper part of Figure \ref{ribbon1/fig}, while the lower part can be
obtained by combining one move $I_2$ with one move $I_3$; Figure \ref{ribbon2/fig}
requires three moves for each vertical disk, one $I_2$, one $I_3$ and one $I_4$;
the upper (resp. lower) part of Figure \ref{ribbon3/fig} can be achieved by one
move $I_2$ (resp. $I_3$) for each vertical disk.

On the other hand, the surfaces of Figure \ref{isotopy/fig} can be easily provided
with adapted handlebody decompositions, in such a way that the relations just
described between moves $I_1, \dots, I_4$ and the above modifications can be
reversed. In fact, only the special cases of those modifications with one vertical
disk are needed.
\end{proof}

\paragraph{Branched coverings}

A non-degenerate PL map $p:M \to N$ between compact PL manifolds of the same
dimension $m$ is called a {\sl branched covering} if there exists an
$(m-2)$-dimensional subcomplex $B_p \subset N$, the {\sl branching set} of $p$, such
that the restriction $p_|: M - p^{-1}(B_p) \to N - B_p$ is an ordinary covering of
finite degree $d$. If $B_p$ is minimal with respect to such property, then we have
$B_p = p(S_p)$, where $S_p$ is the {\sl singular set} of $p$, that is the set of
points at which $p$ is not locally injective. In this case, both $B_p$ and $S_p$, as
well as the {\sl pseudo-singular set} $S'_p = \text{Cl}(p^{-1}(B_p)-S_p)$, are
(possibly empty) homogeneously $(m-2)$-dimensional complexes.

Since $p$ is completely determined, up to PL homeomorphisms, by the ordinary covering
$p_|$ (cf. \cite{Fo57}), we can describe it in terms of its branching set $B_p$ and
its {\sl monodromy} $\omega_p:\pi_1(N-B_p,\ast) \to \Sigma_d$, defined up to
conjugation in $\Sigma_d$, depending on the choice of the base point $\ast$ and on
the numbering of $p^{-1}(\ast)$. 
In particular, the monodromies of the meridians around the $(m-2)$-simplices of
$B_p$ determine the structure of the singularities of $p$. If all such monodromies
are transpositions, then we say that $p$ is {\sl simple}. In this case, every point
in the interior of a $(m-2)$-simplex of $B_p$ is the image of one singular point, at
which $p$ is topologically equivalent to the complex map $z \mapsto z^2$, and $d-2$
pseudo-singular points.

Starting from $B_p \subset N$ and $\omega_p$, we can explicitly reconstruct $M$ and
$p$ by following steps: 1) choose a $(m-1)$-dimensional {\sl splitting complex},
that means a subcomplex $C \subset N - \{\ast\}$ such that $B_p \subset C$ and the
restriction $\omega_{p|}:\pi_1(N - C,\ast) \to \Sigma_d$ vanishes; 2) cut $N$ along
$C$ in such a way that each $(m-1)$-simplex $\sigma$ of $C$ gives raise to 2
simplices $\sigma^-$ and $\sigma^+$; 3) take $d$ copies of the obtained complex
(called the {\sl sheets} of the covering) and denote by $\sigma^\pm_1, \dots,
\sigma^\pm_d$ the corresponding copies of $\sigma^\pm$; 4) identify in pairs the
$\sigma^\pm_i$'s according to the monodromy $\rho = \omega_p(\alpha)$ of a loop
$\alpha$ meeting $C$ transversally at one point of $\sigma$, namely identify
$\sigma^-_i$ with $\sigma^+_{\rho(i)}$. Up to PL homeomorphisms, $M$ is the result
of such identification and $p$ is the map induced by the natural projection of the
sheets onto $N$.

\medskip

A convenient representation of $p$ can be given by labelling each $(m-2)$-simplex of
$B_p$ by the monodromy of a preferred meridian around it and each generator (in a
finite generating set) of $\pi_1(N,\ast)$ by its monodromy, since those loops
together
\mypagebreak
generate $\pi_1(N-B_p,\ast)$. Of course, only the labels on $B_p$ are
needed if $N$ is simply connected. In any case, with a slight abuse of language if
$N$ is not simply connected, we refer to such a representation as a {\sl labelled
branching set}.

\medskip

Two branched coverings $p:M \to N$ and $p':M' \to N$ are called {\sl equivalent}
iff there exists PL homeomorphism $h:N \to N$ isotopic to the identity which
lifts to a PL homeomorphism $k:M \to M'$. By the classical theory of ordinary
coverings and \cite{Fo57}, such a lifting $k$ of $h$ exists iff $h(B_p) = B_{p'}$ and
$\omega_{p'}h_* = \omega_p$ up to conjugation in $\Sigma_d$, where $h_\ast:\pi_1(N -
B_p,\ast) \to \pi_1(N - B_{p'},h(\ast))$ is the homomophism induced by $h$.\break
Therefore, in terms of labelled branching set, the equivalence of branched coverings
can be represented by {\sl labelled isotopy}.

\medskip

By a {\sl covering move}, we mean any non-isotopic modification making a labelled
branching set representing a branched covering $p:M \to N$ into one representing a
different branched covering $p':M \to N$ between the same manifolds (up to PL
homeomorphisms). We call such a move {\sl local}, if the modification takes place
inside a cell and can be performed whatever is the rest of labelled branching set
outside. In the figures depicting local moves, we will draw only the portion
of the labelled branching set inside the relevant cell, assuming everything else to
be fixed.

As a primary source of covering moves, we consider the following two very general
equivalence principles (cf. \cite{PZ03}). Several special cases of these principles
have already appeared in the literature and we can think of them as belonging to the
``folklore'' of branched coverings.

\begin{statement}{Disjoint monodromies crossing}
Subcomplexes of the branching set of a covering that are labelled with disjoint
permutations can be isotoped independently from each other without changing the
covering manifold.
\end{statement}

The reason why this principle holds is quite simple. Namely, being the labelling of
the subcomplexes disjoint, the sheets non-trivially involved by them do not
interact, at least over the region where the isotopy takes place. Hence, the
relative position of such subcomplexes is not relevant in determining the covering
manifold. Typical applications of this principle are the local moves $M_2$ and $R_2$
(cf. Figures \ref{movesM/fig} and \ref{movesR/fig}).

It is worth observing that, abandoning transversality, the disjoint monodromies
crossing principle also gives the special case of the next principle when the
$\sigma_i$'s are disjoint and $L$ is empty.

\begin{statement}{Coherent monodromies merging}
Let $p:M \to N$ be any branched covering with branching set $B_p$ and let $\pi: E \to
K$ be a connected disk bundle imbedded in $N$, in such a way that: 1)~there exists
a (possibly empty) subcomplex $L \subset K$ for which $B_p \cap \pi^{-1}(L) = L$ and
the restriction of $\pi$ to $B_p \cap \pi^{-1}(K-L)$ is an unbranched covering of
$K-L$; 2)~the monodromies $\sigma_1, \dots, \sigma_n$ relative to a fundamental
system $\omega_1, \dots, \omega_n$ for the restriction of $p$ over a given disk
$D = \pi^{-1}(x)$, with $x \in K - L$, are coherent in the sense that $p^{-1}(D)$ is
a disjoint union of disks. Then, by contracting the bundle $E$ fiberwise to $K$, we
get a new branched covering $p':M \to N$, whose branching set $B_{p'}$ is equivalent
to $B_p$, except for the replacement of $B_p \cap \pi^{-1}(K-L)$ by $K-L$, with the
labelling uniquely defined by letting the monodromy of the meridian $\omega =
\omega_1 \dots \omega_n$ be $\sigma = \sigma_1 \dots \sigma_n$.
\end{statement}

We remark that, by connectedness and property 1, the coherence condition required in
2 actually holds for any $x \in K$. Then, we can prove that $p$ and $p'$ have the
same covering manifold, by a straightforward fiberwise application of the Alexander's
trick to the components of the bundle $\pi \circ p: p^{-1}(E) \to K$. A coherence
criterion can be immediately derived from Section 1 of \cite{MP01}.

The coherent monodromy merging principle originated from a classical perturbation
argument in algebraic geometry and appeared in the literature as a way to deform
non-simple coverings between surfaces into simple ones, by going in the opposite
direction from $p'$ to $p$ (cf. \cite{BE79}). In the same way, it can be used
in dimension 3, both for achieving simplicity (cf. \cite{Hr03}) and removing
singularities from the branching set. We will do that in the proof of Theorem
\ref{equiv3g/thm} by means of the moves $S_1$ and $S_2$ of Figure \ref{movesS/fig},
which are straightforward applications of this principle. Actually, analogous
resuts could be proved in dimension 4, but we will not do it here.

The coherent monodromy merging principle, also provides an easy way to verify that
$M_1$ and $R_1$ are local covering moves, as shown in Figures
\ref{moveM1pr/fig} and \ref{moveR1pr/fig}. In both these figures, we apply the
principle for going from \(a) to \(b) and from \(c) to \(d), while \(b) and \(c) are
equivalent up to labelled isotopy.

\begin{Figure}[htb]{moveM1pr/fig}{}{}
\vskip4pt\centerline{\fig{}{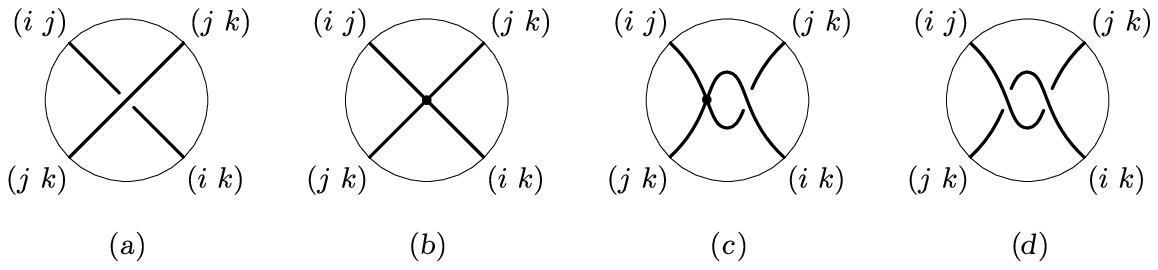}}
\end{Figure}

\begin{Figure}[htb]{moveR1pr/fig}{}{}
\centerline{\fig{}{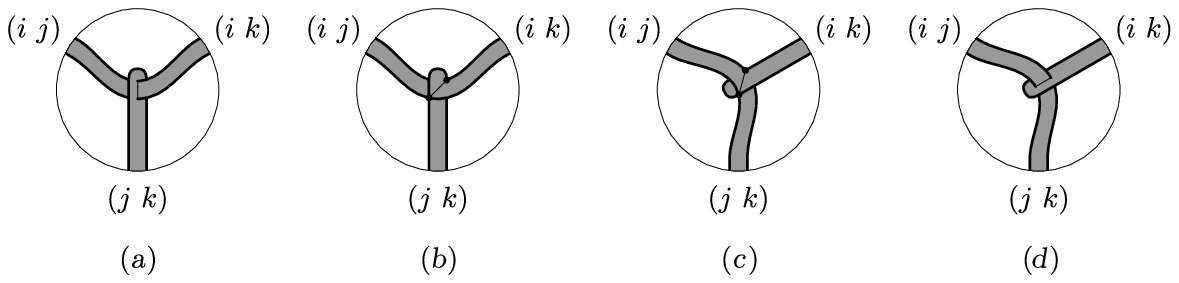}}
\end{Figure}

\medskip

So far we have seen that all the moves presented in the Introduction, except for
moves $T$ and $P_\pm$, are local covering moves. However, we will give a different
proof of that for moves $R_1$ and $R_2$ in Section \ref{todiagram/sec}, by relating
them to 2-deformations of 4-handlebodies.

\medskip

Now, we consider the notion of stabilization that appears in all the equivalence
theorems stated in the Introduction. This is a particular local covering move, which
makes sense only for branched coverings of $S^m$ or $B^m$ and, differently from all
the previous moves, changes the degree of the covering, increasing it by one.

\begin{statement}{Stabilization}
A branched covering $p:M \to S^m$ (resp. $p:M \to B^m$) of degree $d$, can be
stabilized to degree $d + 1$ by adding to the labelled branching set a trivial
separate $(m-2)$-sphere (resp. regularly embedded $(m-2)$-disk) labelled with the
transposition $(i\ d{+}1)$, for some $i = 1, \dots, d$.
\end{statement}

The covering manifold of such a stabilization is still $M$, up to PL homeomorphisms.
In fact, it turns out to be the connected sum (resp. boundary connected sum) of $M$
itself, consisting of the sheets $1, \dots, d$, with the copy of $S^m$ (resp. $B^m$)
given by the extra trivial sheet $d+1$.

By {\sl stabilization to degree $n$} (or {\sl $n$-stabilization}) of a branched
covering $p:M \to S^m$ (resp. $p:M \to B^m$) of degree $d \leq n$ we mean the
branched covering of degree $n$ obtained from it by performing $n - d$
stabilizations as above. In particular, this leaves $p$ unchanged if $d = n$.

\medskip

We conclude this paragraph by focusing on the branched coverings we will deal
with in the following sections, that is coverings of $S^3$ branched over links or
embedded graphs and coverings of $B^4$ branched over ribbon surfaces. We recall
that in this context PL and smooth are interchangeable.

We represent a $d$-fold covering of $p:M \to S^3$ branched over a link $L \subset
S^3$, by a $\Sigma_d$-labelled oriented diagram $D$ of $L$ describing the monodromy
of $p$ in terms of the Wirtinger presentation of $\pi_1(S^3 - L)$ associated to $D$.
Namely, we label each arc of $D$ by the monodromy of the standard positive meridian
around it. Of course, the Wirtinger relations impose constraints on the labelling at
crossings, and each $\Sigma_d$-labelling of $D$ satisfying such constraints do
actually represent a $d$-fold covering of $S^3$ branched over $L$. Then, labelled
isotopy can be realized by means of labelled Reidemeister moves. 

For simple coverings, the orientation of $D$ is clearly unnecessary and there are
three possible ways of labelling the arcs at each crossing: either all with the same
transposition $(i\ j)$ or like at the two crossings in the left side of Figure
\ref{movesM/fig}.

The Montesinos-Hilden-Hirsch representation theorem of closed connected oriented
3-manifolds as branched coverings of $S^3$ (see Introduction), can be formulated in
terms of labelled link diagrams, with labels taken from the three transpositions of
$\Sigma_3$, according to the above labelling rules at crossings.

The extension from branching links to branching embedded graphs is straightforward.
In fact, we only need to take into account extra labelling constraints and labelled
moves at the vertices of the graph.

Finally, let us consider a $d$-fold covering $p:M \to B^4$ branched over a ribbon
surface $F \subset B^4$. Again, we represent the monodromy in terms of the Wirtinger
presentation of $\pi_1(B^4 - F)$ associated to a locally oriented diagram of $F$.
Actually, since we will only consider simple coverings, we will never need local
orientations. 

The same labelling rules as above apply to ribbon intersections (cf. Figure
\ref{movesR/fig}) as well as to ribbon crossings. However, contrary to what happens
for ribbon intersections, when a ribbon crosses under another one, its label changes
only locally (at the undercrossing region). We notice that, if $F \subset B^4$ is a
labelled ribbon surface representing a $d$-fold (simple) covering of $p:M \to B^4$,
then $L = F \cap S^3$ is a labelled link representing the restriction $p_{|\Bd}: \Bd
M \to S^3$. This is still a $d$-fold (simple) covering, having the diagram of $F$ as
a splitting complex.

As mentioned in the Introduction, labelled ribbon surfaces in $B^4$ (that is
coverings of $B^4$ branched over ribbon surfaces) represent all the 4-dimensional
2-handlebodies. By Montesinos \cite{Mo78} (cf. next Section \ref{todiagram/sec}),
for the connected case it suffices to take labels from the three transpositions of
$\Sigma_3$ (that is to consider 3-fold simple coverings).

Though labelled isotopy of branching ribbon surfaces preserves the covering manifold
$M$ up to PL homeomorphisms, we are interested in the (perhaps more restrictive)
notion of {\sl labelled 1-isotopy}, which preserves $M$ up to 2-deformations (cf.
Lemma \ref{1-isotopy/2-equiv/thm}). This can be realized by means of labelled
diagram isotopy and labelled 1-isotopy moves, that is diagram isotopy and
1-isotopy moves of Figure \ref{isotopy/fig}, suitably labelled according to the
above rules.

\paragraph{Kirby diagrams}

A {\sl Kirby diagram} describes an orientable 4-dimensional 2-handlebody $H^0
\cup\break H^1_1 \cup \dots \cup H^1_m \cup H^2_1 \cup \dots \cup H^2_n$ with only
one 0-handle, by encoding 1- and 2-handles in a suitable link $K \subset S^3 \cong
\Bd H^0$. Namely, $K$ has $m$ dotted components spanning disjoint flat disks which
represent the 1-handles and $n$ framed components which determine the attaching maps
of the 2-handles. We refer to \cite{Ki89} or \cite{GS99} for details and basic facts
about Kirby diagrams, limiting ourselves to recall here only the relevant ones for
our purposes.

The assumption of having only one 0-handle is not so restrictive. In fact, given any
connected handlebody, the union of 0- and 1-handles contracts in a natural way to a
connected graph $G$. Then, by choosing a maximal tree $T \subset G$ and fusing all
the 0-handles together with the 1-handles corresponding to the edges of $T$, we get a
new handlebody with only one 0-handle. This fusion process can be performed by
1-handle slidings and 0/1-handle cancellation, so the new handlebody is
equivalent to the original one. As a consequence, different choises of the tree $T$
give raise to handlebodies which are equivalent up 1-handle sliding. This fact
immediately implies that $k$-equivalence between handlebodies having only one
0-handle can be realized without adding any extra 0-handle.

On the other hand, the same assumption of having only one 0-handle, is crucial in
order to make a natural convention on the framings, that allows to express them by
integers fixing as zero the homologically trivial ones.

However, at least in the present context, it seems preferable to renounce this
advantage on the notation for framings in favour of more flexibility in the
representation of multiple 0-handles. The reason is that a $d$-fold covering of $B^4$
branched over a ribbon surface (actually an embedded 2-dimensional 1-handlebody)
turns out to have a natural handlebody structure with $d$ 0-handles.

Of course, the reduction to only one 0-handle is still possible but it must be
performed explicitly. This makes the connection between branched coverings and
ordinary Kirby diagrams more clear and transparent than before.

\medskip

We call a generalized Kirby diagram our representation of an orientable
4-di\-mensional 2-handlebody with multiple 0-handles. It is essentially defined by
overlapping the boundaries of all the 0-handles to let the diagram take place in
$S^3$ and by putting labels in the diagram in order to keep trace of the original
0-handle where each part of it is from. If there is only one 0-handle, the labels
can be omitted and we have an ordinary Kirby diagram.

More precisely, a {\sl generalized Kirby diagram} representing an orientable
4-di\-mensional handledoby $H^0_1 \cup \dots \cup H^0_d \cup H^1_1 \cup \dots \cup
H^1_m \cup H^2_1 \cup \dots \cup H^2_n$ consists of the following data: a boxed
label indicating the number $d$ of 0-handles; $m$ dotted unknots spanning disjoint
flat disks, each side of which has a label from $\{1, \dots, d\}$; $n$ framed
disjoint knots transversal with respect to those disks, with a label from $\{1,
\dots, d\}$ for each component of the complement of the intersections with the disks.
The labelling must be admissible in the sense that all the framed arcs coming out
from one side of a disks have the same label of that side (cf. Figure
\ref{diag1/fig}). This rule makes the labelling redundant and some times we will
omit the superfluous labels. Moreover, being uniquely related to the indexing
of the 0-handles, the labelling must be considered defined up to permutation of $\{1,
\dots, d\}$. Finally, the framings are always drawn as parallel curves, hence no
confusion arises with labels.

\begin{Figure}[htb]{diag1/fig}{}{}
\centerline{\fig{}{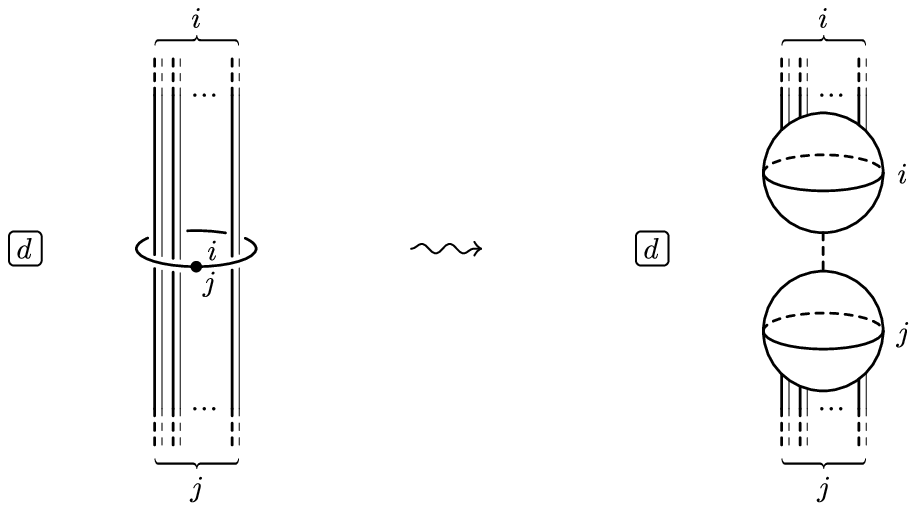}}
\end{Figure}

To establish the relation between a generalized Kirby diagram and the handlebody it
represents, we first convert dot notation for 1-handles into ball notation, as
shown in Figure \ref{diag1/fig}. Here, the two balls, together with the relative
framed arcs, are symmetric with respect to the horizontal plane containing the
disk and squeezing them vertically on the disk we get back the original diagram.
After that, we consider the disjoint union of 0-handles $H^0_1 \cup \dots \cup H^0_d$
and draw on the boundary of each $H^0_i$ the portion of the diagram labelled with
$i$, no matter how we identify such boundary with $S^3$. Then, we attach to $H^0_1
\cup \dots \cup H^0_d$ a 1-handle between each two paired balls (possibly lying in
different 0-handles), according to the diffeomorphism induced by the above symmetry,
so that we can join longitudinally along the handle the corresponding framed arcs.
Of course, the result turns out to be defined only up to 1-handle full twists.
At this point, we have a 1-handlebody $H^0_1 \cup \dots \cup H^0_d \cup H^1_1 \cup
\dots \cup H^1_m$ with $n$ framed loops in its boundary and we use\break such
framed loops as attaching instructions for the 2-handles $H^2_1, \dots, H^2_n$.

We observe that any orientable 4-dimensional 2-handlebody can be represented, up to
isotopy, by a generalized Kirby diagram. In fact, in order to reverse our
construction, we only need that the identification of the boundaries of the
0-handles with $S^3$ is injective on the attaching regions of 1- and 2-handles and
that the attaching maps of the 2-handles run longitudinally along the 1-handles.
These properties can be easily achieved by isotopy.

\medskip

\begin{Figure}[htb]{diag2/fig}{}{}
\centerline{\fig{}{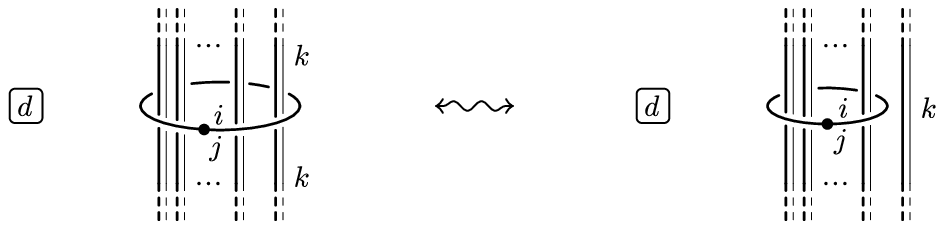}}
\end{Figure}
Sometimes, it will be convenient to derogate from the prescribed labelling rule for
generalized Kirby diagrams, by allowing a framed component with label $k$ to cross a
disk spanned by a dotted component with labels $i$ and $j$, provided that $k \notin
\{i,j\}$. Clearly, such a crossing does not mean that the framed loop goes over the
1-handle corresponding to the dotted one, since it originates from the
identification of different 0-handles. Figure \ref{diag2/fig} depictes the way to
eliminate it.

The above construction gives isotopic handlebody structures if and only if the
starting generalized Kirby diagrams are equivalent up to {\sl labelled isotopy},
generated by labelled diagram isotopy, preserving all the intersections between
loops and disks (as well as labels), and by the three moves described in Figure
\ref{diag3/fig}. Here, we assume $k \neq l$, so that the crossing change at the
bottom of the figure preserves the isotopy class of the framed link in $H^0_1 \cup
\dots \cup H^0_d \cup H^1_1 \cup \dots \cup H^1_m$. 
It is worth remarking\break that, due to this crossing change, the framing convention
usually adopted for ordinary Kirby diagrams cannot be extended to generalized Kirby
diagrams.

\begin{Figure}[htb]{diag3/fig}{}{}
\centerline{\fig{}{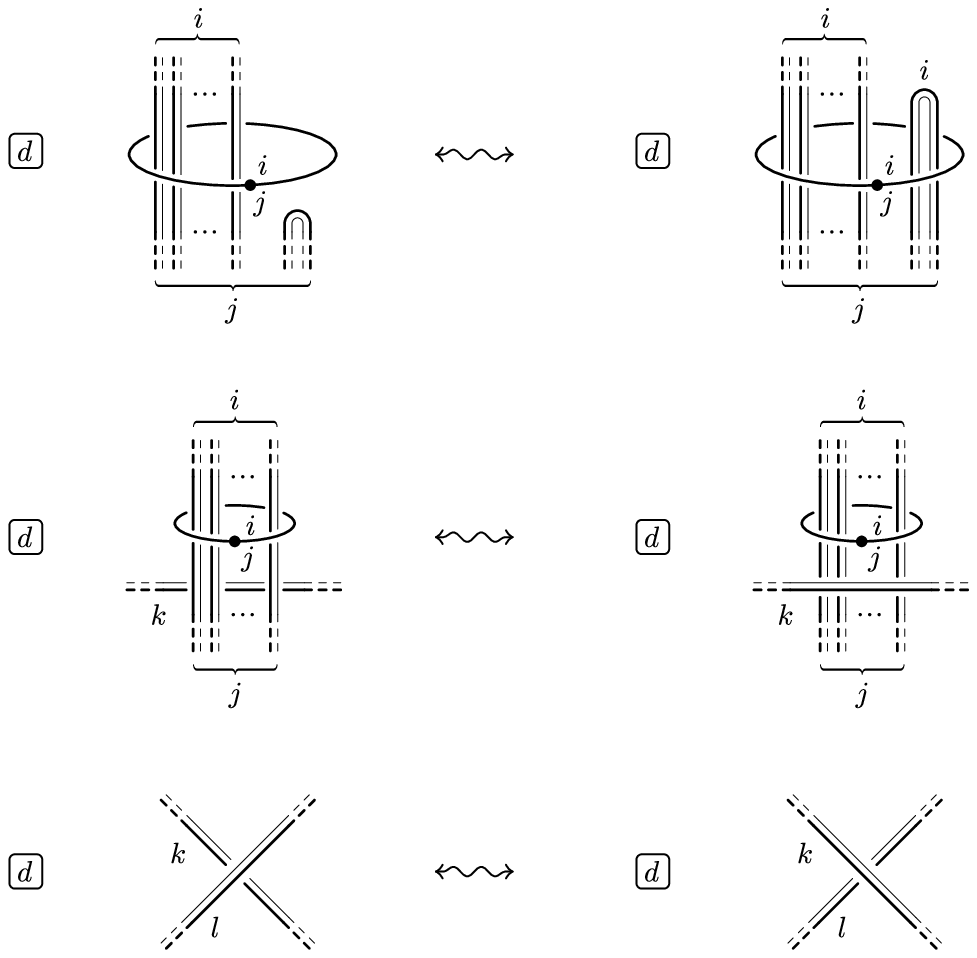}}
\end{Figure}

On the contrary, the other two moves make sense whatever are $i$, $j$ and $k$.
In particular, if $i = j = k$ they reduce to the ordinary ones. Actually, this is the
only relevant case for the second move, usually referred to as ``sliding a 2-handle
over a 1-handle'', being the other cases obtainable by crossing changes.
Moreover, even this ordinary case becomes superfluous in the context of
2-deformations, since it can be realized by addition/deletion of cancelling
1/2-handles and 2-handle sliding\break (cf. \cite{GS99}).

\medskip

The following Figures \ref{diag4/fig} and \ref{diag5/fig} show how to represent
2-deformations of 4-dimensional 2-handlebodies in terms of generalized Kirby
diagrams. Namely, the moves of Figure \ref{diag4/fig} correspond to
addition/deletion of cancelling 0/1-handles (on the right side we assume $i \leq d$)
and 1/2-handles, while the moves of Figure \ref{diag5/fig} correspond to 1- and
2-handle sliding. Except for the addition/deletion of cancelling 0/1-handles, which
does not make sense for ordinary Kirby diagrams, also the rest of the moves reduce
to the ordinary ones if $i=j=k$.

\begin{Figure}[htb]{diag4/fig}{}{}
\vskip6pt\centerline{\fig{}{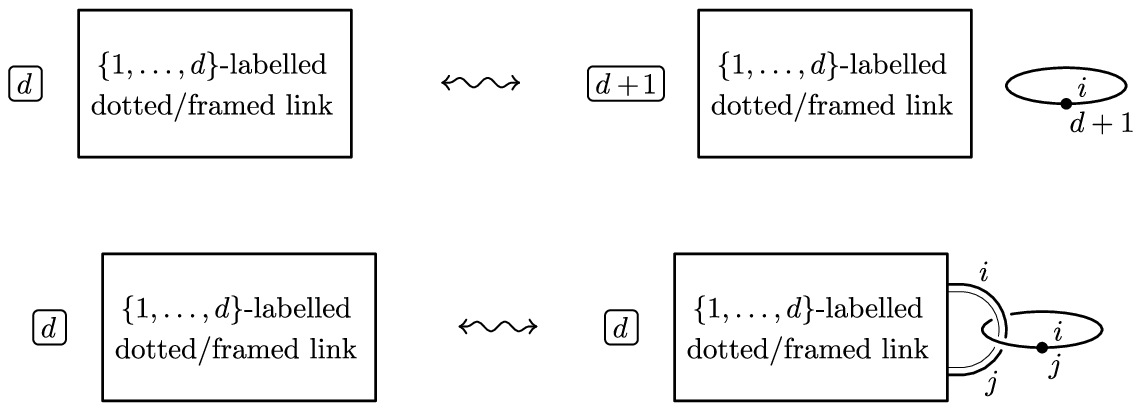}}
\end{Figure}

\begin{Figure}[htb]{diag5/fig}{}{}
\vskip-6pt\centerline{\fig{}{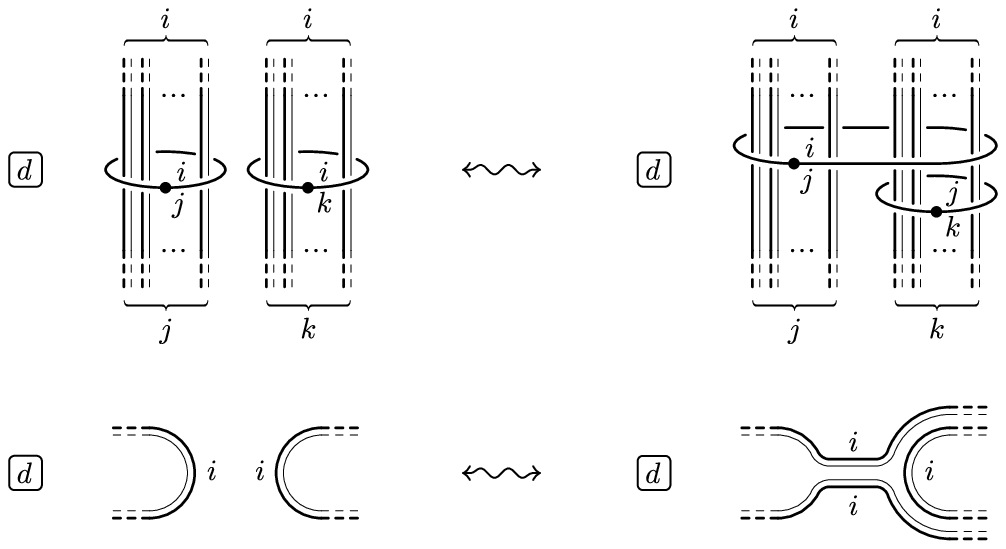}}
\end{Figure}

The 1-handle sliding is included for the sake of completeness, but it can be
generated by addition/deletion of cancelling 1/2-handles and 2-handle sliding, just
like in the ordinary case (cf. \cite{GS99}).

Summing up, two generalized Kirby diagrams represent 2-equivalent 4-dimen\-sional
2-handlebodies if and only if they can be related by the first and third moves of
Figure \ref{diag3/fig} (labelled isotopy), the two moves of Figure \ref{diag4/fig}
(addition/deletion of cancelling handles) and the second move of Figure
\ref{diag5/fig} (2-handle sliding).

Of these, only the first move of Figure \ref{diag3/fig} and the second ones of
Figures \ref{diag4/fig} and \ref{diag5/fig} (for $i = j = d = 1$) make sense in the
case of ordinary Kirby diagrams. Actually, such three moves suffice to realize
2-equivalence of 4-dimensional 2-handlebodies with only one 0-handle, since any
extra 0-handle occurring during a 2-deformation can be eliminated by a suitable
fusion of 0-handles.

\medskip

The main theorem of Kirby calculus \cite{Ki78} asserts that two orientable
4-dimen\-sional 2-handlebodies have diffeomorphic boundaries if and only if they are
related by 2-deformations, blowing up/down and 1/2-handle trading.

In terms of generalized Kirby diagrams these last two modifications can be realized
by the moves of Figure \ref{diag6/fig}. These moves essentially coincide with the
corresponding ones for ordinary Kirby diagrams (with $i = d = 1$), being the involved
labels all the same. 

\begin{Figure}[htb]{diag6/fig}{}{}
\vskip6pt\centerline{\fig{}{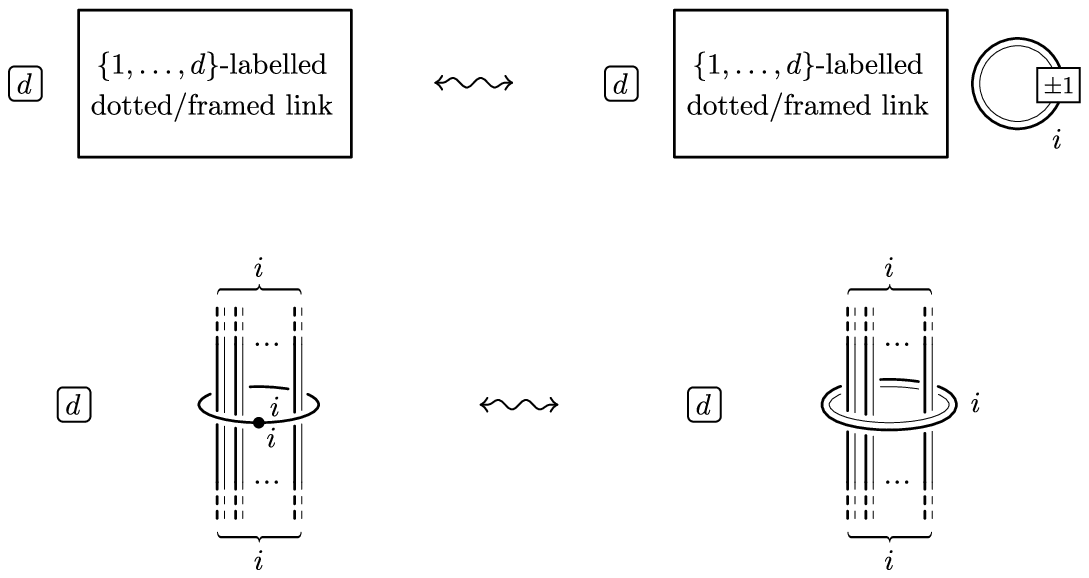}}
\end{Figure}

\begin{Figure}[b]{diag7/fig}{}{}
\vskip6pt\centerline{\fig{}{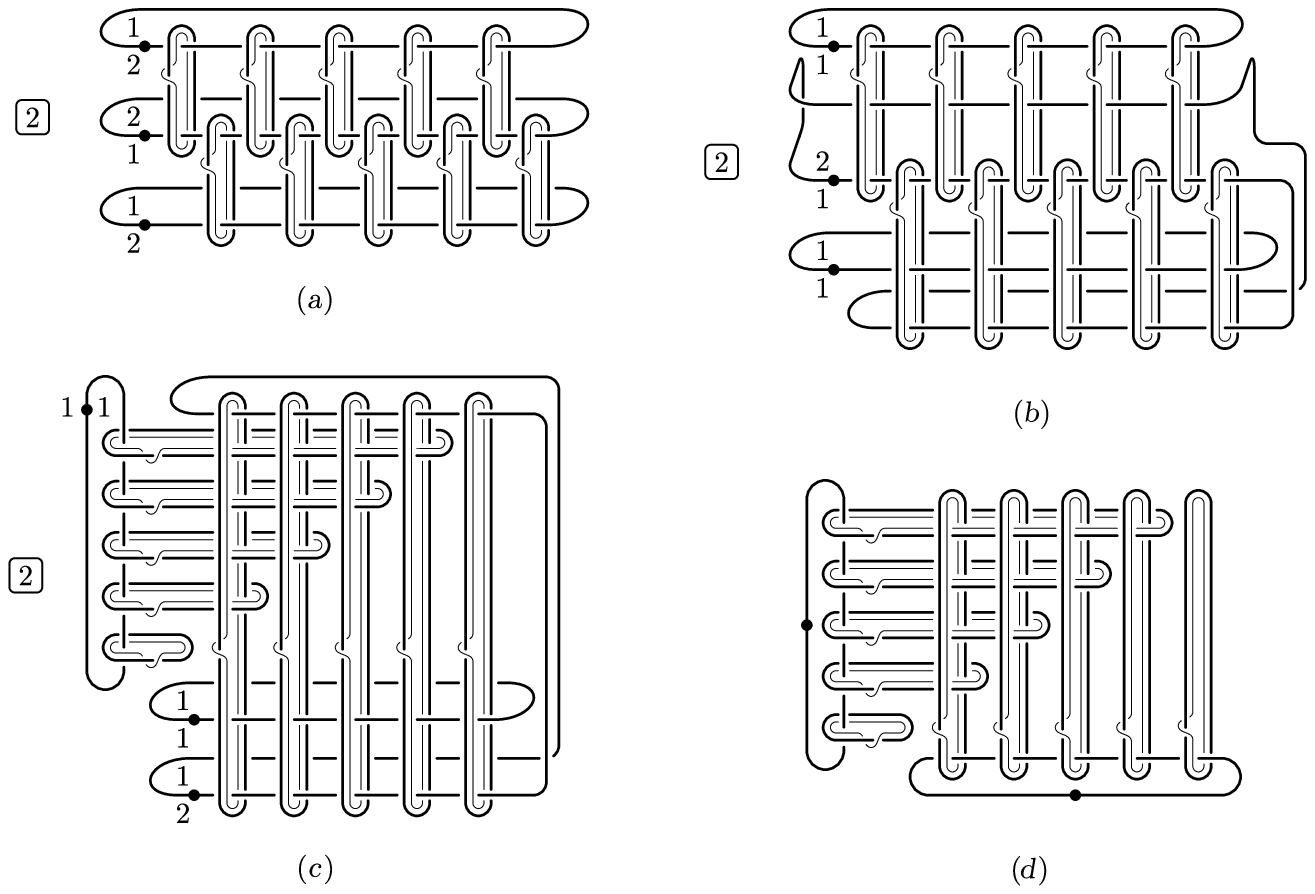}}
\end{Figure}

We conclude this paragraph, by coming back to ordinary Kirby diagrams and in
particular by introducing the standard form that will be used in Section
\ref{toribbon/sec}.

First let us observe that, given any generalized Kirby diagram representing a
connected handlebody, we can use 2-deformation moves to transform it into an
ordinary one, by reducing the number of 0-handles to 1. In fact, assuming $d > 1$, we
can eliminate the $d$-th handle as follows (see Figure \ref{diag7/fig} for an
example with $d = 2$): perform 1-handle sliding in order to leave only one label of
one dotted unknot equal to $d$; untangle such unknot from the rest of the diagram by
labelled isotopy; eliminate the $d$-th 0-handle by 0/1-handle cancellation.

\begin{Figure}[htb]{diagst1/fig}{}{}
\centerline{\fig{}{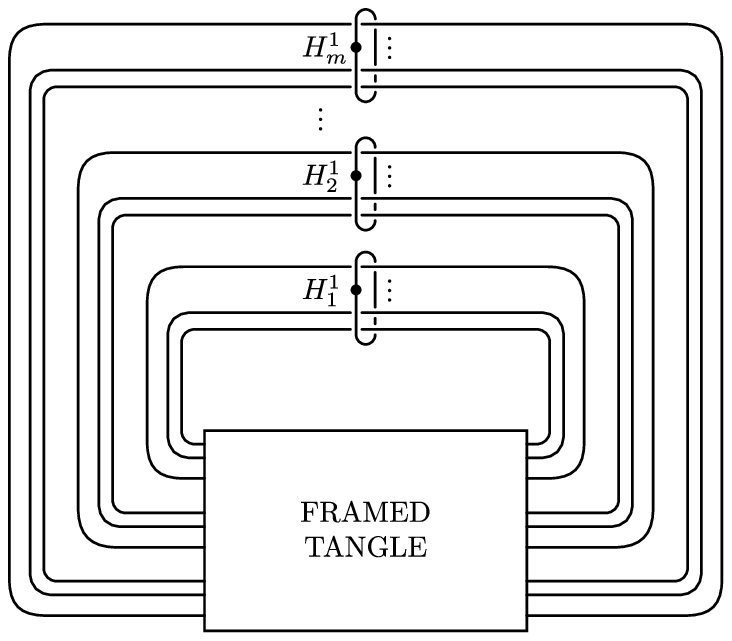}}
\end{Figure}

\begin{Figure}[htb]{diagst2/fig}{}{}
\vskip-6pt\centerline{\fig{}{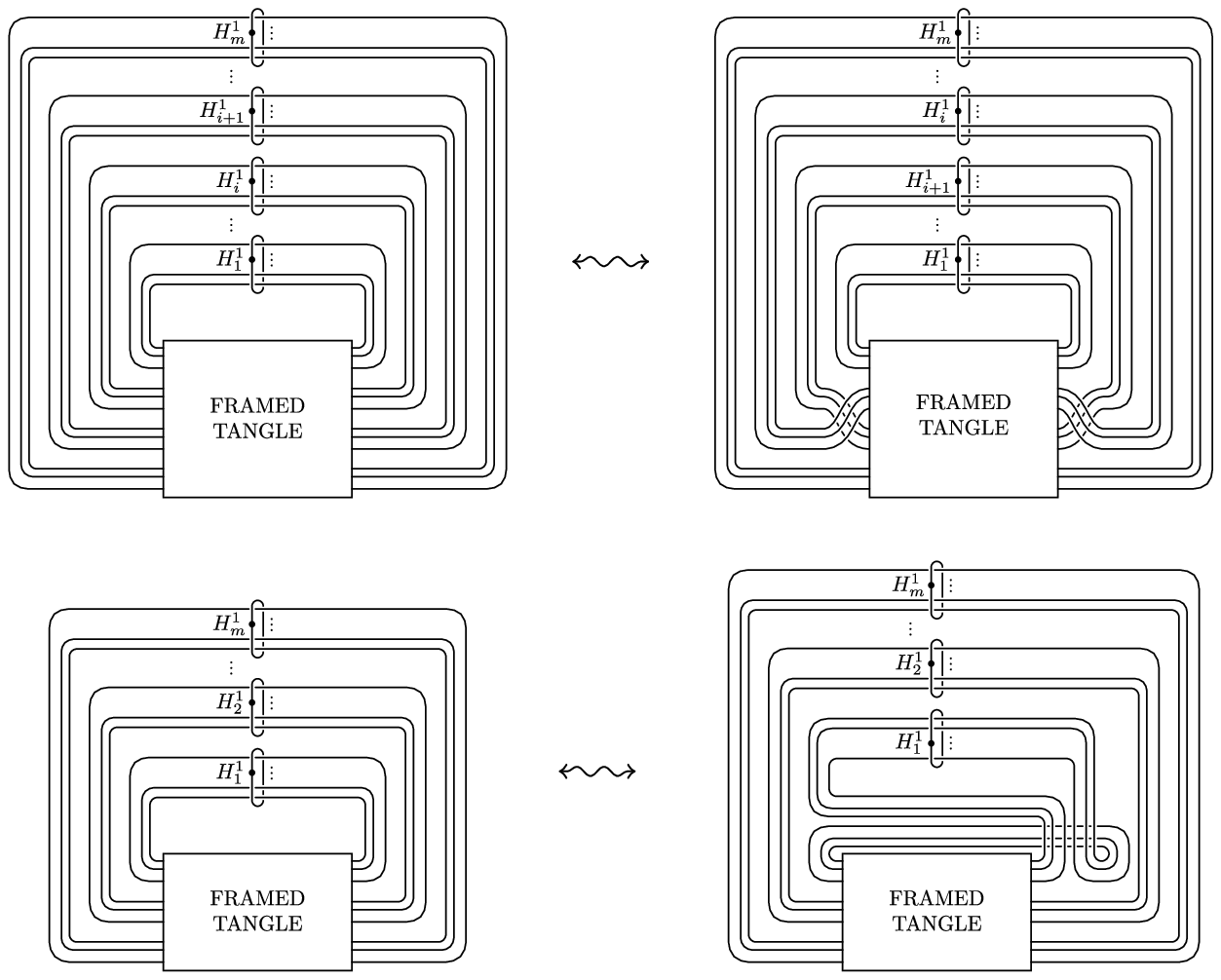}}
\end{Figure}

An ordinary Kirby diagram is said to be in {\sl standard form} if it looks like in
Figure \ref{diagst1/fig}, where all the framings are understood to coincide with the
blackboard one outside the box. Apparently, any ordinary Kirby diagram can be
isotoped into such a standard form. Moreover, isotopy between Kirby diagrams in
standard form consists of framed isotopy fixing the dotted disks together with the
two moves of Figure \ref{diagst2/fig}, which relate different choises for the
vertical order and the orientation of the 1-handles.

It is straightforward to see that, apart from the move of Figure \ref{diagst3/fig}
where an arc of the tangle is isotoped outside the box to pass between two the
dotted disks $H^1_i$ and $H^1_{i+1}$, any isotopy between Kirby diagrams in standard
form which fix the dotted disks can be assumed to move only the framed tangle inside
the box.

\medskip
\begin{Figure}[htb]{diagst3/fig}{}{}
\vskip4pt\centerline{\fig{}{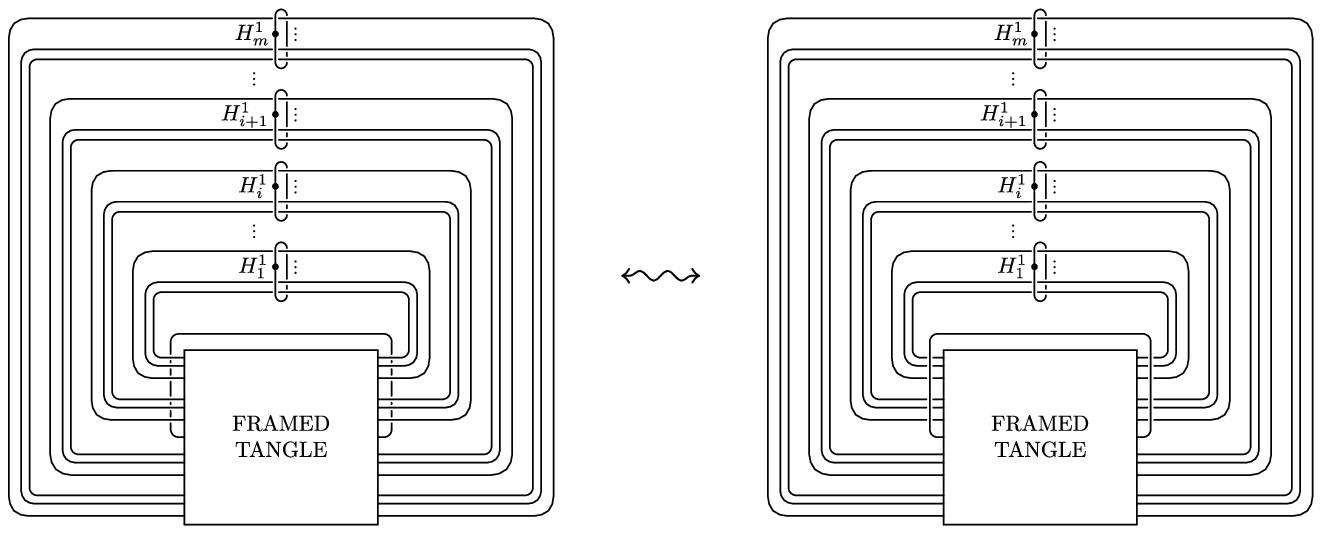}}
\end{Figure}

Finally, in the following proposition, we recapitulate the moves relating ordinary
Kirby diagrams which are 2-equivalent or have diffeomorphic boundaries.

\begin{proposition} \label{stmoves/thm}
Given two ordinary Kirby diagrams $K$ and $K'$ in standard form, denote by $H$
and $H'$ the corresponding 4-dimensional 2-handlebodies. Then:
\begin{itemize}\itemsep0pt
\item[\(a)] \vskip-\lastskip\smallskip
$H$ and $H'$ are 2-equivalent if and only if $K$ and $K'$ are equivalent up to
the first move of Figure \ref{diag3/fig}, the second ones of Figures \ref{diag4/fig}
and \ref{diag5/fig}, the two moves of Figure \ref{diagst2/fig} and framed tangle
isotopy inside the box;
\item[\(b)]
$H$ and $H'$ have diffeomorphic boundaries if and only if $K$ and $K'$ are
equivalent up to the moves listed in \(a), the two moves of Figure \ref{diag6/fig}
and framed tangle isotopy inside the box.
\end{itemize}\vskip-\lastskip\smallskip
All the moves are understood to preserve the standard form (up to diagram isotopy)
in the obvious way and the ones of Figures \ref{diag3/fig}, \ref{diag4/fig},
\ref{diag5/fig} and \ref{diag6/fig} are considered only in the ordinary case, that is
for $i = j = d = 1$.
\end{proposition}

\begin{proof}
It is enough to prove \(a), being \(b) an immediate consequence of it and of the main
theorem of Kirby calculus. To do that, we only need to show that the moves listed in
\(a) allow us to represent any 2-deformation between $H$ and $H'$.

As discussed above, 2-deformations of ordinary Kirby diagrams can be represented
by the ordinary cases of the first move of Figure \ref{diag3/fig} and of the second
ones of Figures \ref{diag4/fig} and \ref{diag5/fig}, together with diagram isotopy.
Moreover, up to diagram isotopy, any of these moves can be performed on a Kirby
diagram in standard form, in such a way that the standard form is preserved.

On the other hand, we already observed that the moves of Figures \ref{diagst2/fig}
and \ref{diagst3/fig}, together with framed tangle isotopy inside the box, allows us
to realize diagram isotopy between Kirby diagrams in standard form.

Then, to conclude the proof, it suffices to notice that the move of Figure
\ref{diagst3/fig} can be obtained by sliding of a 2-handle over some 1-handles (cf.
second move of Figure \ref{diag3/fig}). Hence, in the context of 2-deformations
diagram isotopy between Kirby diagrams in standard form can be realized without it.
\end{proof}

\section{From labelled ribbon surfaces to Kirby diagrams\label{todiagram/sec}}

The aim of this section is to show how any adapted 1-handlebody structure on a
labelled ribbon surface $F$ representing a $d$-fold simple branched covering $p:M
\to B^4$ naturally induces a 2-handlebody structure on $M$ defined up to
2-deformations. 

In this context, naturally means that labelled embedded 1-deformations on $F$ induce
2-deformations on $M$.  Then, by Propositions \ref{1-handles/thm} and
\ref{1-isotopy/thm}, $M$ turns out to be endowed with a 2-handlebody structure,
whose 2-equivalence class is uniquely determined by the labelled 1-isotopy class of
$F$. We denote by $K_F$ the generalized Kirby diagram corresponding to such
4-dimensional 2-handlebody structure (defined up to 2-deformations).

Moreover, we will see that the 2-equivalence class of $K_F$ is also preserved by the
covering moves $R_1$ and $R_2$ of Figure \ref{movesR/fig} and we will discuss some
consequences of this fact. In particular, we will introduce some auxiliary moves
generated by $R_1$ and $R_2$, that will be needed in the next sections. 

\medskip

Let us start with the construction of $K_F$. Given a labelled ribbon surface $F$ as
above with an adapted 1-handlebody decomposition, we can write $F = D_1 \cup \dots
\cup D_m \cup B_1 \cup \dots \cup B_n$, where the $D_h$'s are disjoint flat disks
(the 0-handles of $F$) while the $B_h$'s are disjoint bands attached to $F_0 = D_1
\cup \dots \cup D_m$ (the 1-handles of $F$). 
Looking at the diagram of $F$ in $R^3$ and using for it the same notations as for
$F$ itself, we see that the $D_h$'s, as well as the $B_h$'s, are still disjoint
from each other, while any band $B_h$ may form ribbon intersections with the disks
$D_1, \dots, D_m$.

We denote by $p_0 : M_0 \to B^4$ the simple covering determined by the labelled
surface $F_0 \subset B^4$. The covering manifold $M_0$ turns out to be a
4-dimensional handlebody with $d$ 0-handles and a 1-handle $H^1_h$ for each
disk $D_h$ (cf. \cite{Mo78}). 
A generalized Kirby diagram of $M_0$ can be immediately obtained by replacing any
disk $D_h$ by a dotted unknot coinciding with its boundary, as shown in Figure
\ref{covdiag1/fig}.
Here, there are two possible way to assign the labels $i$ and $j$ to the two faces of
$D_h$. We call such an assignment a {\sl polarization} of the disk $D_h$.

\begin{Figure}[htb]{covdiag1/fig}{}{}
\centerline{\fig{}{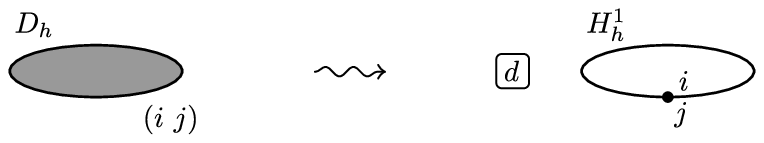}}
\end{Figure}

Now, following \cite{Mo78} (cf. also \cite{IP02}), we have that any band $B_h$
attached to $F_0$ gives rise to a 2-handle $H^2_h$ attached to $M_0$ along the framed
loop given by the unique annular component of $p_0^{-1}(B_h)$. 

In order to describe a labelled framed loop representing $H^2_h$ in the generalized
Kirby diagram, let us call $D_{h_1}$ and $D_{h_2}$ the (possibly coinciding)
disks of $F_0$ at which $B_h$ is attached. Disregarding for the moment the
ribbon intersections of $B_h$ with $F_0$, such framed loop is given by two parallel
copies of $B_h$ lying on opposite sides, joined together to form ribbon
intersections with $D_{h_1}$ and $D_{h_2}$ and labelled consistently with the
polarizations of those disks, as suggested by Figure \ref{covdiag2/fig}.

\begin{Figure}[htb]{covdiag2/fig}{}{}
\centerline{\fig{}{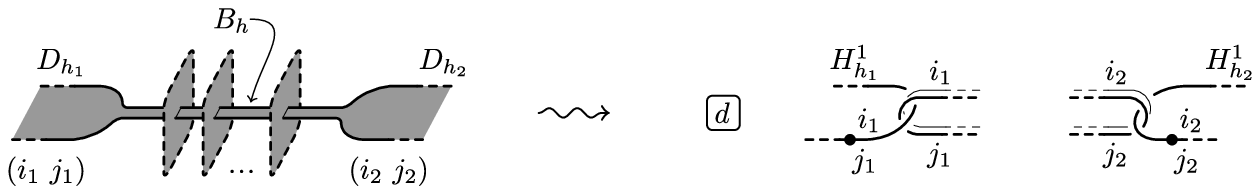}}
\end{Figure}

\begin{Figure}[htb]{covdiag3/fig}{}{}
\vglue0.5pt\centerline{\fig{}{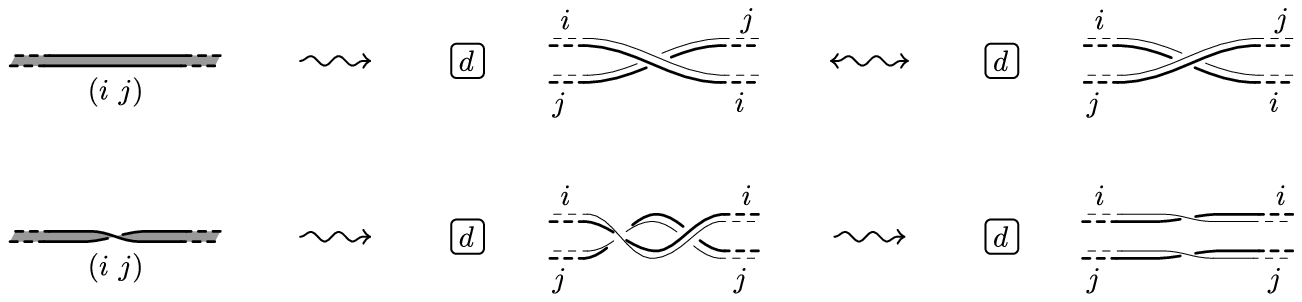}}
\end{Figure}

\begin{Figure}[htb]{covdiag4/fig}{}{}
\centerline{\fig{}{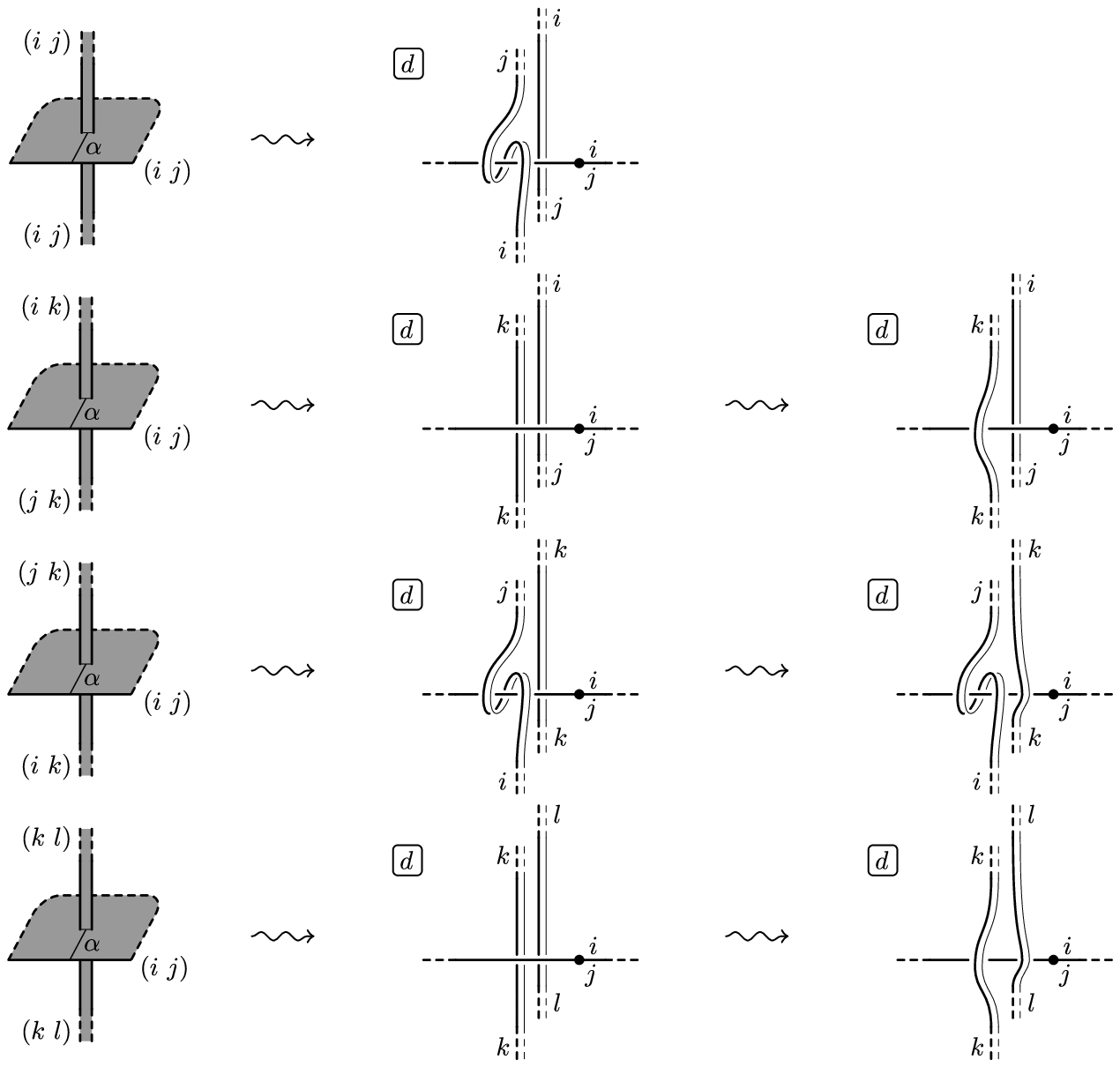}}
\end{Figure}

Actually, to have simultanous labelling consistency at both ends of $B_h$, we may be
forced to interchange the two copies of $B_h$ by a crossover, as in the upper part of
Figure \ref{covdiag3/fig}. The two ways to realize the crossover are equivalent
up to labelled isotopy, since $i\neq j$. Notice that we can perform crossovers
wherever we want along $B_h$, provided their number has the right parity to respect
labelling consistency.
In the lower part of Figure \ref{covdiag3/fig} we see that, up to crossovers, twists
along $B_h$ contribute only to the framing and not to the isotopy type of the
corresponding loop. Namely, each positive (resp. negative) half twist along $B_h$
gives rise to a positive (resp. negative) full twist in the framing.

Figure \ref{covdiag4/fig} explains how to interpret a single ribbon intersection
between $B_h$ and $F_0$ into the generalized Kirby diagram, in the four possible
cases depending on the  monodromies associated to $B_h$ and $F_0$ at that
intersection. Here, we assume that $i$, $j$, $k$ and $l$ are all distinct and use
the notation of Figure \ref{diag2/fig} for the intermediate steps. In all the cases,
the construction is carried out in a regular neighborhood of an arc $\alpha$
contained in $F_0$ and joining the ribbon intersection with $\Bd F_0$.\break 
The labels of the two copies of $B_h$ in the generalized Kirby diagram are
determined by monodromies associated to $B_h$ before and after the intersection 
and by the side from which $\alpha$ approaches $B_h$.
In the first and third cases, we introduce a kink to allow labelling consistency of
both the copies of $B_h$ with respect to the disk (there are two different ways to
realize such a kink, but they are equivalent up to labelled isotopy).
To make all the local labellings at the ribbon intersections fit together with each
other along the two copies of $B_h$ and with the ones already fixed at ends of
$B_h$, we use again crossovers.

We conclude the definition of $K_F$, by specifying that the arcs $\alpha$ related to
dif\-ferent ribbon intersections are assumed to be disjoint, in such a way that the
corresponding constructions do not interact.

Our next aim is to show that $K_F$ is well defined up to 2-deformation moves, in
the sense that the 2-equivalence class of the corresponding 4-dimensional
2-handlebody depends only on the labelled ribbon surface $F$. As a preliminary step,
we prove the following Lemma concerning the choices involved in the construction of
$K_F$ from a 1-handlebody structure of $F$.

\begin{lemma} \label{ribbon/diag/thm}
Let $F \subset B^4$ be a labelled ribbon surface representing a $d$-fold simple
branched covering $p:M \to B^4$. Then, the generalized Kirby diagram $K_F$,
constructed starting from a given adapted 1-handlebody structure on $F$, describes a
4-dimensional 2-handlebody structure on $M$, whatever choices we make for the
polarizations, the crossovers and the arcs $\alpha$. Moreover, such 4-dimensional
2-handlebody structure is uniquely determined up to handle isotopy.
\end{lemma}

\begin{proof} 
Let $F = D_1 \cup \dots \cup D_m \cup B_1 \cup \dots \cup B_n$ be an adapted
1-handlebody decomposition of $F$ as in the definiton of $K_F$ and let us adopt
here all the notations related to it we introduced there.

Then, the lemma immediately follows from \cite{Mo78}, once one has checked that the
framed link of $K_F$ does really represent, up to handle isotopy, the framed link in
$M_0$ consisting of the unique annular component of $p_0^{-1}(B_h)$ for each band
$B_h$ of $F$. Taking into account what we have said above, this is a straightforward 
consequence of the very definition of generalized Kirby diagram.

Nevertheless, for the convenience of the reader, we skecth a direct proof of
the independence of $F_K$, up to handle isotopy, on the choices involved in its
construction.

We have already observed that crossovers are not relevant up to labelled isotopy.
Concerning the arcs $\alpha$, it suffices to prove that the elementary moves of
Figure \ref{covdiag5/fig}, where we replace a single arc $\alpha$ by $\alpha'$,
preserve $K_F$ up to 2-deformation moves. Simple inspection  of all the cases
confirms that once again only labelled isotopy moves are needed.

\begin{Figure}[htb]{covdiag5/fig}{}{}
\centerline{\fig{}{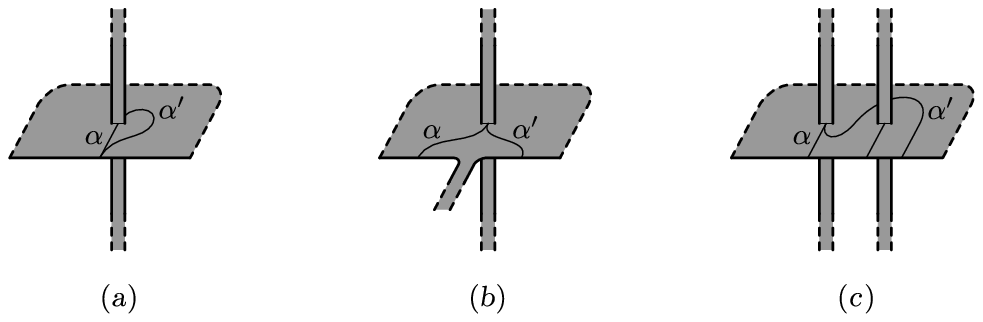}}
\end{Figure}

Thus, it remains to see what happens when we invert the polarization of a disk
$D_h$. The relative dotted unknot with the different labellings giving the two
possible polarizations of $D_h$ is drawn in Figure \ref{covdiag6/fig} \(a) and
\(d).%
\begin{Figure}[htb]{covdiag6/fig}{}{}
\vskip6pt\centerline{\fig{}{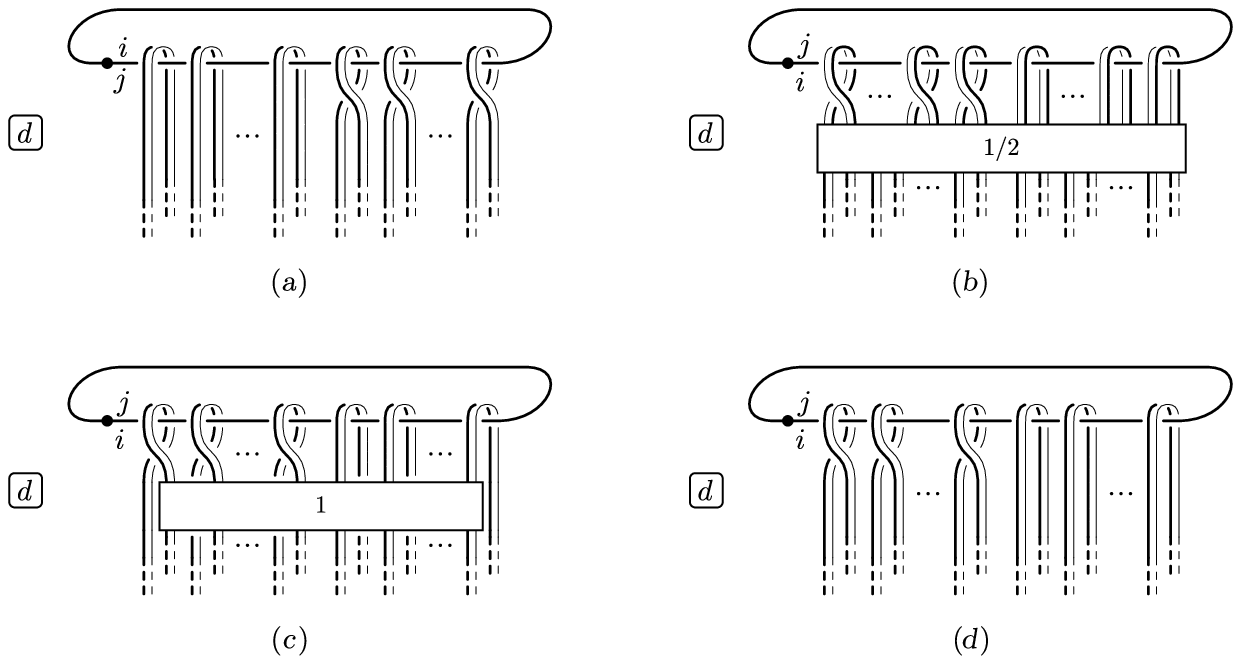}}
\end{Figure}
Here, we assume that the framed arcs passing thought $D_h$, coming either
from bands attached to $D_h$ or from ribbon intersection of bands with $D_h$, have
been isotoped all together into a canonical position. For the sake of clarity, we
sorted such labelled arcs to separate the ones which have been kinked for respecting
labelling consistency.
To see that the diagrams \(a) and \(d) of Figure \ref{covdiag6/fig} are equivalent
up to 2-deformation moves, we consider the other ones as intermediate steps.
We start by isotoping upside down the dotted unknot of \(a) to obtain \(b). Then, we
use labelled isotopy once again to make the arcs labelled by $i$ and the ones
labelled by $j$ form separate positive half twists. These two half twists add up to
give a unique positive full twist in \(c). Finally, we get \(d) by performing a
negative twist on the 1-handle represented by the dotted unknot. Such a 1-handle
twist can be easily realized by the second labelled isotopy move of Figure
\ref{diag3/fig}.
\end{proof}

\begin{proposition} \label{1-equiv/2-equiv/thm}
Let $F \subset B^4$ be a labelled ribbon surface representing a $d$-fold simple
branched covering $p:M \to B^4$. Then, the generalized Kirby diagrams $K_F$
constructed starting from different adapted 1-handlebody structures on $F$, describe
2-equivalent 4-dimensional 2-handlebody structures on $M$. That is, $K_F$ is
uniquely determined by $F$ up to 2-deformation moves.
\end{proposition}

\begin{proof}
We observe that any labelled diagram isotopy (preserving ribbon intersections) on $F$
induces a labelled isotopy on $K_F$ as a generalized Kirby diagram.
Hence, the statement follows from Proposition \ref{1-handles/thm}, once we prove
that performing on $F$ labelled versions of the moves of Figures \ref{ribbon2/fig}
and \ref{ribbon3/fig} without vertical disks corresponds to modifying $K_F$ by
certain 2-deformation moves.

In all the cases we can choose the same polarization for $H^0_i$ and $H^0_j$, since
these can be assumed to be distinct 0-handles (cf. notice after Figure
\ref{ribbon3/fig}). Then, apparently the two moves of Figure \ref{ribbon3/fig}
correspond respectively to addition/deletion of a cancelling pair of 1/2-handles and
to sliding the 2-handle deriving from $H^1_l$ over the one deriving from $H^1_k$.
Similarly, in the case of move of Figure \ref{ribbon2/fig} we have\break two slidings
involving the same 2-handles, one sliding for each of the two parallel copies of
$H^1_l$ forming the framed loop originated from it. We leave to the reader the
straightforward verification of this fact for all the four cases of Figure
\ref{covdiag4/fig}.
\end{proof}

A very simple example of the above construction, without ribbon intersections, is
depicted in Figure \ref{covdiag7/fig}. Here, the adapted 1-handlebody structure of
the labelled ribbon surface on the left is the obvious one with 3 horizontal
0-handles and 10 vertical 1-handles, while the resulting generalized Kirby diagram on
the right is the same of Figure \ref{diag7/fig}. 
We notice that, for a double covering of $B^4$ branched over ribbon surface without
ribbon intersections, such as the one of Figure \ref{covdiag7/fig}, the handlebody
presentation we obtain by our construction coincides, after suitable reduction to
ordinary Kirby diagram, with the one given in \cite{AK80}.

\smallskip
\begin{Figure}[htb]{covdiag7/fig}{}{}
\vskip4pt\centerline{\fig{}{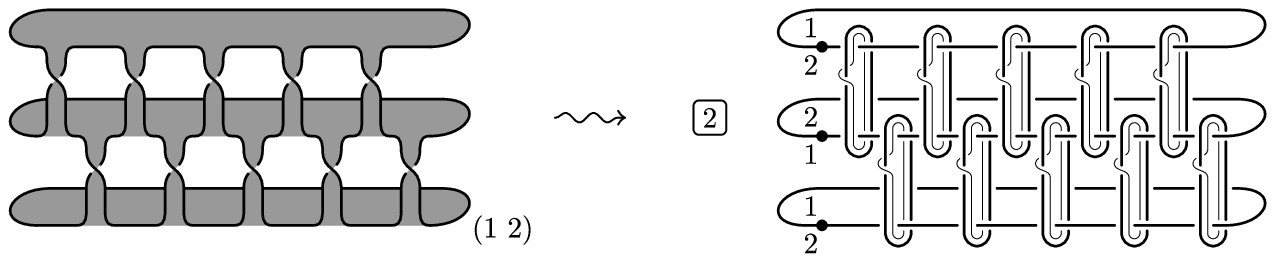}}
\end{Figure}

The following Proposition \ref{movesR/2def/thm} tells us that the 2-equivalence class
of $K_F$ actually depends only on the labelled 1-isotopy class of $F$ and it is also
preserved by stabilization and covering moves $R_1$ and $R_2$. This is
essentially the ``only if'' part of Theorem \ref{equiv4/thm}.

\begin{lemma} \label{1-isotopy/2-equiv/thm}
If the labelled ribbon surfaces $F,F' \subset B^4$, representing $d$-fold simple
branched coverings of $B^4$, are related by labelled 1-isotopy, then the generalized
Kirby diagrams $K_F$ and $K_{F'}$ are equivalent up to 2-deformation moves.
\end{lemma}

\begin{proof}
By Proposition \ref{1-isotopy/thm}, labelled 1-isotopy is generated by labelled
diagram isotopy and the labelled versions of moves $I_1, \dots, I_4$ (cf. Figure
\ref{isotopy/fig}). Since labelled diagram isotopy on $F$ induces labelled isotopy on
$K_F$ as a generalized Kirby diagram, we have only to deal with the moves.

Move $I_1$ admits a unique labelling up to conjugation in $\Sigma_d$. Generalized
Kirby diagrams arising from the labelled ribbon surfaces involved in the
resulting labelled move are depicted in Figure \ref{isotopy1/fig} (we assume the
surfaces endowed with the handlebody structures of the corresponding move of Figure
\ref{ribbon1/fig}). As the reader can easily check, such diagrams are related by
labelled isotopy.

\begin{Figure}[htb]{isotopy1/fig}{}{}
\centerline{\fig{}{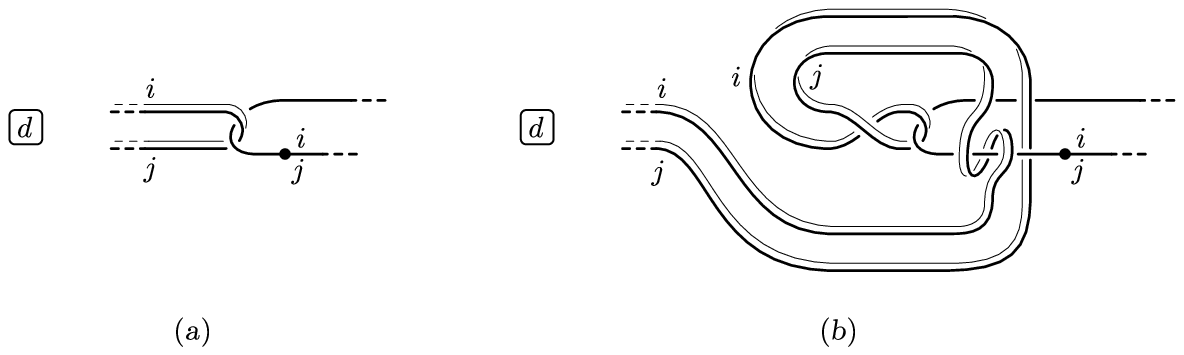}}
\end{Figure}

Moves $I_2$ and $I_3$ admit three distinct labellings up to conjugation in
$\Sigma_d$. Namely, if $(i\ j)$ is the label of the horizontal component, then the
top end of the vertical one can be labelled by $(i\ j)$, $(j\ k)$ or $(k\ l)$. 

\begin{Figure}[b]{isotopy2/fig}{}{}
\centerline{\fig{}{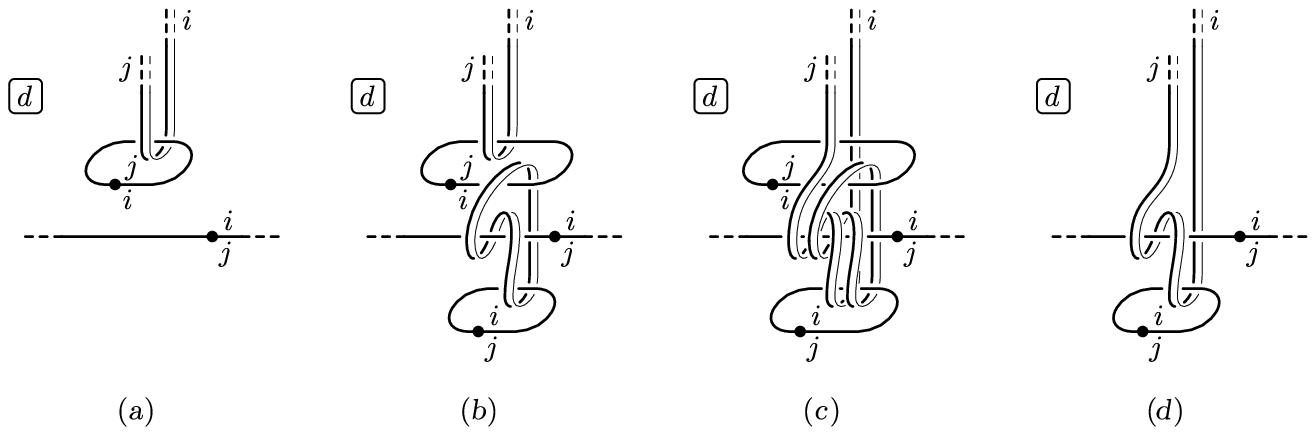}}
\end{Figure}

The first case is considered in Figure \ref{isotopy2/fig} for $I_2$ and Figure
\ref{isotopy3/fig} for $I_3$. Looking at these figures, we have that: \(a) and
\(d) correspond respectively to the surface on the left and right side of the move
with the simplest adapted handlebody structures; \(b) is obtained from \(a) by
1/2-handle addition, followed by 2-handle sliding only in Figure \ref{isotopy3/fig};
\(c) and \(d) are obtained in turn by 2-handle slidings and 1/2-handle cancellation.
The same figures also apply to the second case, after we replace by $k$'s all the
$i$'s in the upper half and the $j$'s in the lower half (except for the labels of the
dotted line in the middle). The third case is trivial and we leave it to the reader.

\break

Finally, let us come to move $I_4$, which requires a bit more work than the other
ones. As above, let $(i\ j)$ be the label of the horizontal band. Then, up to
conjugation in $\Sigma_d$, there are eighteen possible ways to label the move, each
one determined by the transpositions $\lambda$ and $\rho$ labelling respectively the
left and right bottom ends of the diagonal bands.
By direct inspection we see that, excluding the trivial cases when at least two of
the three ribbon intersections involve bands with disjoint monodromies, which are
left to the reader, and taking into account the symmetry of the move with respect to
its inverse, there are only seven relevant cases: 
1)~$\lambda = (i\ j)$ and $\rho = (i\ j)$; 2)~$\lambda = (i\ j)$ and $\rho = (i\ k)$;
3)~$\lambda = (i\ k)$ and $\rho = (i\ j)$; 4)~$\lambda = (i\ k)$ and $\rho = (i\ k)$;
5)~$\lambda = (i\ k)$ and $\rho = (i\ l)$; 6)~$\lambda = (i\ k)$ and $\rho = (j\ l)$;
7)~$\lambda = (i\ k)$ and $\rho = (k\ l)$.

\begin{Figure}[t]{isotopy3/fig}{}{}
\centerline{\fig{}{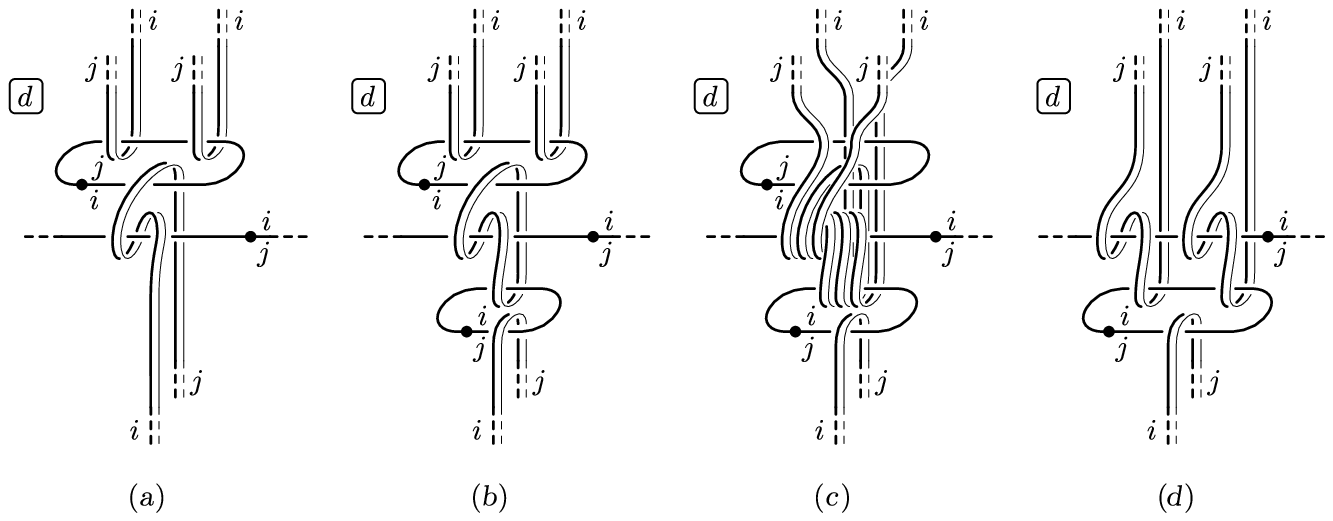}}
\end{Figure}

Figure \ref{isotopy4a/fig} regards case 1. Here, \(a) and \(c) correspond
respectively to the surfaces on the left side and right side of the move with
suitable adapted handlebody structures, while \(b) is related to \(a) by two
2-handle slidings and to \(c) by labelled isotopy. This figure also applies to case
4, after the same label replacement as above.

\begin{Figure}[htb]{isotopy4a/fig}{}{}
\centerline{\fig{}{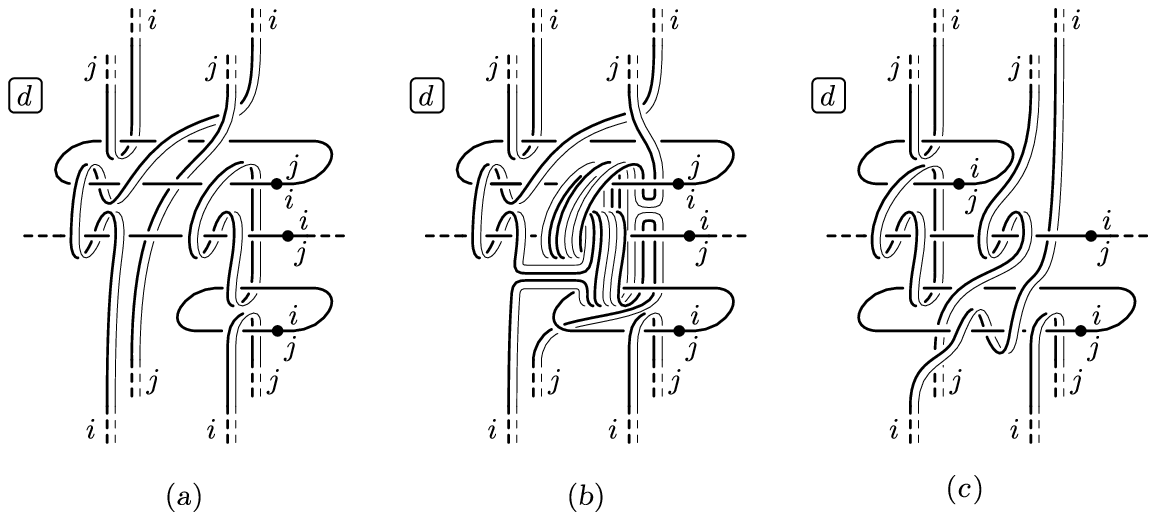}}
\end{Figure}

Similarly, Figure \ref{isotopy4b/fig} concerns with case 2 and, after the
appropriate label replacements, also with cases 3 and 5. This time only one 2-handle
sliding is needed to pass from \(a) to \(b). Figures \ref{isotopy4c/fig} and
\ref{isotopy4d/fig} complete the proof, by dealing with the remaining cases 5 and 7.
The three diagrams of Figure \ref{isotopy4c/fig} are related by 1/2-handle
addition/deletion, while the two diagrams of Figure \ref{isotopy4d/fig} by labelled
isotopy.
\end{proof}

\begin{Figure}[htb]{isotopy4b/fig}{}{}
\centerline{\fig{}{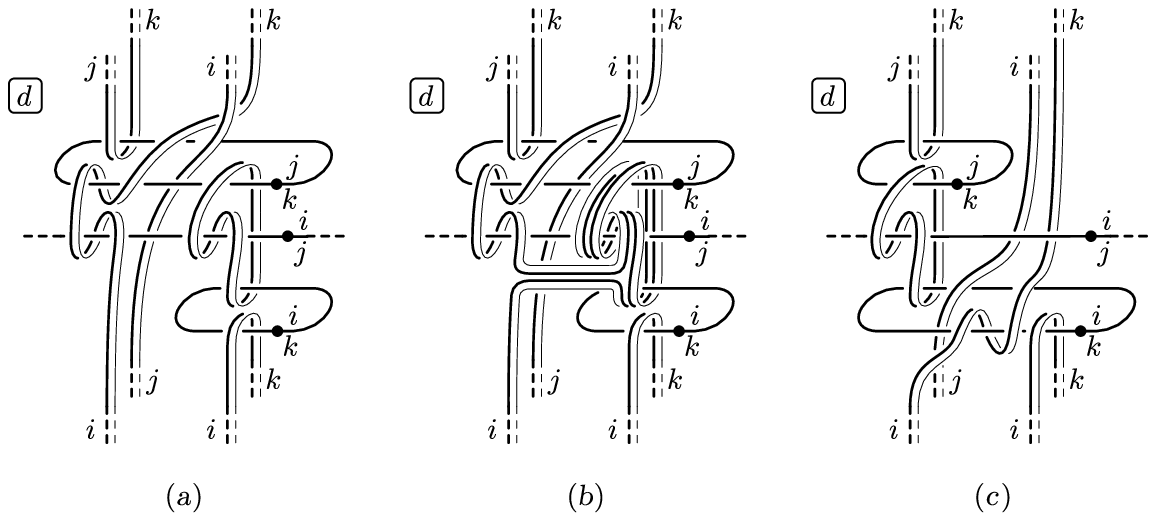}}
\end{Figure}

\begin{Figure}[htb]{isotopy4c/fig}{}{}
\bigskip\centerline{\fig{}{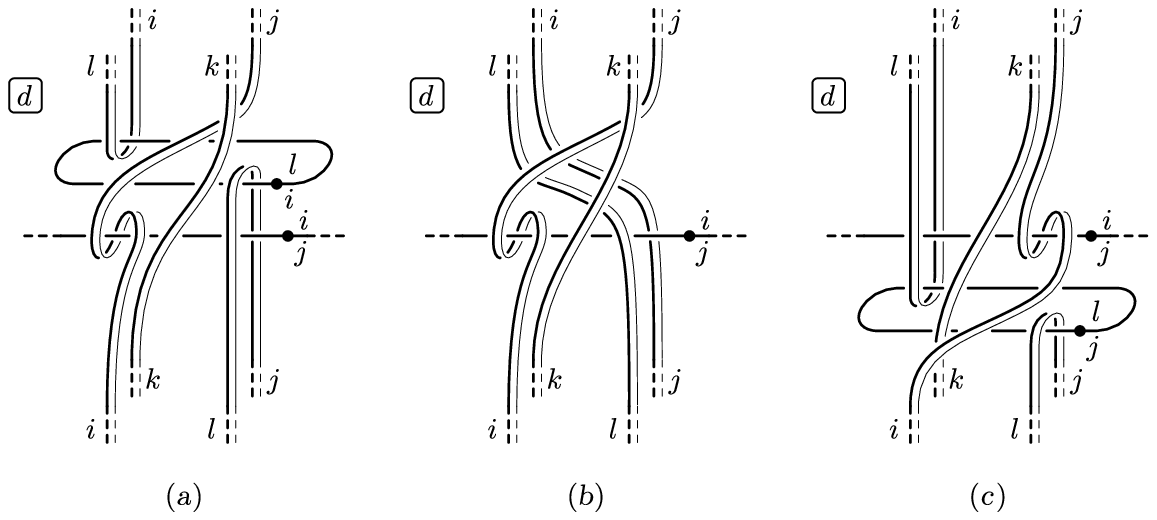}}
\end{Figure}

\begin{Figure}[htb]{isotopy4d/fig}{}{}
\bigskip\centerline{\fig{}{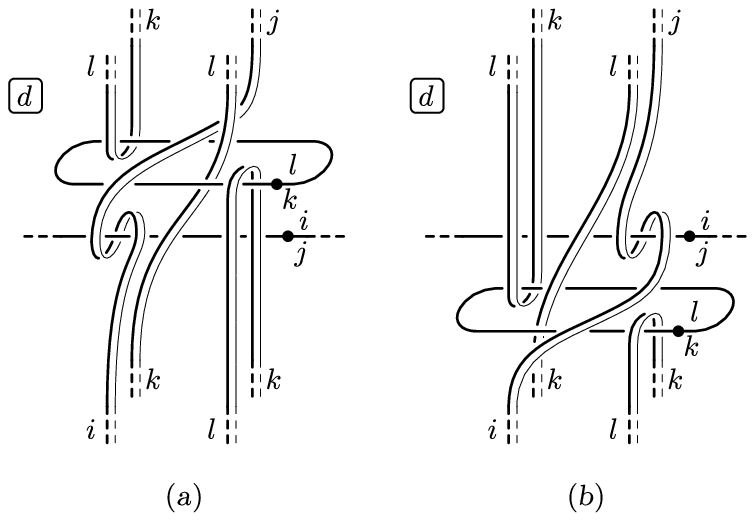}}
\end{Figure}

\break

It is worth remarking that Lemma \ref{1-isotopy/2-equiv/thm} becomes trivial if we
limit ourselves to require that the 4-dimensional 2-handlebodies represented by
$K_F$ and $K_{F'}$ are diffeomorphic, without insisting that they are 2-equivalent.
In fact, labelled isotopy between $F$ and $F'$ (instead of labelled 1-isotopy)
suffices for that, since it induces equivalence between the corresponding branched
coverings, as recalled in Section \ref{prelim/sec}. 
The relation between isotopy and 1-isotopy of ribbon surfaces in $B^4$ on one hand
and diffeomorphism and 2-equivalence of 4-dimensional 2-handlebodies on the other
hand, will be discussed in Section \ref{remarks/sec}.

\begin{proposition} \label{movesR/2def/thm}
If the labelled ribbon surfaces $F,F' \subset B^4$, representing simple branched
coverings of $B^4$, are related by labelled 1-isotopy, stabilization and moves $R_1$
and $R_2$, then the generalized Kirby diagrams $K_F$ and $K_{F'}$ are equivalent up
to 2-deformation moves.
\end{proposition}

\begin{proof}
Labelled 1-isotopy has been already considered in the previous lemma.
From the definitions it is apparent that stabilizing the branched coverings
represented by a labelled ribbon surface $F$ means adding a cancelling pair of
0/1-handles to $K_F$ (cf. Figures \ref{diag4/fig} and \ref{covdiag1/fig}).

Concerning moves $R_1$ and $R_2$, if $F$ and $F'$ differ by such a move, then by
making the right choices in the construction of $K_F$ and $K_{F'}$ we get the same
result up to labelled isotopy.
This is shown in Figure \ref{moveR1def/fig} (to be
compared with Figure \ref{movesR/fig}) for move $R_1$. The analogous and even easier
case of move $R_2$ is left to the reader.
\end{proof}

\begin{Figure}[htb]{moveR1def/fig}{}{}
\centerline{\fig{}{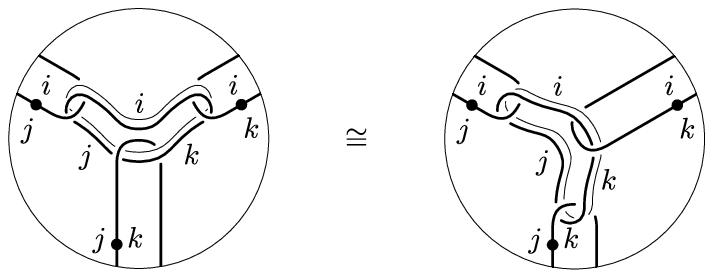}}
\end{Figure}

We conclude this section with some further considerations on the ribbon moves $R_1$
and $R_2$. In particular, we see how they generate, up to labelled 1-isotopy, the
auxiliary moves $R_3$, $R_4$, $R_5$ and $R_6$ described in Figure
\ref{movesRaux/fig}, where $i$, $j$ and $k$ are all distinct. These last moves will
turn out to be useful in the next sections. First, we formalize in the following
proposition the observation made in the Introduction about the inverses of moves
$R_1$ and $R_2$, replacing isotopy by 1-isotopy.

\smallskip
\begin{Figure}[htb]{movesRaux/fig}{}{}
\centerline{\fig{}{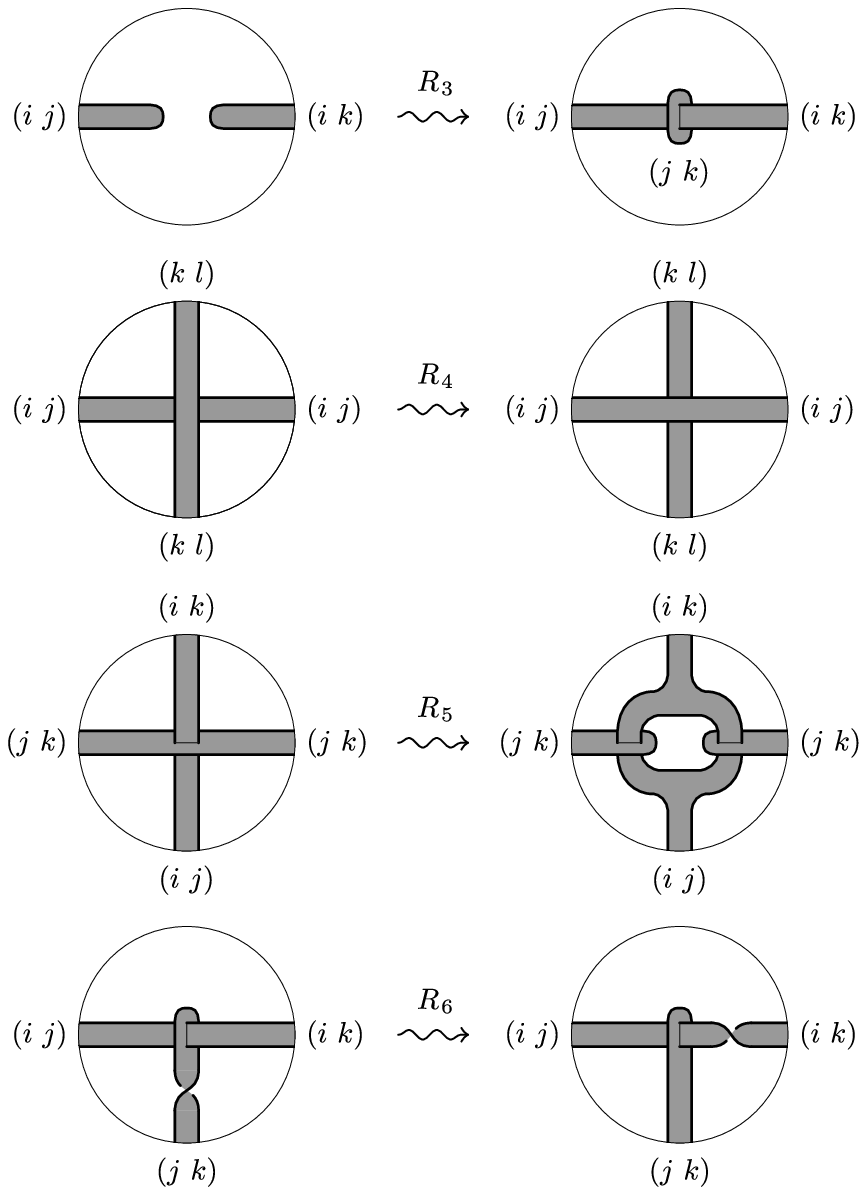}}
\end{Figure}

\begin{proposition} \label{movesInv/thm}
Moves $R_1$ and $R_2$ independently generate their own inverses, up to
labelled 1-isotopy.
\end{proposition}

\begin{proof}
For move $R_1$, the equation $R_1^{-1} = R_1^2$ obtained in the Introduction, by
thinking $R_1$ as a rotation of $120^\circ$, holds also in the present context,
since actually no isotopy is needed. On the other hand, move $R_2^{-1}$ is
equivalent, up to labelled diagram isotopy, to a suitable sequence of three moves of
types $I_2$, $I_3$ and $R_2$ in the order.
\end{proof}

\begin{proposition} \label{movesAux/thm}
Moves $R_1$ and $R_2$ generate moves $R_3$, $R_4$, $R_5$ and $R_6$, as well as their
inverses, up to labelled 1-isotopy.
\end{proposition}

\begin{proof} By Proposition \ref{movesInv/thm} we do not need to worry about
inverses. Move $R_4$ can be easily obtained as the composition of one move
$R_2^{-1}$ and one move $R_2$. Figures \ref{moveR3pr/fig}, \ref{moveR5pr/fig} and
\ref{moveR6pr/fig}              
\begin{Figure}[htb]{moveR3pr/fig}{}{}
\centerline{\fig{}{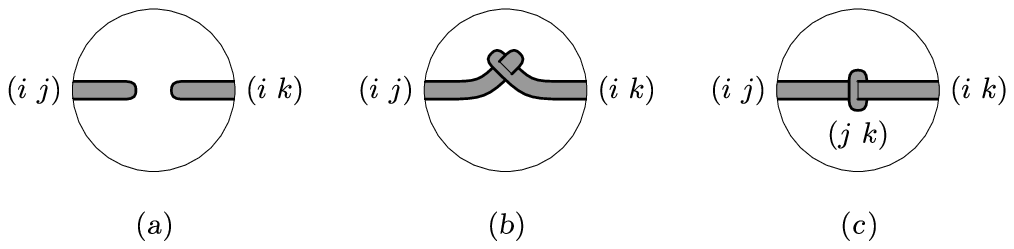}}
\end{Figure}%
\begin{Figure}[htb]{moveR5pr/fig}{}{}
\centerline{\fig{}{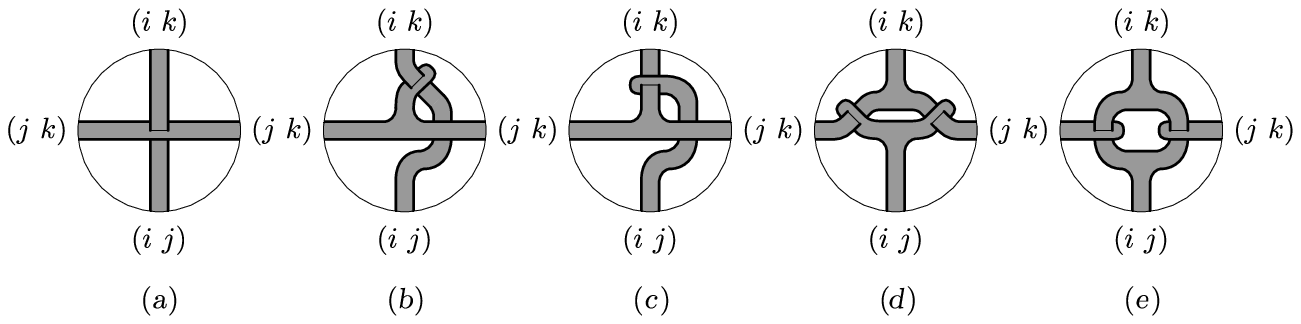}}
\end{Figure}%
\begin{Figure}[htb]{moveR6pr/fig}{}{}
\centerline{\fig{}{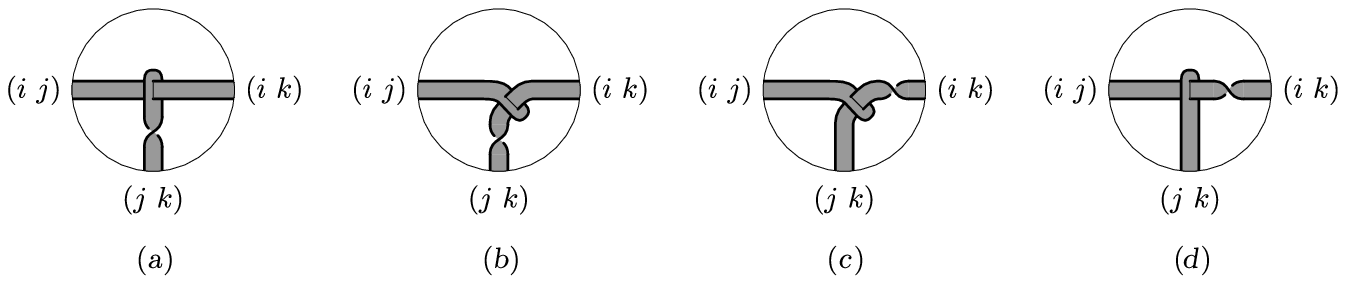}}
\end{Figure}
respectively shows how to get moves $R_3$,
$R_5$ and $R_6$ in terms of labelled 1-isotopy and moves $R_1$. In Figure
\ref{moveR3pr/fig}, we pass from \(a) to \(b) by one move $I_2$ and from \(b) to
\(c) by one move $R_1$. In Figure \ref{moveR5pr/fig}, \(b) is equivalent to \(a) up
to labelled diagram isotopy, then we perform respectively one move $R_1$, one move
$I_3$ and one pair of moves $R_1$ and $R_1^{-1}$ to obtain in the order \(c), \(d)
and \(e). In Figure \ref{moveR6pr/fig} we see that, up to conjugation by move $R_1$,
the twist transfer of move $R_6$ can be realized by the labelled diagram isotopy
between \(b) and \(c).
\end{proof}

\begin{remark} \label{orient/rem}
By labelled 1-isotopy and moves $R_1$ and $R_2$, any labelled ribbon surface
representing a connected simple branched covering of $B^4$ can be made orientable,
without changing the 2-equivalence class of the covering 4-dimensional
2-handlebody. In fact, twist transfer allows us to eliminate non-orientable
bands as shown in Figure \ref{movesRori/fig}. Here, assuming $(i\ j)$ and $(j\ k)$
distinct, we pass from \(a) to \(b) by two moves of types $I_2$ and $I_3$, and from
\(b) to \(c) by one move $R_6^{-1}$.
\end{remark}

\begin{Figure}[htb]{movesRori/fig}{}{}
\centerline{\fig{}{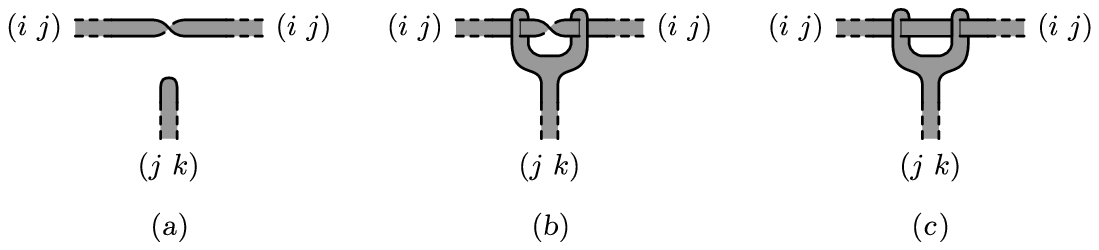}}
\end{Figure}

\section{From Kirby diagrams to labelled ribbon surfaces\label{toribbon/sec}}

In this section we prove the surjectivity of the map defined in the previous one,
which associates to each simply labelled ribbon surface $F$ the 2-equivalence class
of the generalized Kirby diagrams $K_F$. Since everything can be done
componentwise and any generalized Kirby diagram of a connected 4-dimensional
2-handlebody is 2-equivalent to an ordinary one (cf. Section \ref{prelim/sec}), we
will focus on ordinary Kirby diagrams and we will come back to the general case in
the last Proposition \ref{surjectivity/thm}.

Namely, for any ordinary Kirby diagram $K$, we construct a labelled ribbon surface
$F_K$ which represents the 2-equivalence class of the corresponding 4-dimensional
2-handlebody as a 3-fold simple branched covering of $B^4$ (cf. Proposition
\ref{diag/ribbon/diag/thm}). Such a construction is canonical in the sense that
the 4-stabilization of $F_K$ is uniquely determined up to labelled 1-isotopy and
covering moves $R_1$ and $R_2$. In this sense, $F_K$ actually depends only on the
2-equivalence class of $K$. 

In the light of the results of Section \ref{equivalence/sec}, we will relax the
vertical triviality condition for the link $L'$ in step \(c) of our construction to
just triviality (see Remark~\ref{trivialstate/rem}).
However, we temporarily impose the vertical triviality condition, in order to be
able to prove that $F_K$ does not depend (in the above sense) on the choice of $L'$,
without resorting to the results of the next section.

\medskip

Given an ordinary Kirby diagram $K$ describing a 4-dimensional 2-handlebody $H^0 \cup
H^1_1 \cup \dots \cup H^1_m \cup H^2_1 \cup \dots \cup H^2_n$, we let $C_1,
\dots, C_m$ be the disjoint flat disks\break spanned by the unknots corresponding to
the 1-handles and $L_1, \dots, L_n$ be the framed loops corresponding to the
2-handles. Moreover, we put $L = L_1 \cup \dots \cup L_n$ and think of it
indistinctly as a link or as a link diagram. Then, the construction of the labelled
ribbon surface $F_K$ is accomplished by the following steps:
\begin{itemize}\itemsep3pt

\item[\(a)]\vskip-\lastskip\smallskip
isotope $K$ into a standard form (cf. Figure \ref{diagst1/fig});

\item[\(b)] add to $K$ two standard disks $A_0,B_0 \subset R^3$ as shown in Figure
\ref{cover1/fig}; $A_0$ and the bottom part of $B_0$ (the one parallel to $A_0$ in
the diagram) are supposed to lie in a horizontal plane such that $K$ is entirely
contained in the half space above it; the height of the rest of $B_0$ and the $C_i$'s
varies according to the height of the arcs of $L$ passing through them (cf. step \(c)
below and Figure \ref{disks2/fig});

\begin{Figure}[htb]{cover1/fig}{}{}
\centerline{\fig{}{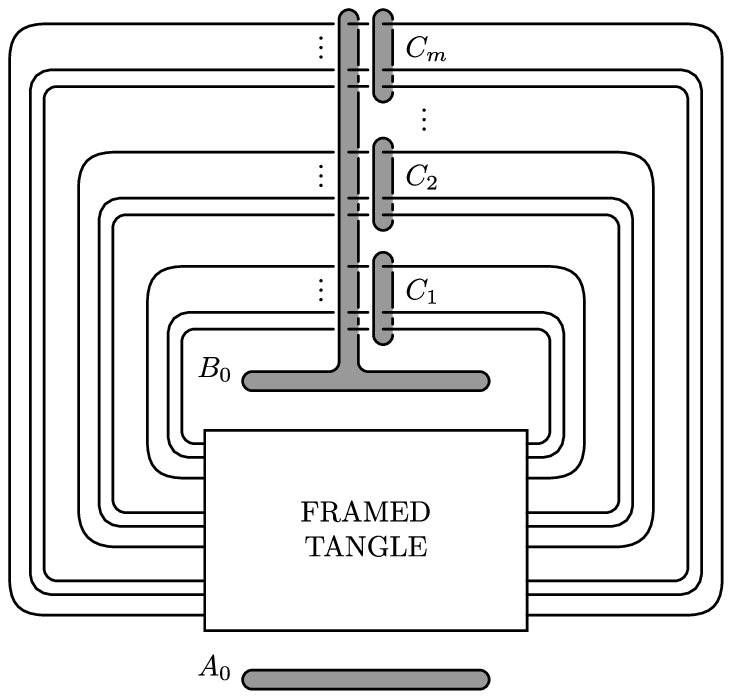}}
\end{Figure}

\begin{Figure}[b]{crossing1/fig}{}{}
\centerline{\fig{}{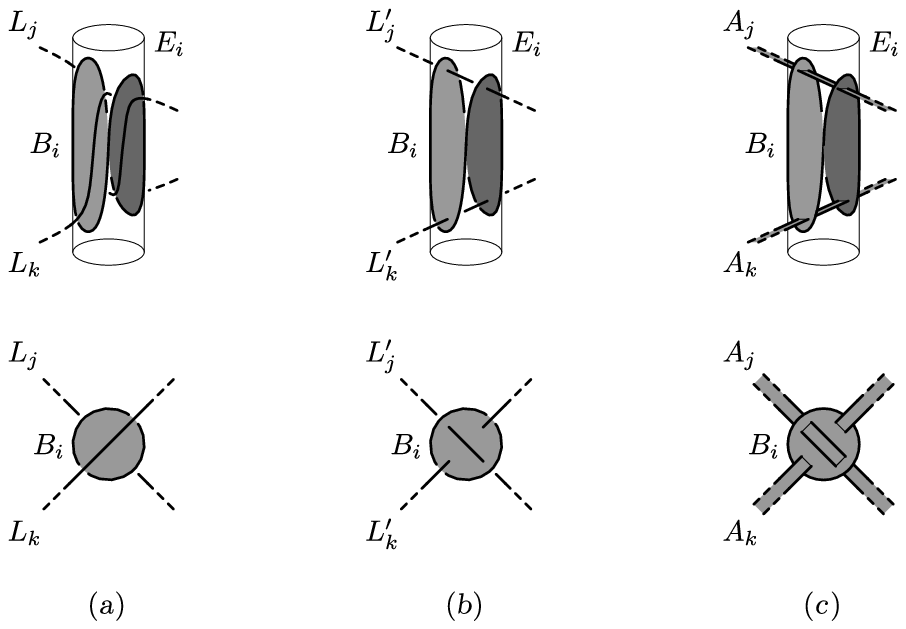}}
\end{Figure}

\item[\(c)] \label{stepc} choose a vertically trivial state $L'$ of $L$ (cf. Remark
\ref{trivialstate/rem}) and call $L'_i$ the component of $L'$ corresponding to $L_i$;
we think of $L'$ as a vertically trivial link which coincides with $L$ outside $E_1
\cup \dots \cup E_l$, where each $E_i$ is a cylinder projecting\break onto a small
circular neighborhood of a changing crossing inside the tangle box of Figure
\ref{cover1/fig} (of course, this is possible only after having suitably vertically
isotoped $L$); a cylinder $E_i$, together with the relative portion of diagram, is
depicted in Figure \ref{crossing1/fig} \(a) and \(b), where $j$ and $k$ may or may
not be distinct; here $B_i \subset E_i$ is a regularly embedded disk without vertical
tangencies, separating the two arcs of $L \cap E_i$ and forming four transversal
intersection with $L'$; furthermore, we assume that the height function is strictly
monotone and nearly constant on each arc of $L'$ outside the tangle box and that it
separates any two of such arcs, that is different arcs have disjoint height
intervals;

\item[\(d)] consider disjoint (possibly non orientable) narrow ribbons $A_1, \dots,
A_n \subset R^3$, such that each $A_i$ has $L'_i$ as its core and is obtained by a
regular vertical homotopy from a ribbon representing half the framing of $L_i$ plus
one positive (resp. negative) full twist for each positive (resp. negative)
crossing of $L_i$ inverted to get $L'_i$; each $A_i$ is assumed to be disjoint from
$A_0$ and to form with the $B_i$'s and the $C_j$'s only the ribbon intersections
shown in Figures \ref{crossing1/fig} \(c) and \ref{cover2/fig};

\begin{Figure}[htb]{cover2/fig}{}{}
\centerline{\fig{}{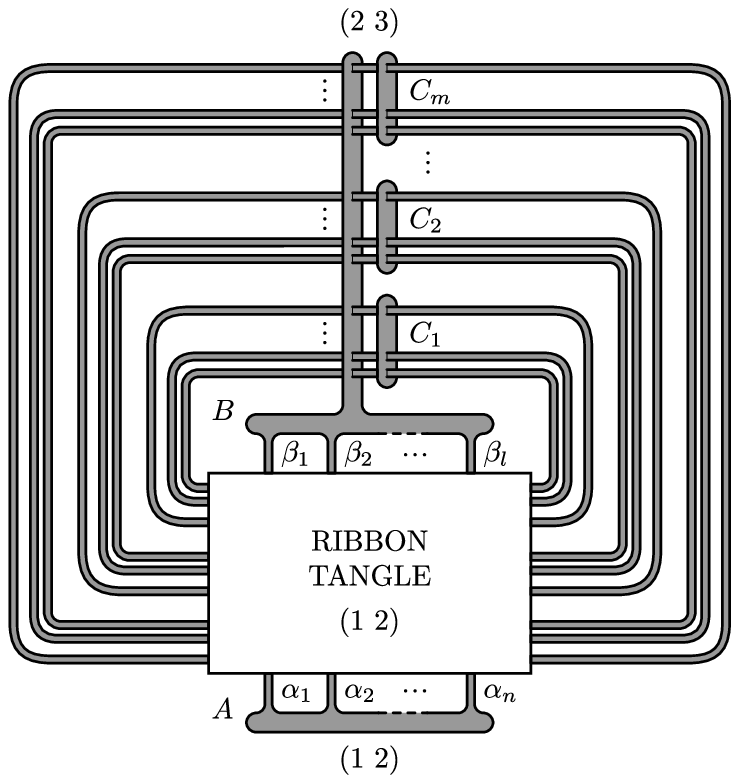}}
\end{Figure}

\item[\(e)] attach to $A_0 \cup A_1 \cup \dots \cup A_n$ disjoint narrow bands
$\alpha_1, \dots, \alpha_n \subset R^3$ to get a connected non singular surface $A$;
each $\alpha_i$ is an embedded 1-handle between $A_0$ and $A_i$, which is assumed
to be disjoint from the $E_j$'s and to have standard projection outside the tangle
box as in Figure \ref{cover2/fig}; we constrain the intersection of the $\alpha_i$'s
with the tangle box to assume height values disjoint from the ones of the link $L'$;
more precisely, if $[a_i,b_i]$ is the height interval of $L'_i$, then inside the
tangle box $\alpha_i$ takes height values just below $a_i$; by the proof of
Proposition \ref{diag/ribbon/thm}, this last assumption is clearly much more that we
really need, we nevertheless make it in order to keep things simpler;

\item[\(f)] attach to $B_0 \cup B_1 \cup \dots \cup B_l$ disjoint narrow bands
$\beta_1, \dots, \beta_l \subset R^3$ to get a connected non singular surface $B$;
each $\beta_i$ is an embedded 1-handle between $B_0$ and $B_i$, which is assumed to
be disjoint from $A$ and from the interiors of the $E_j$'s and to have standard
projection outside the tangle box as in Figure \ref{cover2/fig}; 

\item[\(g)] define $F_K \subset B^4$ to be the ribbon surface having $A \cup B
\cup C_1 \cup \dots \cup C_m$ as 3-dimensional diagram, with the unique labelling
assigning to $A_0$ the transposition $(1\ 2)$ and to $B_0, C_1, \dots, C_m$ the
transposition $(2\ 3)$ (cf. Figure \ref{cover2/fig}); in particular, Figure
\ref{crossing2/fig} shows such labelling in a neighborhood of $B_i$ (here, as in
Figure \ref{crossing1/fig},\break $j$ and $k$ may or may not be distict).

\end{itemize}\vskip-\lastskip\smallskip

\begin{Figure}[htb]{crossing2/fig}{}{}
\centerline{\fig{}{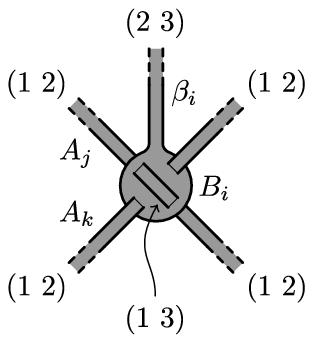}}
\end{Figure}

It is worth noting that, by Remark \ref{orient/rem} one more step could be added in
order to make the ribbon surface $F_K$ orientable and even more to make its diagram
blackboard parallel. However, we omit such additional step, since we will not need
those properties in what follows.


\begin{proposition} \label{diag/ribbon/thm}
Let $K$ an ordinary Kirby diagram. Then the 4-stabilization of the labelled ribbon
surface $F_K$ constructed above is uniquely determined by $K$, that is it does not 
depend on the choices involved in the construction, up to labelled 1-isotopy and
moves $R_1$ and $R_2$.
\end{proposition}

\begin{proof}
First of all, as a preliminary, we add some extra structure to the above
construction of $F_K$. Namely, we consider disjoint disks $D_1, \dots, D_n \subset
R^3$ respectively spanned by $L'_1, \dots, L'_n$, such that the intersection of $D_1
\cup \dots \cup D_n$ with any hori\-zontal plane is either empty, one point or one
arc. We constrain all these arcs to project onto the line segments drawn in 
\begin{Figure}[htb]{disks1/fig}{}{}
\vskip6pt\centerline{\fig{}{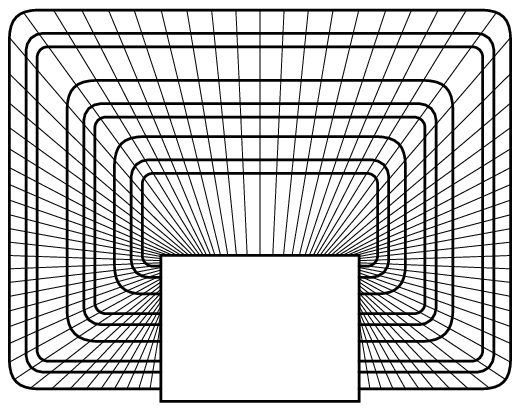}}
\end{Figure}%
Figure \ref{disks1/fig} outside the tangle box of Figure \ref{cover1/fig}.\break
Then, the portion of the $D_i$'s outside the tangle box turns out to be nearly
horizontal and completely determined by $L'$. Moveover, a suitable choice of the
height function of $B_0$ and of the $C_j$'s, ensures that the $D_i$'s form with $B_0
\cup C_1 \cup \dots \cup C_m$ only clasps\break and ribbon intersections, like the
ones depicted in Figure \ref{disks2/fig}, where the 3-dimensional view on the left
side is compared with the corresponding diagram on the right side. Here, the clasps
arise from the intersections of the $L_i$'s with $B_0 \cup C_1 \cup \dots \cup C_m$
and the ribbon intersections are due to the height variations of $B_0$ and $C_j$'s
needed to let the $L'_i$'s pass through them.

\begin{Figure}[htb]{disks2/fig}{}{}
\centerline{\fig{}{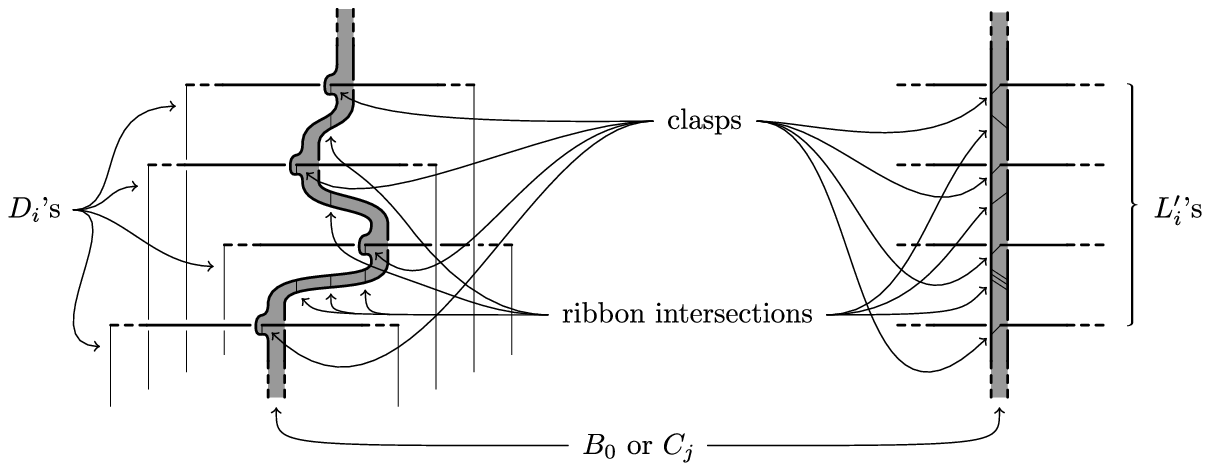}}
\end{Figure}

Without loss of generality, up to small perturbations, we can require that also inside
the tangle box the $D_i$'s form only clasps and ribbon intersections with $F_K$.
Namely, we assume that: 1) $A_j \cap D_j$ consists of a certain number of disjoint
clasps connecting $L'_j$ with the boundary of $A_j$, in such a way that $D_j \cup A_j$
is collapsible; 2) each $B_j$ forms with the $D_i$'s four clasps and some (possibly
none) ribbon intersections, as shown in Figure \ref{disks3/fig}; 3) the $\beta_j$'s may
pass through the $D_i$'s, forming ribbon intersections with them. Finally, we observe
that, by construction, each $D_i$ is disjoint from the $A_j$'s with $j \neq i$ and from
all the $\alpha_j$'s.

\begin{Figure}[htb]{disks3/fig}{}{}
\centerline{\fig{}{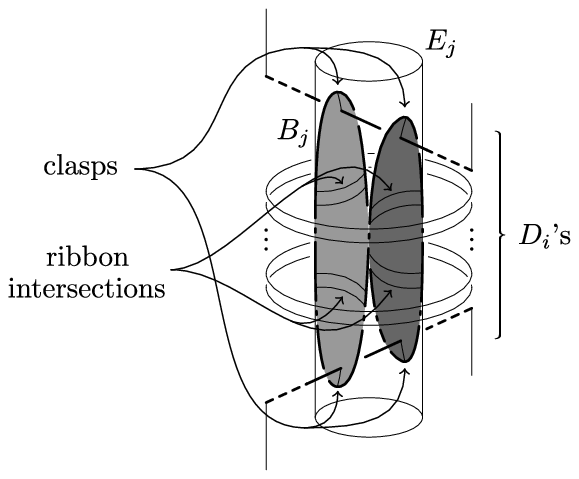}}
\end{Figure}

At this point we pass to the core of the proof. Given an ordinary Kirby diagram $K$,
the relevant choices occurring in the construction of $F_K$ are, in order:\break
1)~the standard form $K$; 2)~the vertically trivial state $L'$; 3)~the bands
$\alpha_1, \dots, \alpha_n$; 4)~the bands $\beta_1, \dots, \beta_l$. In fact, the
$A_i$'s and the $B_i$'s are uniquely determined up to diagram isotopy.

We prove that the 4-stabilization of $F_K$ is independent on these choices, up to
labelled 1-isotopy and moves $R_1$ and $R_2$, by proceeding in the reverse order and
assuming each time that all previous choices have been fixed. By Propositions
\ref{movesInv/thm} and \ref{movesAux/thm}, in addition to moves $R_1$ and $R_2$, we
can use also the moves $R_3$, $R_4$, $R_5$ and $R_6$ introduced in the previous
section, as well as the inverses of all such moves.

Concerning the $\alpha_i$'s and the $\beta_i$'s, it suffices to prove that, in
presence of a stabilizing disk, labelled 1-isotopy and the moves above enable us to
change them one by one.

Figure \ref{indepBeta/fig} shows how to deal with the band $\beta_i$ of Figure
\ref{crossing2/fig}. The small disk with label $(3\ 4)$ in \(a) is the stabilizing
disk. This can be moved to form one ribbon intersection with $\beta_i$ as in \(d),
by labelled 1-isotopy and four moves $R_2$. Parts \(b) and \(c) of the figure
represent 1-isotopic intermediate steps. Looking at the diagram, we can realize such
a modification by an isotopy $H:B^2 \times [0,1] \to R^3$ between the two disks
labelled $(3\ 4)$ in \(a) and \(d), which is a suitable homeomorphism of $B^2 \times
[0,1]$ onto a regular neighborhood of $B_i$ whose boundary forms four ribbon
intersections with $A_i$ and $A_j$. Finally, we perform one move $R_3^{-1}$ on \(d)
to cut the band $\beta_i$.\break The result is clearly independent of $\beta_i$ up to
diagram isotopy, so we are done.

\begin{Figure}[htb]{indepBeta/fig}{}{}
\centerline{\fig{}{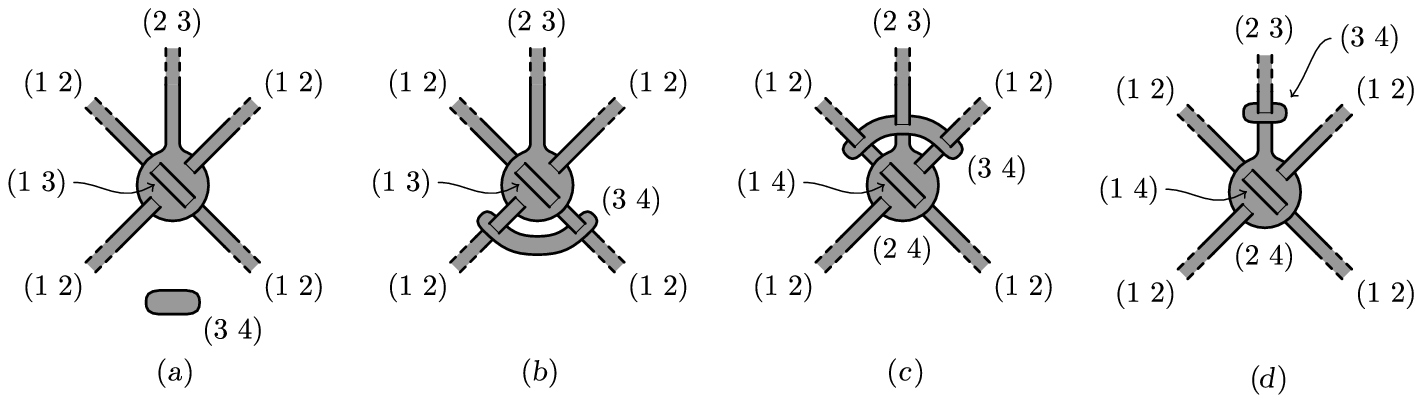}}
\end{Figure}

\begin{Figure}[b]{indepAlpha/fig}{}{}
\centerline{\fig{}{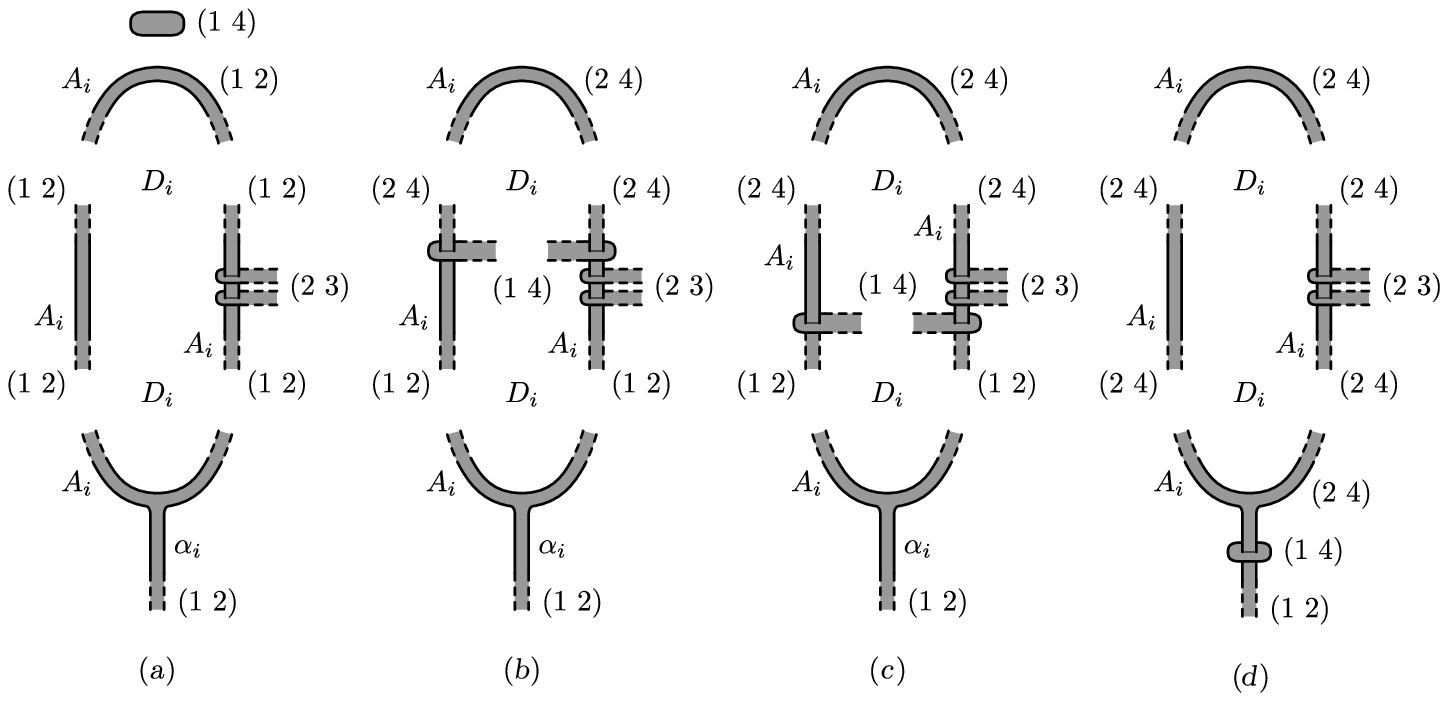}}
\end{Figure}

The same idea also applies to the band $\alpha_i$, as suggested by Figure
\ref{indepAlpha/fig}. We warn the reader that Figure \ref{indepAlpha/fig} is much
more sketchy than Figure \ref{indepBeta/fig}. In fact, in place of the disk $B_i$ we
have here the complex $D_i \cup A_i$, which can be large and complicated, although
still collapsible. Before of starting the process, we let the stabilizing disk of
Figure \ref{indepBeta/fig} \(a) pass first through $B_0$ and then through $A_0$, in
such a way that its label becomes $(1\ 4)$, as in Figure \ref{indepAlpha/fig} \(a).
This time, the need for move $R_2$ is due to the ribbon intersections that $A_i$ may
form passing through the $B_j$'s and the $C_j$'s and that the $B_j$'s (including
$B_0$), the $C_j$'s and the $\beta_j$'s may form passing through the interior
of $D_i$ (cf. Figures \ref{disks2/fig} and \ref{disks3/fig}). In particular, the
ribbon intersections along $A_i$ always appear in pairs, each pair being formed with
$B_0$ and one of the $C_j$'s (cf. Figure \ref{cover2/fig}) or with one of the $B_j$'s
(cf. Figure \ref{crossing1/fig}). Any such pair looks like the one pictured in
Figure \ref{indepAlpha/fig} \(a). Comparing steps \(b) and \(c), one see how the
stabilizing disk can be pushed beyond this pair, by using labelled 1-isotopy and
moves $R_2$.\break On the other hand, only one move $R_4$ suffices to go beyond each
one of the ribbon intersections in the interior of $D_i$. Eventually, the
stabilizing disk reach the position of step \(d), so we can conclude as above by
cutting the band $\alpha_i$.

Now we pass on to the vertically trivial state $L'$. Recall that we are thinking
of it as a vertically trivial link, that is a vertically trivial diagram together
with a compatible height funtion. Of course, course different choices of the height
function compatible with the same diagram are related by a vertical diagram
isotopy. Such an ambient isotopy can be used to relate the entire resulting
surfaces except for the bands $\alpha_1, \dots, \alpha_n$. In fact, each band
$\alpha_i$ is forced, by the vertical constraints we imposed in its definition (cf.
step \(e) of the construction of $F_K$), to be attached to the corresponding
$A_i$ near to the point of minimum height. However, this problem can be overcome by
the above proof of independence on the $\alpha_i$'s. Actually, this delicate point
is the only obstruction to the naive solution of the problem of the independence on
the $\alpha_i$'s consisting in fixing a standard form for them. 

Thus, having settled the problem of the height function, we are left with the
modifications of Proposition \ref{vertstat/thm} on the diagram. First we address the
change of ordering of the link components. Of course, it is enough to deal with the
transposition in the vertical order of any two components $L'_i$ and $L'_j$, that
means to simultaneously change all the crossing involving both $L'_i$ and $L'_j$.
Without loss of generality we assume that $L'_i$ lies under $L'_j$.
We begin as above, operating with the stabilizing disk around the $D_i$ to get the
configuration of Figure \ref{indepAlpha/fig} \(d). As a result, we have a global
labelling change on the ribbon $A_i$ from $(1\ 2)$ to $(2\ 4)$, while the labelling
of $A_j$ is left unchanged. Then, we can perform the crossing changes as described
in Figure \ref{indepOrder/fig}. Here, apart from 1-isotopy, we have only one move
$R_4$ relating \(b) and \(c). To be precise, Figure \ref{indepOrder/fig} covers
only one of the two possible cases, the other one being covered by the same steps in
the reverse order with the roles of $A_i$ and $A_j$ exchanged. After all the crossing
changes have been performed, we bring the stabilizing disk back to the original
position, by reversing the process of Figure \ref{indepAlpha/fig}.

\begin{Figure}[htb]{indepOrder/fig}{}{}
\medskip\centerline{\fig{}{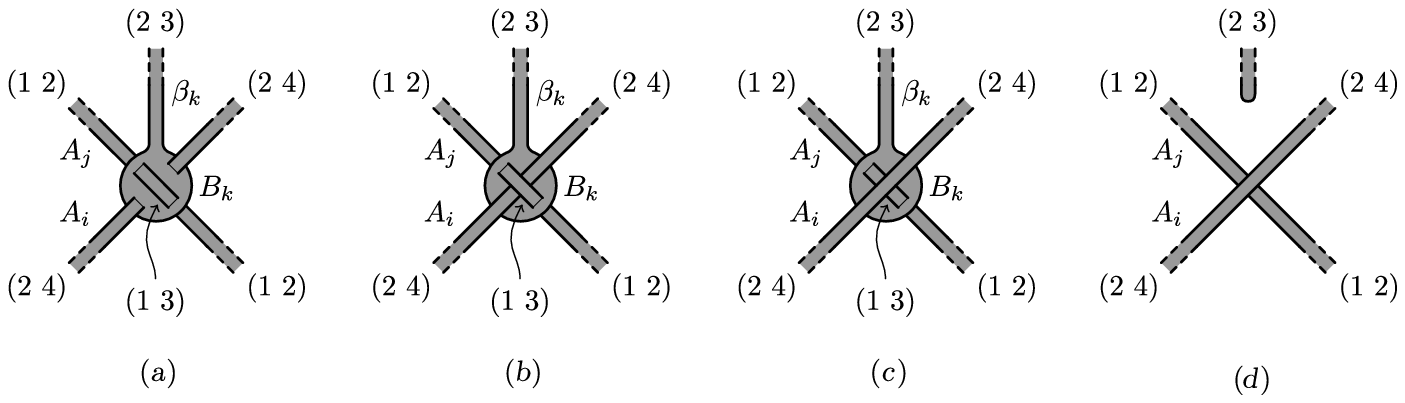}}
\end{Figure}

Concerning single crossing changes making a vertically trivial component $L'_i$
into a different vertically trivial state of $L_i$, there are four cases to be
considered, depending on sign of the crossing and on whether $L'_i$ coincides with
$L_i$ at that crossing or not.
In all cases, it is not restrictive to assume that the two points of $L'_i$
projecting to the crossing are distinct from the unique minimum height point $p_i$
of $L'_i$ and that the vertical segment joining them is contained in $D_i$. Such a
segment divide $D_i$ into two disks. We call $D'_i$ the one which does not contain
$p_i$ and assume that the attaching arc of the band $\alpha_i$ to $A_i$ is disjoint
from the part of $A_i$ running along the boundary of $D'_i$.

Figure \ref{indepState/fig} indicates how to realize the crossing change in one of
the four cases. For the other three cases it suffices to apply a mirror symmetry
to all the stages and/or reverse their order.
First, we pass from \(a) to \(b) by the same the process described in Figure
\ref{indepAlpha/fig}, with the only difference that here we have $D'_i$ in place
of $D_i$. We continue that process one more step to get \(c), by pushing the
stabilizing disk beyond one of the two pairs of ribbon intersections between $A_i$
and $B_j$. Then, we obtain \(d) from \(c) in the same way we obtained \(c) from \(a)
in Figure \ref{indepOrder/fig}. Finally, we push back the stabilizing disk through
the same ribbon intersections as above to achieve \(e). At this point, we bring back
the stabilizing disk in the original position, by reversing the process from \(a) to
\(b), after having transferred the full twist present on it to $A_i$ by move $R_6$.
We remind the reader that the additional full twist on $A_i$ compensates the change
of crossing (cf. definition of $A_i$ in step \(d) of the construction of
$F_K$).

\begin{Figure}[htb]{indepState/fig}{}{}
\centerline{\fig{}{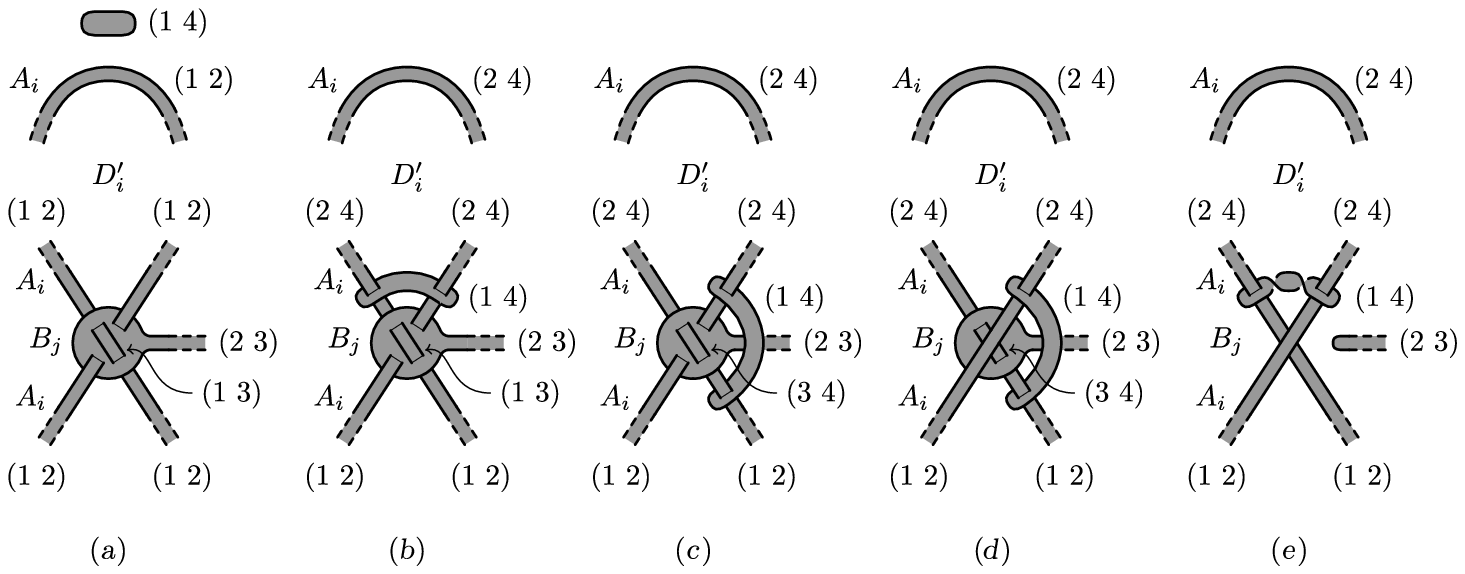}}
\end{Figure}

It remains to deal with the standard form of $K$. As discussed in Section
\ref{prelim/sec}, different choices of such standard form are related by a finite
sequence of Reidemaister moves inside the tangle box and moves of the types depicted
in Figures \ref{diagst2/fig} and \ref{diagst3/fig}. Hence, we have to prove the
invariance of the 4-stabilization of $F_K$ with respect to those moves.

We observe that any Reidemeister move on $L$ induces the same move on $L'$ and so
just a diagram isotopy on $F_K$, provided that none of the involved crossings
(before as well as after the move) has been changed when passing from $L$ to $L'$.
The reason is that in this case the two links coincide inside a small 3-cell where
the move takes place and such 3-cell is free from the $B_i$'s. We leave to the reader
the straightforward verification that a vertically trivial state $L'$ of $L$ with
the required property can be always achieved by a suitable application of the naive
unknotting procedure described in Section \ref{prelim/sec} (with height function
on each component as in Figure \ref{height/fig} \(a) or (c), depending on the move).

On the other hand, by relaxing a little bit the constraints for the standard form of
$K$ outside the tangle box, the moves of Figures \ref{diagst2/fig} and
\ref{diagst3/fig} reduce to the ones given in Figure \ref{indepDiag1/fig}.

\begin{Figure}[htb]{indepDiag1/fig}{}{}
\vskip-6pt\centerline{\fig{}{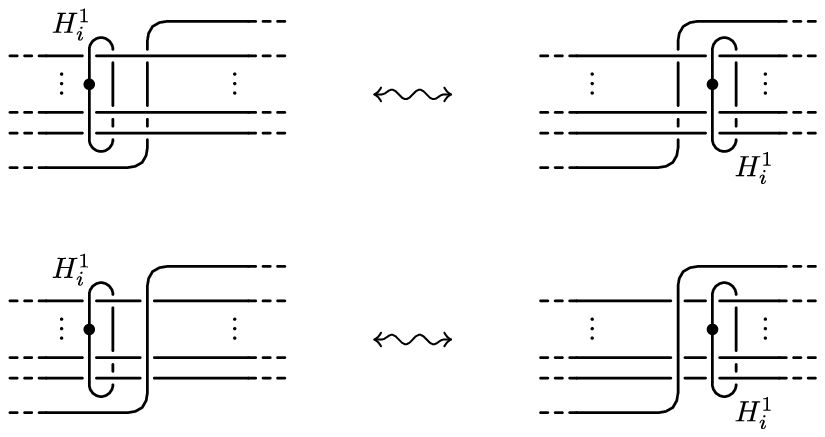}}
\end{Figure}

\begin{Figure}[b]{indepDiag2/fig}{}{}
\centerline{\fig{}{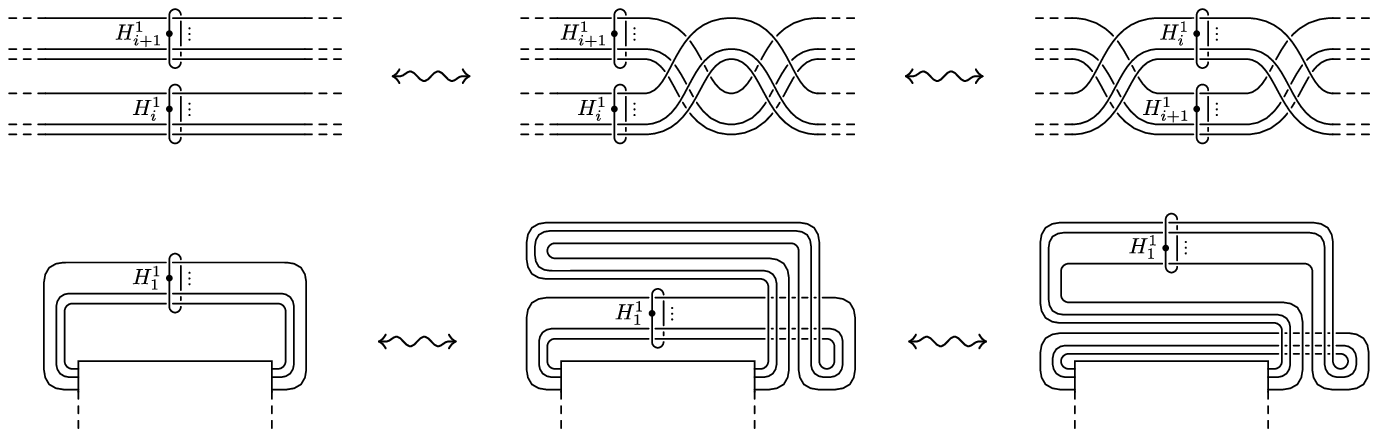}}
\end{Figure}

\begin{Figure}[b]{indepDiag3/fig}{}{}
\vskip9pt\centerline{\fig{}{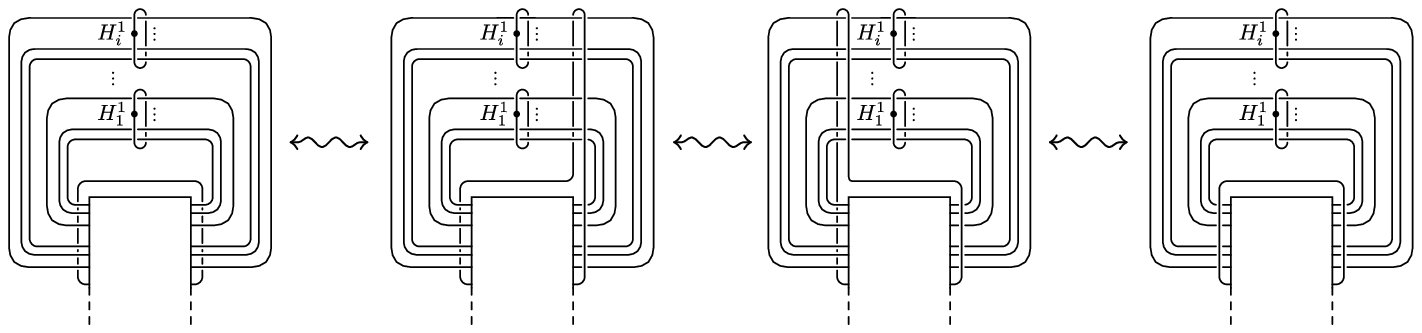}}
\end{Figure}

To be more precise, let us look at the following Figures \ref{indepDiag2/fig} and
\ref{indepDiag3/fig}, which regard respectively the moves of Figures
\ref{diagst2/fig} and \ref{diagst3/fig}. Here, only the relevant part of the
diagram is drawn, however it should be clear what is happening. We notice that, in
addition to Reidemeister moves inside the tangle box, we just need to bring some
crossings out of the tangle box by isotopy of the planar link diagram and to move
certain dotted components through them by a sequence of moves as in Figure
\ref{indepDiag1/fig} (for the sake of simplicity, the figures are drawn as if the
Reidemeister moves were performed directly outside the tangle box).

To conclude the proof, we have only to interpret the moves of Figure
\ref{indepDiag1/fig} on the diagram $K$ in terms of moves $R_i$ on the corresponding
surface $F_K$. In fact, we already know how to handle Reidemeister moves, while
isotopy of the planar link diagram of $K$ obviously induces diagram isotopy on
$F_K$.

We first operate on the disks $C_i, \dots, C_{i+h}$ corresponding to the 1-handles
we want to move, as described in Figure \ref{indepDiag4/fig}. Namely, apart from the
diagram isotopy relating \(a) and \(b), we perform a certain number of moves $R_2$
to let the stabilizing disk labelled by $(3\ 4)$ reach the position of \(c) and then
we create the new component $B_0'$ with monodromy $(2\ 4)$ that appears in \(d) by
one move $R_3$. We emphasize that also the monodromy of the disks $C_i, \dots,
C_{i+h}$ is now $(2\ 4)$. 
Actually, we could limit ourselves to one single disk at a time, but for our purpose
it is more convenient to operate simultaneously on all the disks involved by the
moves of Figures \ref{indepDiag2/fig} and \ref{indepDiag3/fig}.

\begin{Figure}[htb]{indepDiag4/fig}{}{}
\centerline{\fig{}{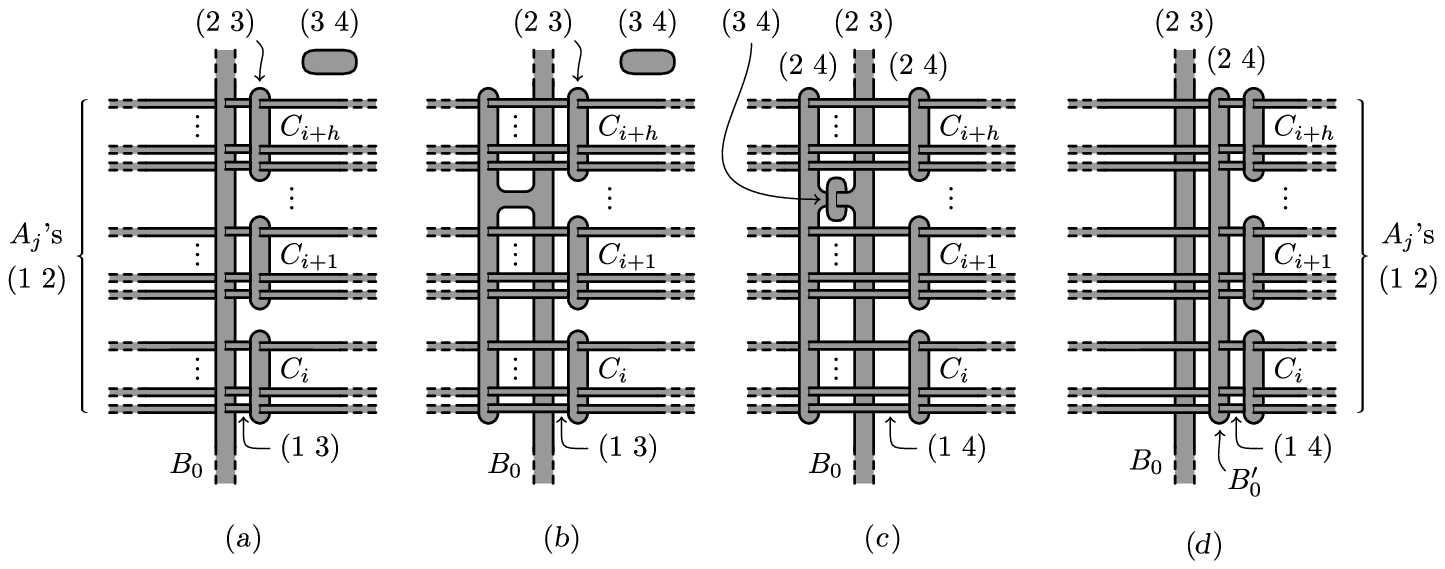}}
\end{Figure}

Then, we observe that the moves of Figure \ref{indepDiag1/fig} would correspond to a
diagram isotopy of $F_K$, if $L'$ coincided with $L$ at all the involved link
crossings. In fact, if this were the case it would be enough to slide $B_0'$ and
$C_i$ over/under all such crossings. So, we only have problems with the disks
$B_k$ at the crossings we changed to get the vertically trivial link $L'$. Figure
\ref{indepDiag5/fig} shows how to handle these crossings. Here, apart from
1-isotopy, we apply four moves $R_4$ to pass from \(b) to \(c).
\end{proof}

\begin{Figure}[htb]{indepDiag5/fig}{}{}
\vglue2mm\centerline{\fig{}{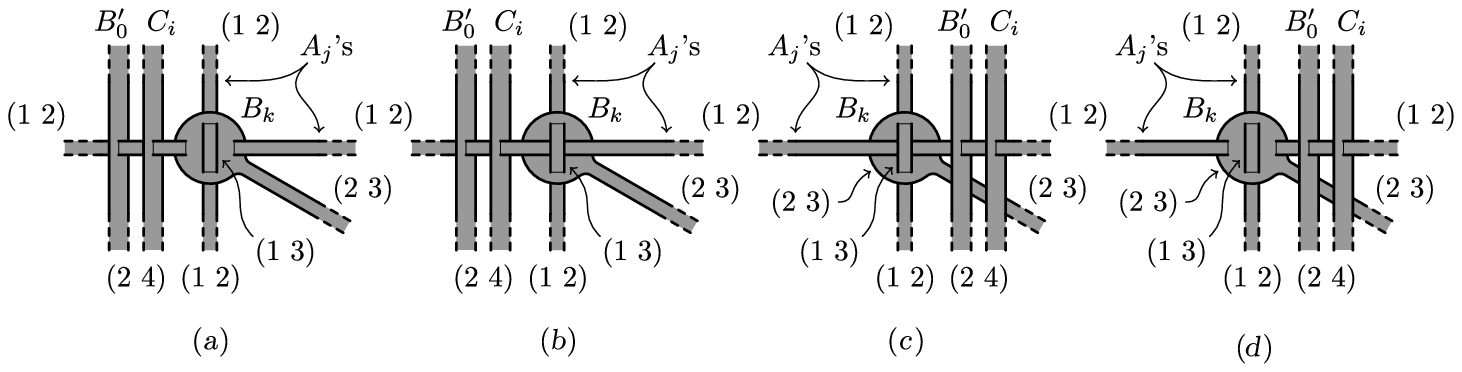}}
\end{Figure}

\begin{remark} \label{diag/ribbon/rem}
The last part of the proof of Proposition \ref{diag/ribbon/thm} suggests that the
above construction of the labelled ribbon surface $F_K$ could be easily adapted to
work directly with any ordinary Kirby diagram, avoiding the neeed for the standard
form. A little more effort would enable us to include also generalized Kirby
diagrams, of course by using labellings of arbitrary degree (cf. proof of
Proposition \ref{surjectivity/thm}).\break
The reason for focusing on the ordinary case is that we want to control the degree
of the converings, in view of the results of the next section. On the contrary, the
choice of working with the standard form is motivated only by the sake of
convenience.
\end{remark}

So far we have proved that the labelled ribbon surface $F_K$ associated to a Kirby
diagram $K$, is well defined up to 4-stabilization, labelled 1-isotopy and ribbon
moves. The next proposition addresses the invariance of $F_K$ under 2-deformations
of $K$.

\begin{proposition} \label{2def/moves/thm}
If $K$ and $K'$ are 2-equivalent ordinary Kirby diagrams, then the 4-stabilizations
of the labelled ribbon surfaces $F_K$ and $F_{K'}$ are equivalent up to labelled
1-isotopy and moves $R_1$ and $R_2$.
\end{proposition}

\begin{proof}
By Proposition \ref{diag/ribbon/thm} we can assume $K$ and $K'$ to be in standard
form and the 2-deformation between them to be described in terms of moves as in
Proposition \ref{stmoves/thm}. Moreover, we do not need to worry about handle
isotopy, since diagram isotopy has already been treated in the proof of Proposition
\ref{diag/ribbon/thm}, while the first move of Figure \ref{diag3/fig} performed on
ordinary Kirby diagrams in standard form trivially induces labelled isotopy on the
corresponding labelled ribbon surfaces.

Hence, it remains to consider the cases when $K$ and $K'$ differ by a pair of
cancelling 1/2-handles (cf. Figure \ref{diag4/fig}) and by a 2-handle sliding
(cf. Figure \ref{diag5/fig}).

The first case is quite easy. In fact, if the disk $C_i$ and the loop $L_j$ of $K$
represent two cancelling handles $H^1_i$ and $H^2_j$, then in $F_K$ only the ribbon
$A_j$ passes once through $C_i$. Therefore, by a move $R_3$ we can remove
$C_i$ and break $A_j$ into two long tounges. At this point, labelled 1-isotopy
allows us to completely retract such tounges and after that also the ones $\beta_k
\cup B_k$ related to crossings involving $L_j$.

\begin{Figure}[b]{sliding1/fig}{}{}
\centerline{\fig{}{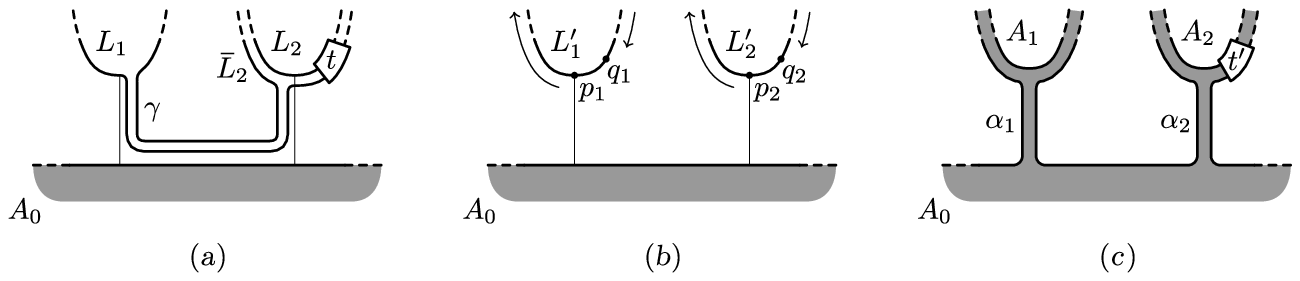}}
\end{Figure}

The case of a 2-handle sliding requires some preliminaries. First of all, let us
renumber the 2-handles starting from the two ones involved in the sliding, in such a
way that $H^2_1$ slides over $H^2_2$. In terms of Kirby diagram, this means to
replace $L_1$ with the band connected sum $L_1 \#_\gamma \bL_2$, where $\bL_2$ is a
parallel copy of $L_2$ realizing its framing, and $\gamma$ is a band connecting $L_1$
to $\bL_2$. Up to isotopy, $\gamma$ can be assumed to be a blackboard parallel band
which does not form any crossing with the $L_i$'s, as in Figure~\ref{sliding1/fig}
\(a). We also assume that the vertically trivial status $L'$ choosen to construct
$F_K$ satisfies the following properties: 1) the vertical order of the components is
the one given the numbering, that is $L'_i$ lies under $L'_j$ for any $i < j$;\break
2) the minimum point $p_1$ (resp. $p_2$) and the maximum point $q_1$ (resp. $q_2$)
of the height function on $L'_1$ (resp. $L'_2$) coincide with the end points of
(resp. are close to) the attaching arc of $\gamma$ to $L_1$ (resp. $\bL_2$), as
in Figure \ref{sliding1/fig} \(b). Here, the arrows indicate the orientations that
we will use in the framing computation at the end of the proof, so they are not
relevant for the moment. Finally, we choose $\alpha_1$ and $\alpha_2$ to be
blackboard parallel bands, such that $\gamma$ can be thought to run parallel to them
and to the part of $\Bd A_0$ between them, as in Figure \ref{sliding1/fig} \(c). For
the sake of convenience, the framing of $L_2$ and the ribbon $A_2$ are assumed to be
blackboard parallel outside the twist boxes $t$ and $t'$ respectively in Figure
\ref{sliding1/fig} \(a) and \(c). We warn the reader that the number of twists
inside such boxes is not the same in \(a) and \(c), accordingly to step \(d) of the
construction on $F_K$.

Once it has been set up in this way, the sliding can be interpreted in terms of
ribbon moves on the 4-stabilization of $F_K$ as sketched in Figure
\ref{sliding2/fig}. We think of $A_1$ as a 1-handle attached to $\Cl(A - A_1)$
and slide one of its attaching arcs along the boundary of $\Cl(A - A_1)$ as
indicated by the arrows in \(a) and \(b). Before of reaching the twist box $t$, this
sliding can be entirely realized by labelled diagram isotopy, except for the labelled
1-isotopy moves (of types $I_2$ and $I_3$) needed to pass through the disks $B_i$
encountered by $A_2$. Each time a disk $B_i$ is passed through, two new ribbon
intersections appear as shown in part \(a) of Figure \ref{sliding3/fig}. Then, we
use again 1-isotopy to split $B_i$ into two twin disks similar to the original one,
as suggested by the remaining parts of Figure \ref{sliding3/fig}.

\begin{Figure}[htb]{sliding2/fig}{}{}
\centerline{\fig{}{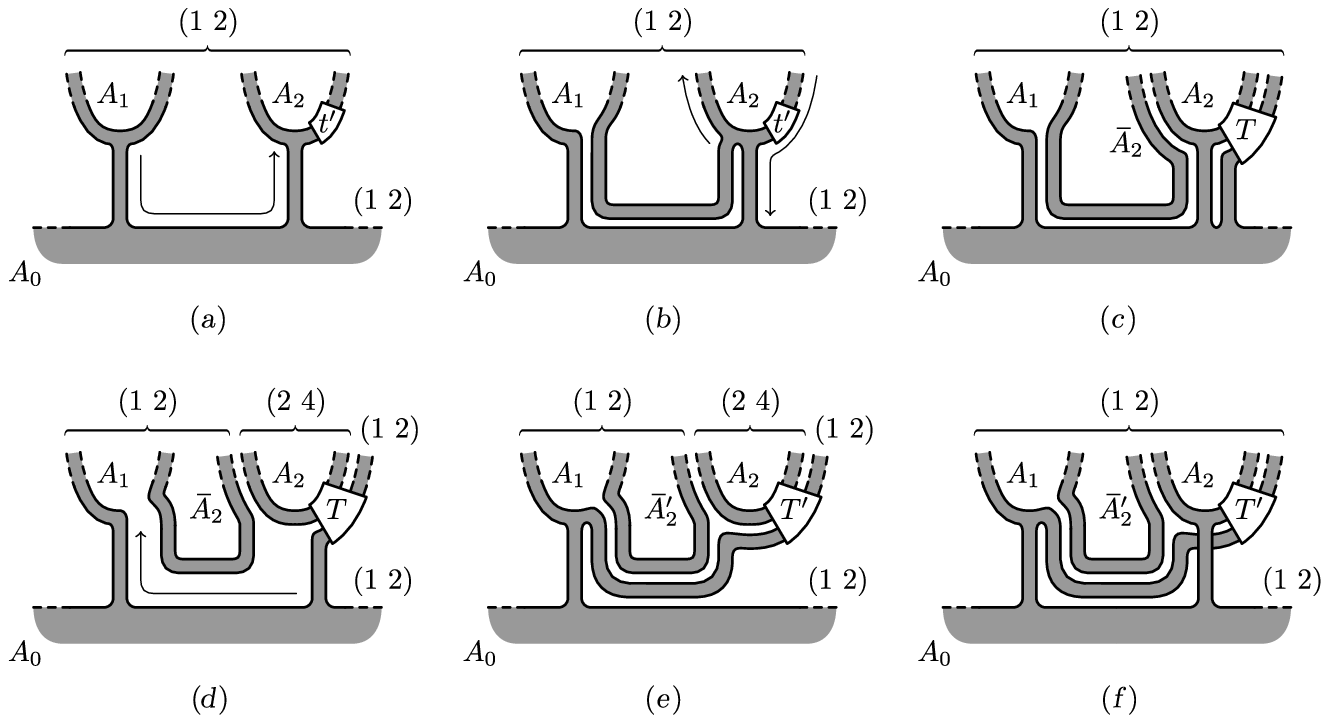}}
\end{Figure}

\begin{Figure}[htb]{sliding3/fig}{}{}
\vglue2mm\centerline{\fig{}{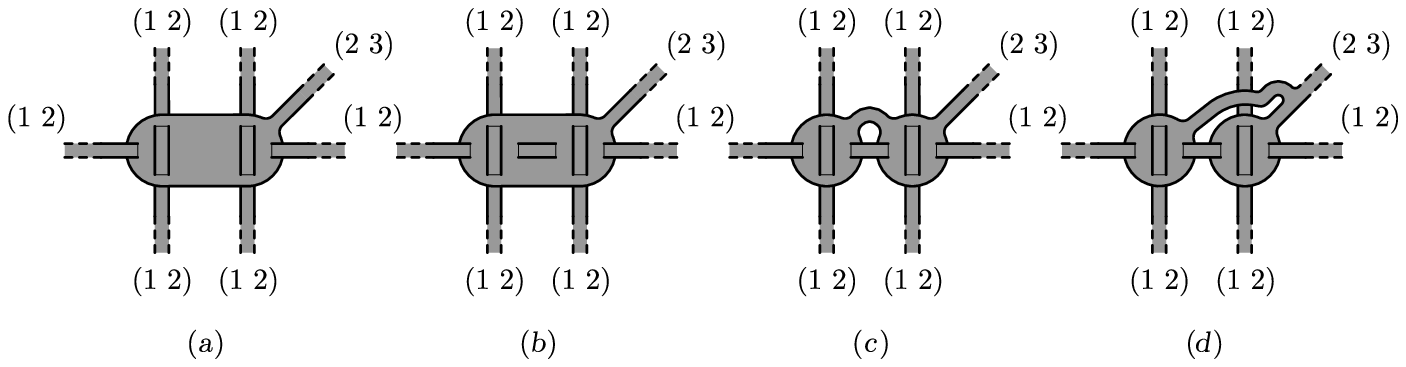}}
\end{Figure}

When traversing the twist box $t'$ in \(b) to get the twist box $T$ in \(c), after
having followed all the twists of $A_2$, we add some further crossings between $A_2$
and the parallel closed ribbon $\bA_2$ (together with further $B_i$'s), in order to
make them unlinked. Actually, $\bA_2$ itself is always well defined thanks to these
additional crossings, having their number the same parity of the number of half
twists in $t'$. Figure \ref{sliding4/fig} shows how to add a positive crossing; for
a negative one it suffices to mirror the figure. Here, some moves other than
1-isotopy are needed: one move $R_5$ from \(a) to \(b); two opposite moves $R_6$
(twist transfers) from \(c) to \(d); two moves $R_1$ from \(d) to \(e).

\begin{Figure}[htb]{sliding4/fig}{}{}
\medskip\centerline{\fig{}{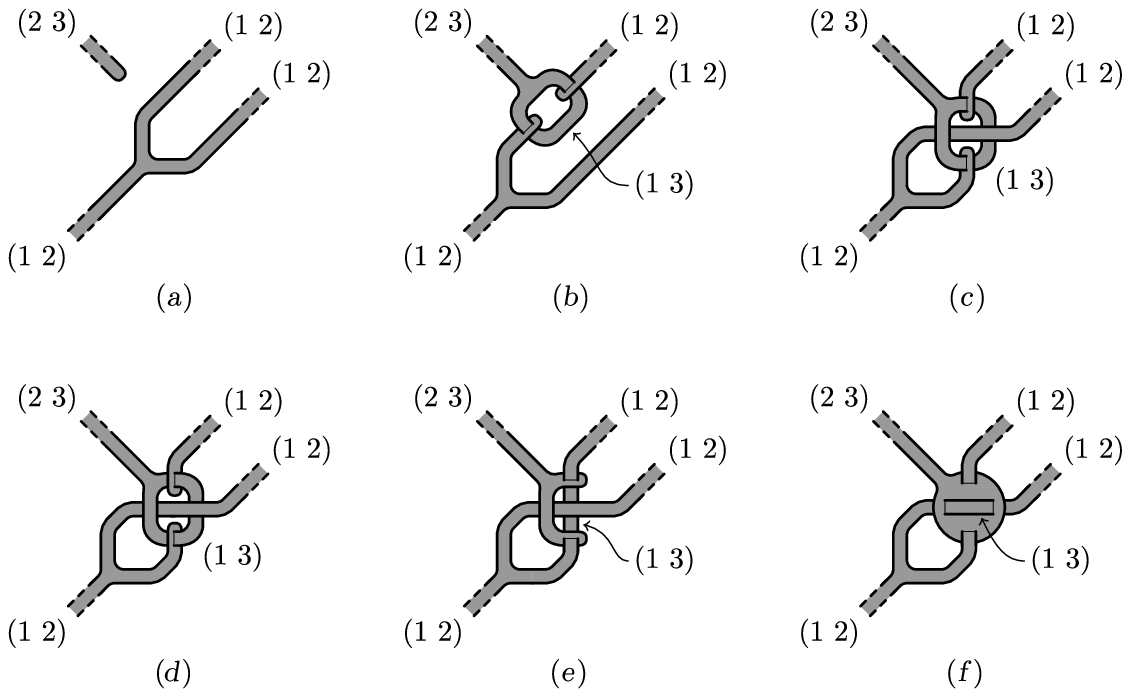}}
\end{Figure}

Then, we consider a disk $D_2$ spanned by $L'_2$ as in the proof of Proposition
\ref{diag/ribbon/thm} and perturb it near to $A_2$ in such a way that it becomes
disjoint from $\bA_2$, while remaining disjoint from all the other $A_i$'s and
continuing to form only clasps and ribbon intersections with the rest of the
surface. Such a perturbed disk can be used in place of the original one in the
process of Figure \ref{indepAlpha/fig}, to bring the stabilizing disk around
$\alpha_2$. After that, we cut $\alpha_2$ by a move $R_3^{-1}$ to obtain the
labelled surface of Figure \ref{sliding2/fig} \(d). At this point, we can continue
the sliding as indicated by the arrow and we can change all the crossings where
$\bA_2$ passes over $A_2$, by operating as in Figure \ref{indepOrder/fig}. 
In this way we get \(e), where $\bA'_2$ and $T'$ differ from $\bA_2$ and $T$ only by
the performed crossing changes. To end up with \(f), we first observe that $\bA'_2$
crosses always under $A_2$, so it can be pushed down below the plane $z = a_2$
(recall that $[a_2,b_2]$ is the eight interval of $L'_2$). Hence, after having
restored the band $\alpha_2$ by a move $R_3$, we can use the process of Figure
\ref{indepAlpha/fig} in the opposite direction, this time with the original disk
$D_2$, to take back the stabilizing disk.

Finally, we want to verify that the labelled ribbon surface resulting from all
the above modifications coincides with $F_{K'}$, where $K'$ is the ordinary Kirby
diagram obtained from $K$ by replacing $L_1$ with $L_1 \#_\gamma \bL_2$. 

Looking at Figure \ref{sliding2/fig} \(f), we call $\bL'_2$ the core of $\bA'_2$
and observe that here the original ribbon $A_1$ has been replaced by the ribbon $A_1
\#_\gamma \bA'_2$ with core $L'_1 \#_\gamma \bL'_2$. 
\mypagebreak
Taking into account the choices
made at the beginning about the height function of $L'$ and taking the care of
preserving the vertical triviality of $\bL'_2$ when pushing down $\bA'_2$, we can
assume that $L'_1 \#_\gamma \bL'_2, L'_2, \dots, L'_n$ form a vertically trivial
link.

We claim that, up to isotopy, this is a vertically trivial status of the link formed
by $L_1 \#_\gamma \bL_2, L_2, \dots, L_n$ and that the crossings at which the two
links differ are exactly the ones marked by the presence of a disk $B_i$. 
In fact, it is clear from the construction that, by inverting such crossings in the
vertically trivial link formed by $L'_1 \#_\gamma \bL'_2, L'_2, \dots, L'_n$, we get
a link of components $L_1 \#_\gamma \widehat L_2, L_2, \dots, L_n$, where $\widehat
L_2$ is a certain parallel copy of $L_2$. Then, our claim reduces to asserting that
$\bL_2$ and $\widehat L_2$ represent the same framing of $L_2$, that is $\Lk(L_2,
\bL_2) = \Lk(L_2, \widehat L_2)$. This equality between linking numbers follows from
some easy computations involving the writhes $w_2 = \Wr(L_2)$ and $w_2' = \Wr(L'_2)$
and the signed number $c_2 = (w_2 - w'_2)/2$ of the crossings of $L_2$ inverted to
get $L'_2$. Denoting by $f_2 = \Lk(L_2,\bL_2)$ the framing of $L_2$ in $K$,
we have $f_2 - w_2$ full twists inside the twist box $t$ of Figure \ref{sliding1/fig}
and $f_2 + 2c_2 - 2w_2$\break half twists inside the twist box $t'$ of Figure
\ref{sliding2/fig} (see step \(d) in the definition of $F_K$). As a consequence, the
additional crossings we inserted inside the twist box $T$ of Figure
\ref{sliding2/fig} is $- 2w'_2 - (f_2 + 2c_2 - 2w_2) = 2c_2 - f_2$. 
Then, the signed number of the crossings between $A_2$ and $\bA_2$ marked by the
$B_i$'s is $-f_2$, being $-2c_2$ the signed number of such crossings outside the
twist box $T$. Since both $\bA_2$ and $\bA'_2$ are unlinked from $A_2$, this number
of crossings remains unchanged if we replace $\bA_2$ with $\bA'_2$ and we can
conclude that $\Lk(L_2,\widehat L_2) = f_2$.

It remains to check that the ribbon $A_1 \#_\gamma \bA'_2$ in Figure
\ref{sliding2/fig} \(f) represents the right half integer framing of $L'_1
\#_\gamma \bL'_2$. To do that, we orient $L_1$, $L_2$, $L'_1$ and $L'_2$ accordingly
to Figure \ref{sliding1/fig} \(b). Then, the signed number of crossings of $L_1
\#_\gamma \bL_2$ to be inverted in order to get $L_1' \#_\gamma \bL'_2$ is $c_1 +
c_2 + \Lk(L_1,L_2)$, where $c_1$ is defined analogously to $c_2$.\break
On the other hand, the framing of $L_1 \#_\gamma \bL_2$ in $K'$ is $f_1 + f_2 +
2\Lk(L_1,L_2)$, where $f_1$ is the framing of $L_1$ in $K$, and so $A_1 \#_\gamma
\bA'_2$ should be equivalent up to vertical regular homotopy to a ribbon
representing the half integer framing $f_1/2 + f_2/2 + c_1 + c_2 + 2 \Lk(L_1,L_2)$.
The reader can easily realize that this is the case, taking into account that
$\Wr(L_1 \#_\gamma \bL_2) = w_1 + w_2 + 2 \Lk(L_1,L_2)$, where $w_1 = \Wr(L_1)$.
\end{proof}

\begin{remark} \label{smalldisks/rem}
Let us recall that any crossing in a Kirby diagram $K$ can be inverted, up to
2-deformation, by adding a suitable pair of 1/2-handles, as shown in Figure
\ref{crossing3/fig}. In the light of the preceding proposition, any disk $B_i$ in
the 4-stabilization of $F_K$ can interpreted, up to labelled 1-isotopy and ribbon
moves, as such a pair of 1/2-handles. A direct proof of this fact is provided in
Figure \ref{crossing4/fig}. Here, apart from labelled 1-isotopy, we perform two twist
transfers (moves $R_6^\pm$) to get \(c) from \(b) and one move $R_5$ followed by one
more twist transfer to get \(d) from \(c).
\end{remark}

\begin{Figure}[htb]{crossing3/fig}{}{}
\centerline{\fig{}{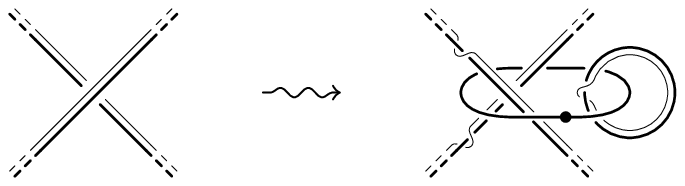}}
\end{Figure}

\begin{Figure}[htb]{crossing4/fig}{}{}
\centerline{\fig{}{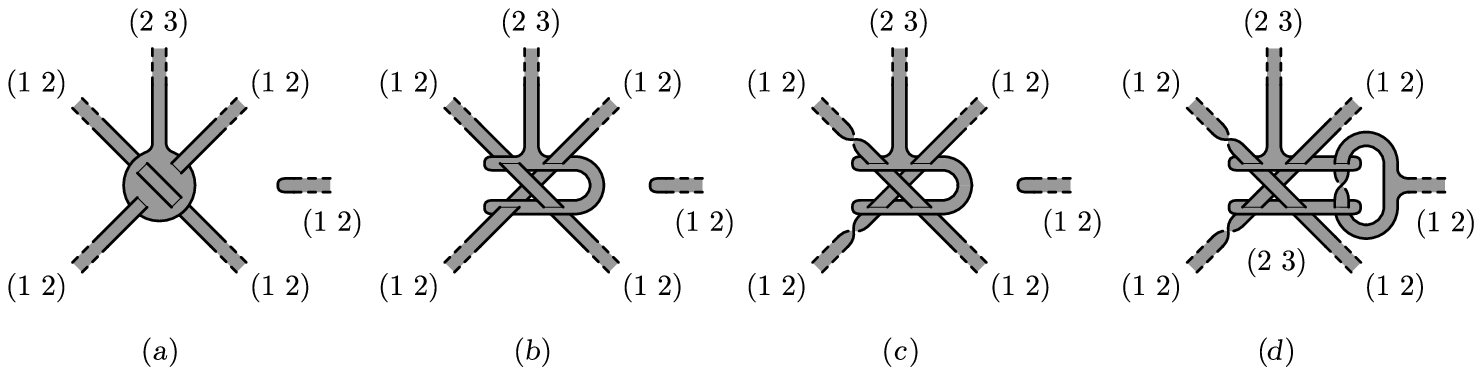}}
\end{Figure}

\break

To complete this section we are left with showing that, up to 2-deformations, the
construction of $F_K$ given above is inverted by that one of $K_F$ given in the
previous section. As observed at the beginning of the section, this implies that any
2-equivalence class of 4-dimensional 2-handlebodies can be represented as a simple
branched covering of $B^4$ (cf. Proposition \ref{surjectivity/thm}). In particular,
we can insist that the covering has degree 3 in the case of connected handlebodies.

\begin{proposition} \label{diag/ribbon/diag/thm}
Let $K$ be an ordinary Kirby diagram and $F = F_K$ be the corresponding labelled
ribbon surface. Then, the generalized Kirby diagram $K_F$ is equivalent to $K$ up
to 2-deformation moves.
\end{proposition}

\begin{proof}
Recall that for constructing $K_F$ one need first to choose an adapted 1-handlebody
structure on $F$, even if the 2-equivalence class of $K_F$ is independent on this
choice by Proposition \ref{1-equiv/2-equiv/thm}.

We claim that there exists an adapted 1-handlebody structure on $F$, naturally
related to the above construction, such that the generalized Kirby diagram $K_F$
constructed starting from it is equivalent to $K$ up to labelled isotopy, 1-handle
slidings and deletion of cancelling 0/1-handles.

Without loss of generality, we suppose that $K$ is already in the standard form
choosen in step \(a) of the construction of $F_K$. Moreover, we adopt all the
notations introduced during that construction.

To specify the claimed adapted handlebody structure of $F$, we first decompose each
ribbon $A_i \subset F$ as $A_i^0 \cup A_i^1$, where $A_i^0$ is a small 0-handle
containing the attach\-ing arc of $\alpha_i$ and $A_i^1$ is a 1-handle. Then, we
consider the adapted 1-handlebody structure of $F$ whose 0-handles are $A_0 \cup
\alpha_1 \cup \dots \cup \alpha_n \cup A_1^0 \cup \dots \cup A_n^0$, $B$,
$C_1, \dots, C_m$\break and whose 1-handles are $A_1^1, \dots, A_n^1$ (cf. Figure
\ref{cover2/fig}). By a suitable choice of the $\alpha_i$'s and the $\beta_i$'s
we can assume that all the 0-handles are blackboard parallel.

The generalized Kirby diagram $K_F$ constructed starting from this handlebody
structure is sketched in Figure \ref{diagcovdiag1/fig}. Here, as well as in all
the figures of this proof, we omit to draw the framings for the sake of readability.
Some further details of $K_F$ are shown in Figures \ref{diagcovdiag2/fig} and
\ref{diagcovdiag3/fig}. The labelled isotopy modifications described there are
performed at all the $\alpha_i$'s and at all the crossings between the $A_i$'s.

\begin{Figure}[htb]{diagcovdiag1/fig}{}{}
\centerline{\fig{}{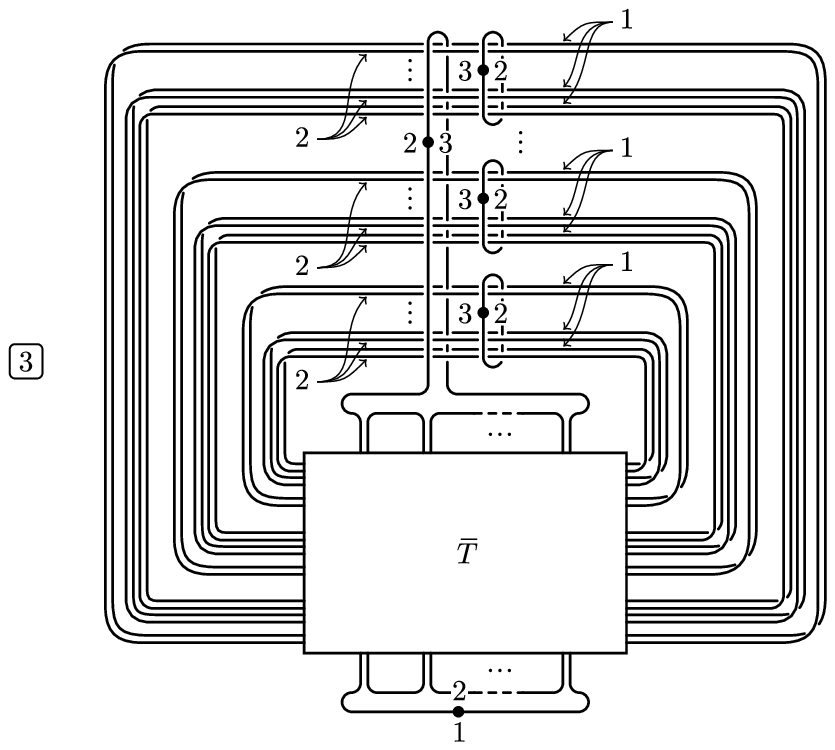}}
\end{Figure}

\begin{Figure}[htb]{diagcovdiag2/fig}{}{}
\smallskip\centerline{\fig{}{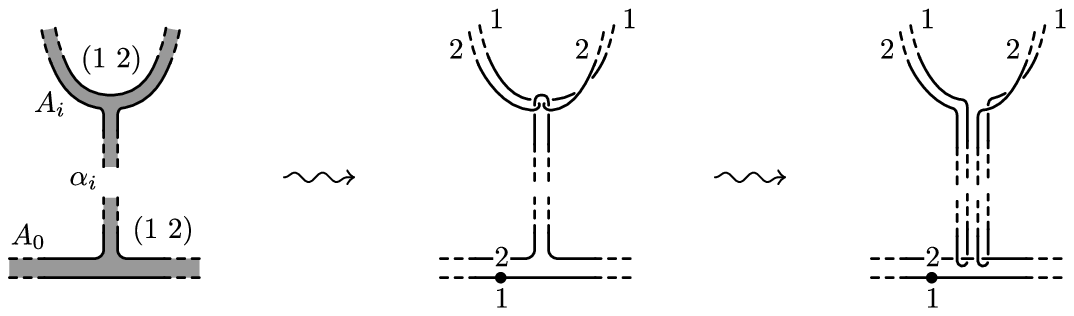}}
\end{Figure}

\begin{Figure}[htb]{diagcovdiag3/fig}{}{}
\smallskip\centerline{\fig{}{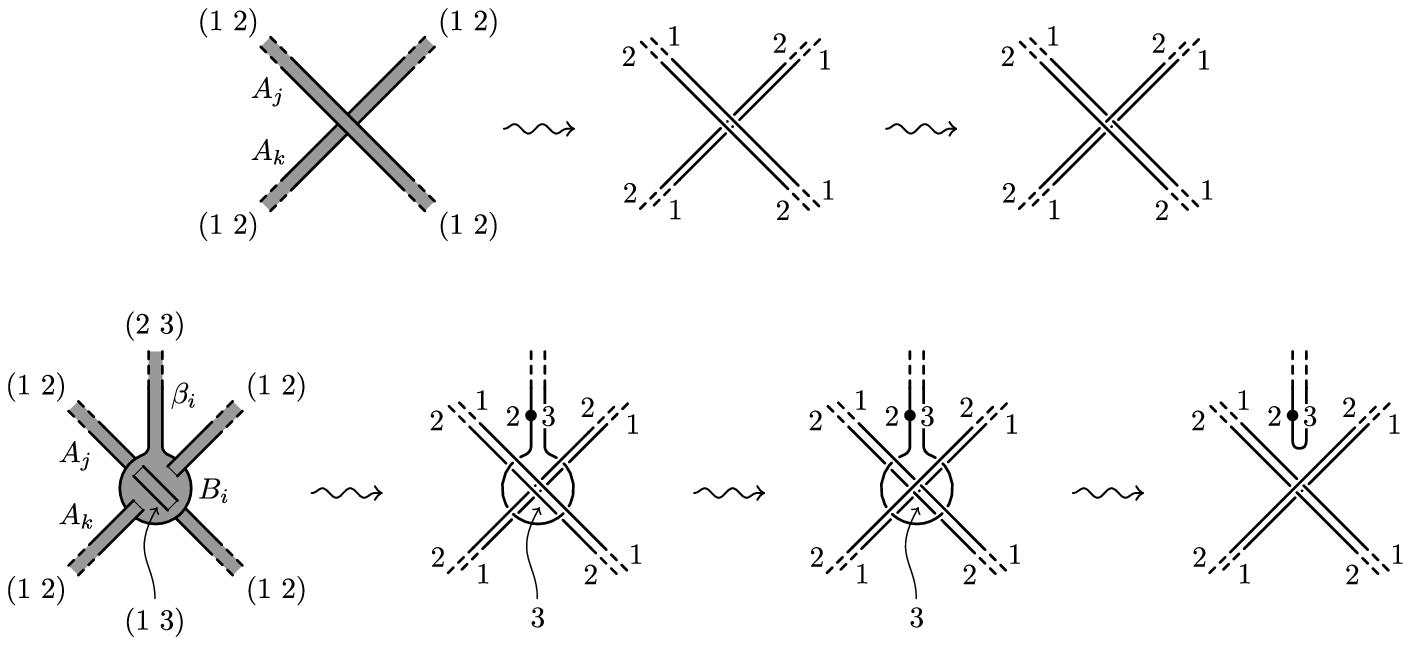}}
\end{Figure}

Once such modifications have been performed, we slide all the 1-handles corresponding
to the $C_i$'s over the one corresponding to $B$, in such a way that, up to diagram
isotopy, we are left with the diagram of Figure \ref{diagcovdiag4/fig}. Here we have
two overlapping but vertically separated tangle boxes $T$ (in front) and $T'$ (in
back), respectively labelled by 2 and 1.

Disregarding for the moment the framings, the link formed by the undotted components
in Figure \ref{diagcovdiag4/fig} is the componentwise band connected sum of the
original link $L$ and a parallel copy $L''$ of its vertically trivial state $L'$
pushed down to cross under everything else (including the dotted components). Each
component $L_i$ of $L$ 
\mypagebreak
is connected to the corresponding component $L_i''$ of
$L''$ by a band $\gamma_i$ running back and forth on the two sides of $\alpha_i$.
Since $L''$ is trivial and unlinked from the rest of the diagram and the bands
$\gamma_i$ can be assumed to be disjoint from a set of trivializing disks for $L''$,
we can isotope the diagram to get back $K$ entirely labelled by 2 with two extra
dotted components labelled by 1,2 and 2,3 separated from it. Finally, such dotted
components can be eliminated by 0/1-handle cancellation.

\begin{Figure}[htb]{diagcovdiag4/fig}{}{}
\centerline{\fig{}{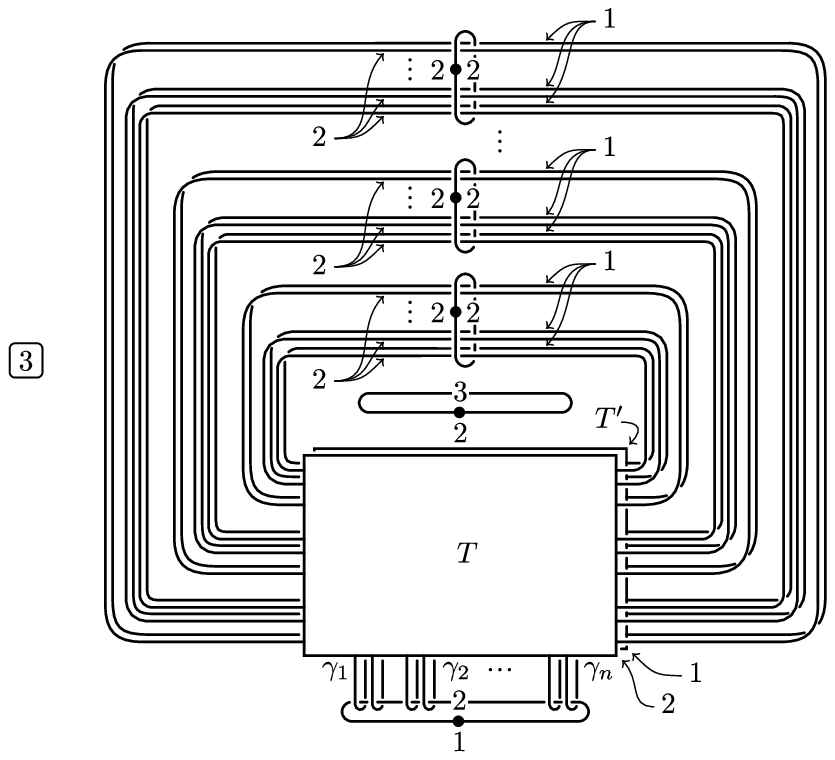}}
\end{Figure}

Consider now the framings. To verify that the final framings we obtain coincide with
the original ones, we proceed like in the last part of the proof of Proposition
\ref{diag/ribbon/thm}.\break Let $f_i$ be the framing of $L_i$ in $K$ and $c_i$ be
the signed number of the crossings of $L_i$ inverted to get $L'_i$, that is $c_i =
(w_i - w'_i)/2$ where $w_i = \Wr(L_i)$ and $w'_i = \Wr(L'_i)$. We observe that
the framing of $L_i \#_{\gamma_i} L''_i$ in the diagram of Figure
\ref{diagcovdiag4/fig} is the band connected sum of two half integer framings along
$L_i$ and $L''_i$, both of which differ from the blackboard framing by $f_i + 2 c_i
- 2 w_i$ half twists. Hence, for $L_i \#_{\gamma_i} L''_i$ we have $f_i + 2c_i
- 2w_i$ full twists added to the blackboard framing. After we have performed the
0/1-handle cancellations to reduce the diagram to an ordinary one, the blackboard
framing of $L_i \#_{\gamma_i} L''_i$ can be encoded in the usual way by the
integer\break $w_i + w'_i$. So we are done, since $(f_i + 2c_i - 2w_i) + (w_i + w'_i)
= f_i$.
\end{proof}

For future reference (see Remark \ref{trivialstate/rem}), we observe that the proof
of the above Proposition \ref{diag/ribbon/diag/thm} still works if we assume that
the link $L'$ used in the construction of $F_K$ is any trivial state of $L$ and not
necessarily a vertically trivial state of it.\break In fact, the vertical triviality
of $L'$ has been used only to conclude that $L''$ is trivial and for that the
triviality of $L'$ sufficies.

\begin{proposition} \label{surjectivity/thm}
Any orientable 4-dimensional 2-handlebody $H$ with $c$ connected components is
2-equivalent to a special one having generalized Kirby diagram of the form $K_F$,
for some labelled (orientable) ribbon surface $F \subset B^4$ representing $H$ as a
simple branched covering of $B^4$ of degree 3$c$.
\end{proposition}

\begin{proof}
Up to 1-handle sliding and deletion of cancelling 0/1-handles, we can assume that $H$
has one 0-handle in every component. So, it can be represented by a generalized
Kirby diagram $K$ which is the disjoint union of $c$ ordinary Kirby diagrams $K_1,
\dots, K_c$, such that each $K_i$ is separated from all the others and is entirely
labelled by $i$. Disregarding these labels, we construct the labelled ribbon surfaces
$F_{K_1}, \dots, F_{K_c}$. Then, we put $F = F_{K_1} \sqcup \dots \sqcup F_{K_c}$,
after the labels $1,2,3$\break of each $F_{K_i}$ has been replaced respectively by
$3i - 2, 3i-1,3i$.
The 4-dimensional 2-handlebody represented by $F$ as a $3c$-fold branched covering
of $B^4$ can be proved to be 2-equivalent to $H$, by applying Proposition
\ref{diag/ribbon/diag/thm} componentwisely.
For the orientability of $F$, we refer to Remark \ref{orient/rem}.
\end{proof}

\begin{remark} \label{surjectivity/rem}
Notice that in both the above Propositions \ref{diag/ribbon/diag/thm} and
\ref{surjectivity/thm}, once a suitable 1-handlebody structure is fixed on $F$,
only 1-handle sliding and addition/deletion of cancelling 0/1-handles are needed to
get the wanted 4-dimensional 2-handlebody from $K_F$, up to handle isotopy. 
\end{remark}

\section{The equivalence theorems\label{equivalence/sec}}

This section completes the proof of the four equivalence theorems stated in the
Introduction. The first and main step, is to show that 4-dimensional 2-handlebodies
up to 2-deformation are bijectively represented, through the map $F \mapsto K_F$, by
simply labelled ribbon surfaces up to labelled 1-isotopy, stabilization and ribbon
moves $R_1$ and $R_2$, besides labelling conjugation (remember that labelling is
actually defined only up to conjugation in $\Sigma_d$).

Like in the previous section, we first restrict our attention to the connected
case and then come back to the general case with Proposition \ref{bijectivity/thm}.
Recall that, for the connected case, we also have the map $K \mapsto F_K$, which
associates to each ordinary Kirby diagram $K$ a labelled ribbon surface $F_K$
(defined up to labelled 1-isotopy and ribbon moves) representing its 2-equivalence
class as a 3-fold simple branched covering of $B^4$. In the light of Proposition
\ref{diag/ribbon/diag/thm}, we will be done once we prove that such map is
surjective up to labelled 1-isotopy, stabilization and ribbon moves $R_1$ and $R_2$.
This is the aim of Propositions \ref{to3fold/thm} and \ref{tospecial/thm}.

\medskip

To begin with, we notice that a $d$-fold simple branched covering of $B^4$
represented by a labelled ribbon surface $F \subset B^4$ is connected if and only if
the transpositions which appear as labels of any diagram of $F$ generate a
transitive subgroup of the symmetric group $\Sigma_d$. This is trivially equivalent
to say that they generate all $\Sigma_d$.

In particular, in this case we can use labelled 1-isotopy move $I_2$ to expand from
$F$ a tongue which, after a suitable sequence of ribbon intersections, is labelled
with any given transposition $\tau \in \Sigma_d$ on its tip. Passing all the rest of
the diagram through the tip of such a tongue and putting everything back in the
original position, has the same effect as conjugating all the labels by $\tau$.
Hence, any labelling conjugation can be obtained by a suitable labelled 1-isotopy.
This is the reason why labelling conjugation does not appear in the statements of our
equivalence theorems concerning connected coverings, while it does in Proposition
\ref{bijectivity/thm}.

Before going on, we also introduce the following notion of {\sl special position} for
a labelled ribbon surface $F \subset B^4$ representing a (possibly disconnected)
simple branched covering of $B^4$. We say that $F$ is in special position if its
diagram is entirely contained in the projection plane except for a finite number of
ribbon half twists and of ribbon intersections and crossings as the ones depicted in
Figure \ref{special1/fig} with $i$, $j$, $k$ and $l$ all distinct and $i < j < k$.

\begin{Figure}[htb]{special1/fig}{}{}
\centerline{\fig{}{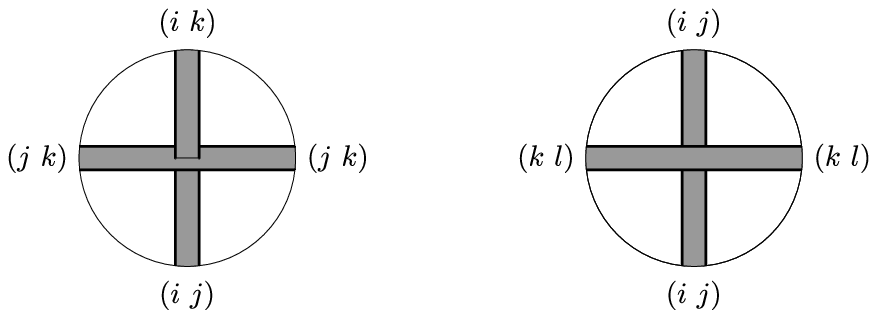}}
\end{Figure}

Labelled ribbon surfaces in special position have some remarkable properties that
will be useful in the next proofs. Namely, any such $F$ is the disjoint union of
subsurfaces $F_1, \dots, F_{d-1}$, where $d$ is the degree of the covering, such
that: 1) the labels attached to $F_i$ are all of the type $(i\ j)$ with $j = i+1,
\dots, d$, and so $F_{d-1}$ is entirely labelled by $(d{-}1\ d)$; 2) $F_i$ does not
form ribbon self intersections or self crossings, that is its diagram can be
considered planar except for ribbon half twists; 3) all the ribbon intersections of
$F$ consist of a ribbon of $F_i$ which pass through a ribbon of $F_j$ with $i < j$,
hence $F_1$ is nowhere passed through by any other $F_i$.

Clearly, special position is quite restrictive. For example, even the very peculiar
labelled ribbon surfaces $F_K$ are not in special position, due to the ribbon
crossings inside the tagle box and to the disks $B_i$. Nevertheless, the next lemma
tells us that things are different if we reason up to ribbon moves.

\begin{lemma}\label{specialposition/thm}
Any labelled ribbon surface representing a connected simple branched covering of
$B^4$ of degree $d \geq 3$ can be put in special position through labelled 1-isotopy
and moves $R_1$ and $R_2$.
\end{lemma}

\begin{proof}
Let $F \subset B^4$ be a labelled ribbon surface as in the statement. Forgetting
the labelling restrictions of Figure \ref{special1/fig}, labelled diagram isotopy
allows us to make the diagram of $F$ entirely contained in the projection plane
except for a finite number of ribbon half twists and of ribbon intersections and
crossings. We omit the details of this essentially trivial step and focus on the
task of eliminating the ribbon intersections and crossings which do not satisfy the
above labelling restrictions.

We change any ribbon intersection between ribbons with disjoint monodromies into a
crossing, by a move $R_2$. Moreover, we change any crossing between ribbons with
non-disjoint monodromies into two ribbon intersections, by the first labelled
1-isotopy move of Figure \ref{special2/fig}, where $k$ may or may not be equal to
$j$. Then, we apply the second labelled 1-isotopy move of Figure
\ref{special2/fig}, where $k \not\in \{i,j\}$, to eliminate all the ribbon
intersections between ribbons with the same monodromy. Here, we use the hypotheses
that the covering is connected and has degree $d \geq 3$, to get the tongue labelled
by $(i\ k)$ on its tip. We choose such a tongue to minimize the number of ribbon
intersections and crossings, so that none of these is formed with a ribbon having the
same monodromy. As above, we replace any crossing with a ribbon having non-disjoint
monodromy by two ribbon intersections and any ribbon intersection with a ribbon
having disjoint monodromy by a crossing.

\begin{Figure}[htb]{special2/fig}{}{}
\centerline{\fig{}{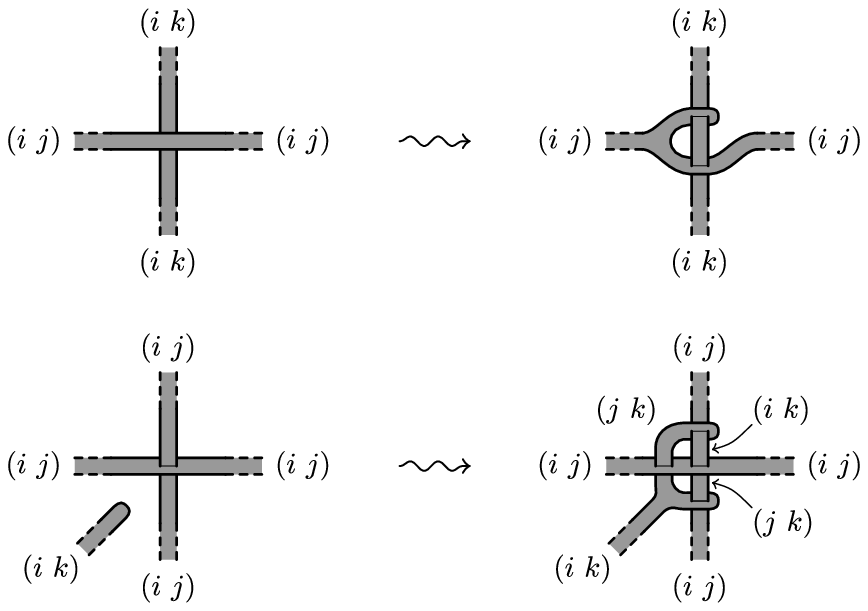}}
\end{Figure}

Thus, we are left only with ribbon intersections and crossings as in Figure
\ref{special1/fig}, with $i$, $j$, $k$ and $l$ all distinct. The ones which do not
satisfy the inequalities $i < j < k$ can be eliminated, by performing one move
$R_5$ followed by two moves $R_1^{\pm1}$ for the ribbon intersections and just one
move $R_4$ for the crossings.
\end{proof}

Let us now pass on to the announced Propositions \ref{to3fold/thm} and
\ref{tospecial/thm}.

\begin{proposition}\label{to3fold/thm}
Up to labelled 1-isotopy and moves $R_1$ and $R_2$, any labelled ribbon surface $F
\subset B^4$ representing a connected simple branched covering of $B^4$ of degree $d
\geq 3$ is equivalent to the $d$-stabilization of a labelled ribbon surface $F'
\subset B^4$ representing a simple 3-fold branched covering of $B^4$.
\end{proposition}

\begin{proof}
We proceed by induction on $d$. For $d = 3$ there is nothing to prove. Given $F$ 
as in the statement with $d > 3$, we prove that it is equivalent to  the\break
$d$-stabilization of a labelled ribbon surface representing a simple branched
covering of $B^4$ of degree $d-1$.

To prove the inductive step, we first put $F$ in special position, by applying Lemma
\ref{specialposition/thm}, and modify it in such a way that the label $(1\ d)$ does
not appear anymore in its diagram. 

Notice that, all the labels $(1\ d)$ of $F$ are attached to the subsurface $F_1
\subset F$, consisting of the pieces of $F$ labelled by $(1\ i)$, with $i = 2, \dots,
d$. As we said after the definition of special position, $F_1$ does not form
ribbon self intersections or self crossings and is nowhere passed through by any
other component of $F$. Moverover, no piece of $F_1$ labelled by $(1\ d)$ is crossed
over by any ribbon.

Consider an adapted 1-handlebody decomposition on $F$ such that crossings and half
twists only occur along 1-handles. On the 0-handles of $F_1$ which are labelled by
$(1\ d)$, we operate as in Figure \ref{special3/fig}, where $1 < i < d$. By choosing
the tongue labelled $(i\ d)$ to minimize the number of ribbon intersections and
crossings, we can preserve special position (cf. proof of Lemma
\ref{specialposition/thm}).

\begin{Figure}[htb]{special3/fig}{}{}
\centerline{\fig{}{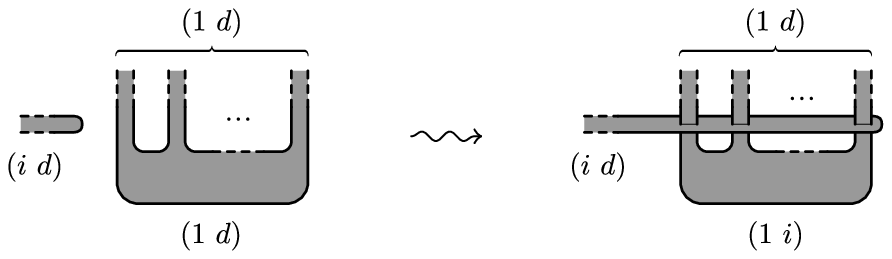}}
\end{Figure}

\begin{Figure}[htb]{special4/fig}{}{}
\centerline{\fig{}{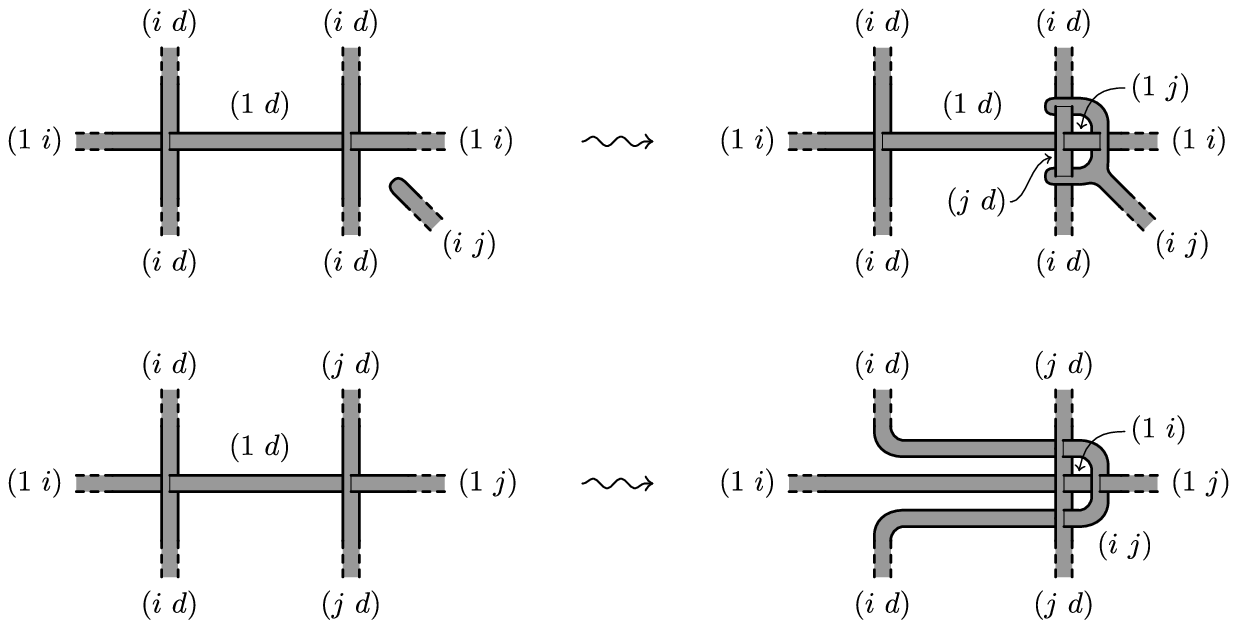}}
\end{Figure}

After that, only some segments of 1-handles delimited by ribbon intersections are
still labelled by $(1\ d)$, as sketched on the left side of Figure \ref{special4/fig}, where $1 < i < j < d$. Here, we have two cases, depending on whether the two delimiting ribbons
\mypagebreak
have the same
label or not. The upper part of the Figure shows how to reduce the first case to the
second, while the lower part tells us how to eliminate the label $(1\ d)$ in this
second case. In both cases, we leave to the reader to restore special position and to
check that no problem arise with ribbons which possibly cross under the tract labelled
by $(1\ d)$.

Once the label $(1\ d)$ has been eliminated from the diagram of $F$, while
preserving special position, we push $F_1$ down below all the rest of $F$, except
for some tongue terminating at a ribbon intersection, as suggested by right side of
Figure \ref{special5/fig}. This can be done by vertical diagram isotopy and moves
$R_4$ at the ribbon crossings where $F_1$ crosses above $F - F_1$. Then, we slide
$F_1$ horizontally under $F - F_1$ to make the diagram as in Figure
\ref{special5/fig}, where $F_1$ is contained in the lower box and $F - F_1$ in the
upper one, apart from the ribbons connecting the two boxes.\break Notice that the
labels in the upper (resp. lower) box do not involve 1 (resp. $d$), while the labels
of the connecting ribbons do not involve both 1 and $d$.

\begin{Figure}[htb]{special5/fig}{}{}
\centerline{\fig{}{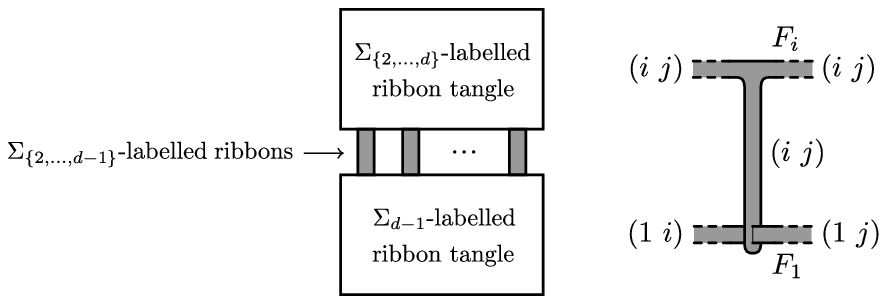}}
\end{Figure}

Finally, the modifications described in Figure \ref{special6/fig} allows us to
isolate a stabilizing disk labelled by $(1\ d)$, by removing $d$ from all the other
labels.

\begin{Figure}[htb]{special6/fig}{}{}
\centerline{\fig{}{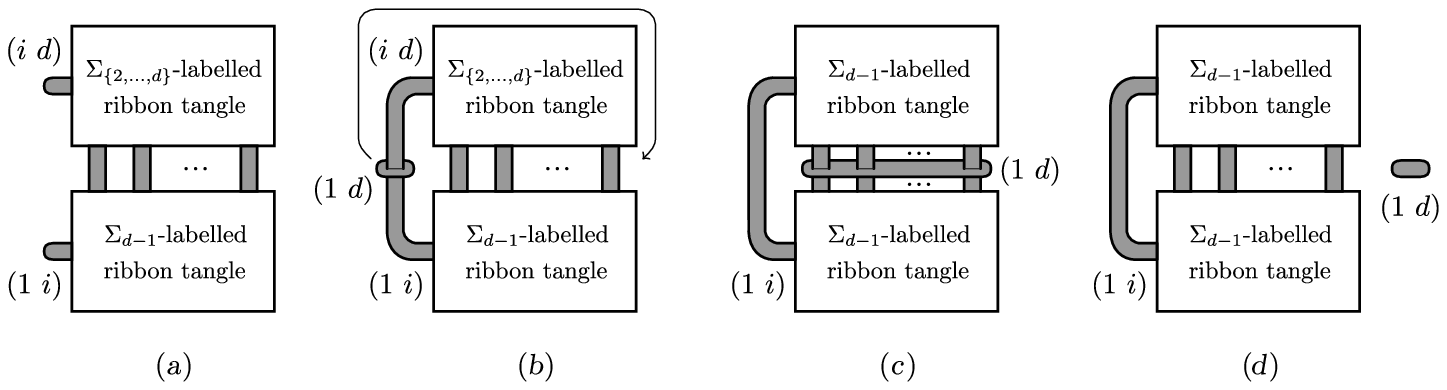}}
\end{Figure}

Namely, we expand from the boxes two tongues labelled $(1\ i)$ and $(i\ d)$ for some
$i = 2, \dots, d-1$, as in \(a). This can be always done, possibly after having
expanded some other $\Sigma_{\{2,\dots,d-1\}}$-labelled tongues connecting the two
boxes, in order to make the traspositions in the upper (resp. lower) box generate
all the symmetric group $\Sigma_{\{2,\dots,d\}}$ (resp. $\Sigma_{d-1}$). Then, we
connect the tips of the two above tongues by a move $R_3$ and use labelled 1-isotopy
to move the resulting new disk with label $(1\ d)$ as indicated by the arrow in \(b).
Eventually, we get the diagram in \(c), where also the upper box takes labels in
$\Sigma_{d-1}$, as well as the lower one, so that the only label involving $d$ is
the one of the disk between the two boxes. Such disk can be disentagled from the
ribbons connecting the boxes by using move $R_2$, to get \(d).
\end{proof}

\begin{proposition}\label{tospecial/thm}
For any labelled ribbon surface $F \subset B^4$ representing a connected 3-fold
simple branched covering of $B^4$, there exists an ordinary Kirby diagram $K$ such
that the 4-stabilizations of $F$ and $F_K$ are equivalent up to labelled 1-isotopy
and moves $R_1$ and $R_2$.
\end{proposition}

\begin{proof}
By Lemma \ref{specialposition/thm}, we can suppose $F$ to be in special position.
In this case, as we said after the definition of special position, $F$ is the
disjoint union of two non-empty subsurfaces $F_1$ and $F_2$, the first of which takes
labels $(1\ 2)$ and $(1\ 3)$, while the second one is entirely labelled by $(2\ 3)$.
Moreover, the diagram of $F$ cannot have any ribbon crossing, since there are no
disjoint transpositions in $\Sigma_3$, and all the ribbon intersections are formed
by $F_1$ passing through $F_2$. These can be polarized to have planar projection as
in the left side of Figure \ref{special1/fig} with $i = 1$, $j = 2$ and $k = 3$, up
to labelled diagram isotopy which locally half twists the horizontal ribbon.

\begin{Figure}[htb]{ribbondiag1/fig}{}{}
\centerline{\fig{}{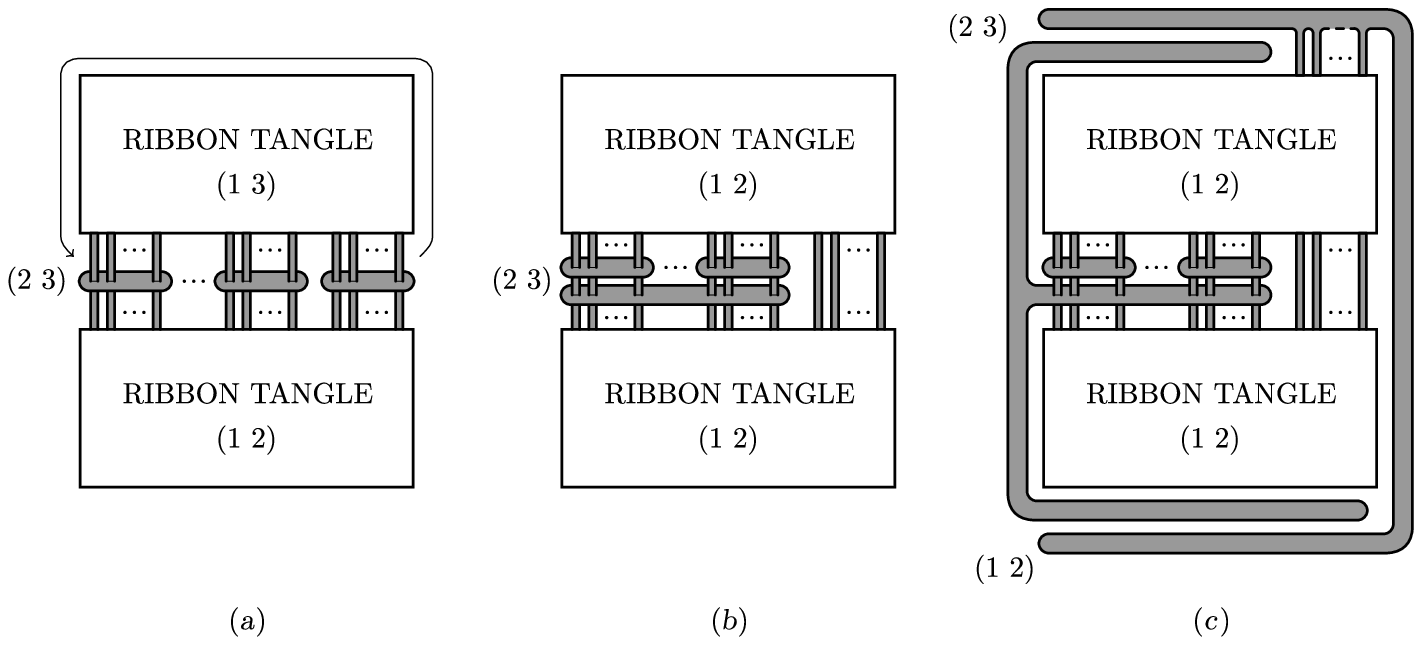}}
\end{Figure}

Consider an adapted 1-handlebody decomposition of $F$ such that half twists only
occur along 1-handles. By move $I_2$ and the tongue technique already seen in the
previous proofs, we insert a ribbon intersection along each 1-handle of $F_2$,
taking care that special position is preserved. Then, we apply a move $R_5$ at
every ribbon intersection of $F$. After that, $F_2$ is a disjoint union of disks and
we can use move $R_6$ to flatten its diagram into the projection plane, still
preserving special position and the above polarization of the ribbon intersections. 
Finally, a labelled diagram isotopy suffices to put $F$ into the form depicted in
Figure \ref{ribbondiag1/fig} \(a). Such an isotopy can be realized in two steps: 1)
lift all the $(1\ 3)$-labelled parts of $F_1$ above the projection plane and push
all the $(1\ 2)$-labelled ones below it, by a vertical isotopy fixing $F_2$; 2)~move
the planar diagram of $F$ to the wanted form, by a suitable horizontal labelled
isotopy. Of course, this last step does not preserve any more the special position.

Let us assume that both tangle boxes in Figure \ref{ribbondiag1/fig} \(a) are
non-empty and that there are at least two $(2\ 3)$-labelled disks between them. We
leave to the reader to see that such assumption can be made without loss of
generality.

By labelled 1-isotopy, we move the rightmost $(2\ 3)$-labelled disk as suggested by
the arrow, to form a long bar under the other ones like in \(b). During this process
all the labels in the upper box are changed in $(1\ 2)$. Then, we obtain the four
bars at top and bottom which appear in \(c) by labelled diagram isotopy. In
particular, the ones labelled by $(1\ 2)$ are expanded from a 0-handle of $F$
picked up from the upper box.

We warn the reader that the groupings of the vertical bands at different levels
in Figure \ref{ribbondiag1/fig} \(c), as well as in Figure \ref{ribbondiag4/fig}
below, are totally uncorrelated. Their apparent correspondence in the diagrams has
only a pictorial value.

\begin{Figure}[b]{ribbondiag2/fig}{}{}
\centerline{\fig{}{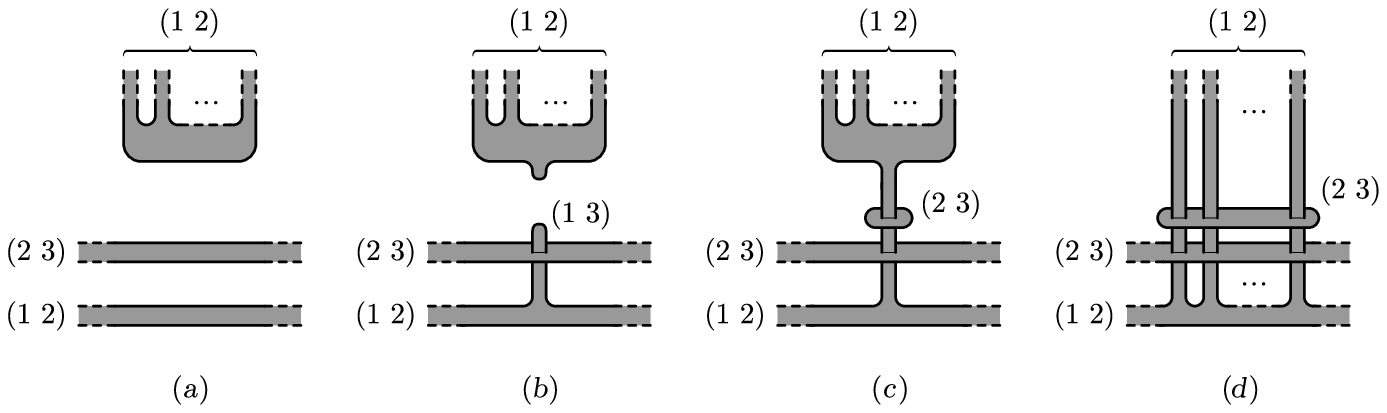}}
\end{Figure}

The following Figure \ref{ribbondiag2/fig} shows how to incorporate all the
0-handles in the lower box of Figure \ref{ribbondiag1/fig} \(c) into the $(1\
2)$-labelled bar at bottom. Here, apart from labelled 1-isotopy, only one move $R_3$
occurs between \(b) and \(c). Similarly, all the 0-handles in the upper box can be
incorporated into the $(1\ 2)$-labelled bar at top. After that, the two ribbon
tangles consist of a certain number of bands which are attached directly to the
top/bottom of the $(1\ 2)$-labelled bar. We subdivide such bands, by inserting new
0-handles at the intermediate minima and maxima, in such a way that each one of the
resulting pieces runs monotonically with respect to the vertical direction of the
diagram plane. By labelled diagram isotopy, all the 0-handles corresponding to
minima (resp. maxima) inside upper (resp. lower) box can be moved to the lower
(resp. upper) one.\break Then, also the new 0-handles can be incorporated into the
bars at top and bottom to get a diagram as in Figure \ref{ribbondiag4/fig} \(a),
where the ribbon tangles of Figure \ref{ribbondiag1/fig} \(c) are replaced by ribbon
braids. Denote by $X$ and $Y$ the corresponding ordinary braids, disregarding the
ribbon half twists (cf. diagram \(b) of Figure \ref{ribbondiag4/fig}).

\begin{Figure}[htb]{ribbondiag4/fig}{}{}
\vskip3pt\centerline{\fig{}{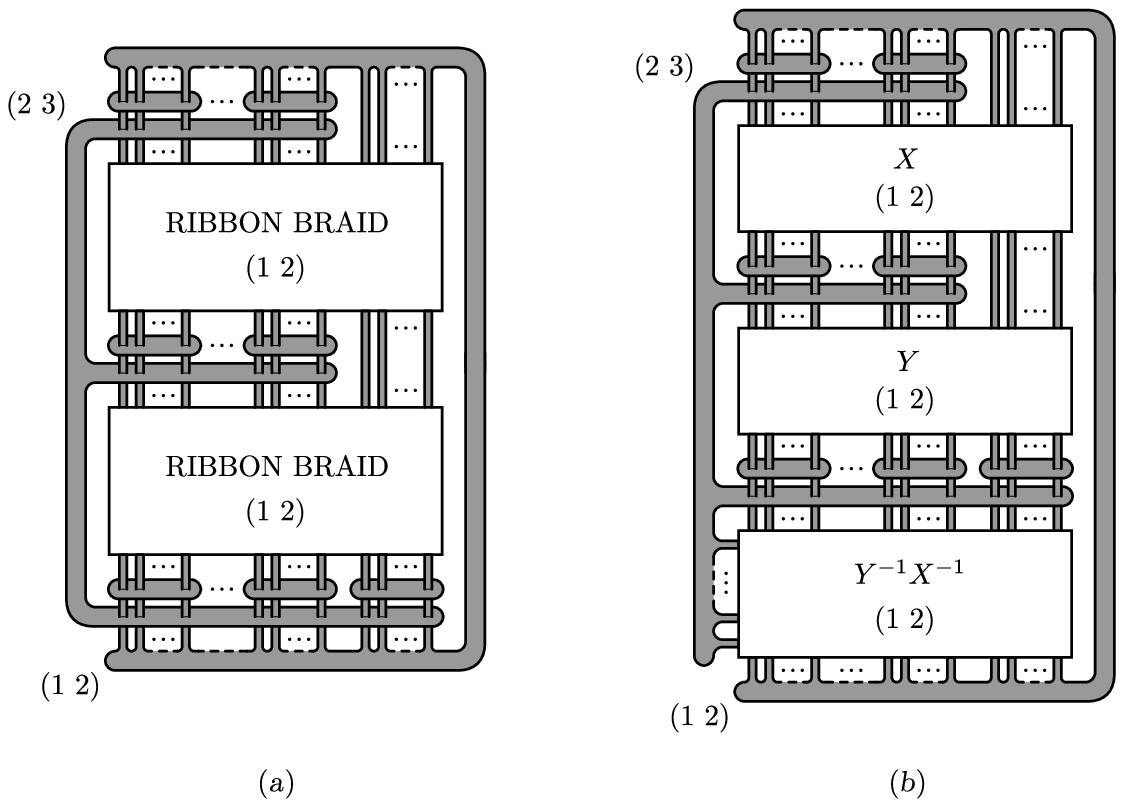}}
\end{Figure}

Our next goal is to insert in the diagram a third box with a ribbon braid
representing the blackboard framing of $Y^{-1}X^{-1}$, as in Figure
\ref{ribbondiag4/fig} \(b). The ribbon crossing relative to a standard generator of
the braid group can be added just above the bottom bar, together with a small disk
expanded from the $(2\ 3)$-labelled vertical bar on the left side, as shown in Figure
\ref{ribbondiag3/fig}. Such a modification essentially coincides with the one
described in Figure \ref{sliding4/fig}, thus we already know how to realize it in
terms of labelled 1-isotopy and ribbon moves. The inverse generator can be dealt
with similarly. That is enough to get Figure \ref{ribbondiag4/fig} \(b).

\begin{Figure}[htb]{ribbondiag3/fig}{}{}
\centerline{\fig{}{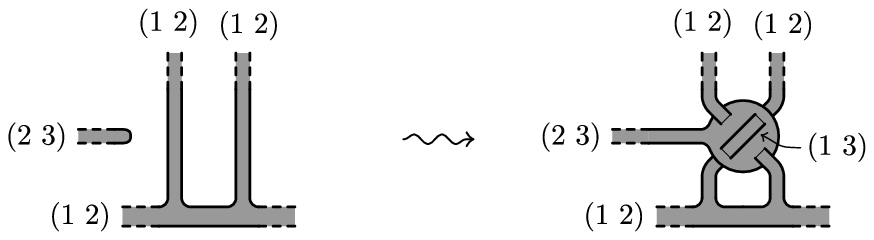}}
\end{Figure}

Now, by a labelled diagram isotopy, we align the $(2\ 3)$-labelled horizontal bars as
sketched in Figure \ref{ribbondiag5/fig}. Here, the bars are numbered to make
clear the isotopy and, for the sake of readability, only one of the bands forming
the ribbon braid of Figure \ref{ribbondiag4/fig} \(b) is drawn. During the
isotopy, all the other bands are kept parallel to this one outside the boxes.
Looking at the right side of Figure \ref{ribbondiag4/fig}, we see that the isotopy
destroys the braid structure, by introducing self-crossings along the bands. However,
such self-crossings are at most six for each band and satisfy a property that will be
crucial in the following: going from bottom to top, at each self-crossing the band
passes first under and then over.

\begin{Figure}[htb]{ribbondiag5/fig}{}{}
\centerline{\fig{}{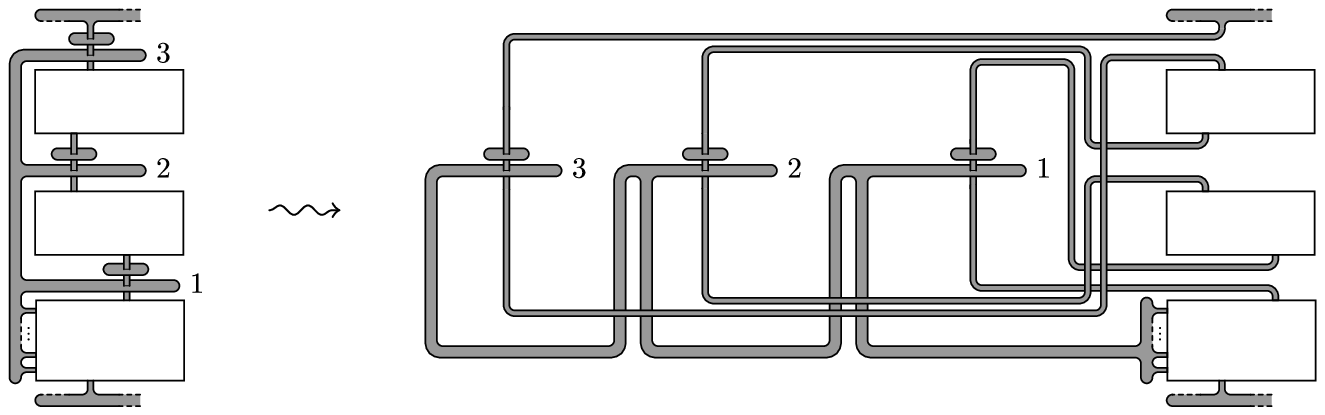}}
\end{Figure}

Using the stabilization disk as in Figure \ref{indepDiag4/fig}, the bars we have
just aligned can be disconnected in turn from the rest of the diagram and then
reconnected differently to form a single long bar like in Figure
\ref{ribbondiag6/fig} \(a). The tangle box in this last diagram includes the three
braid boxes of the previous one together with the obvious part of the bands outside
them. To get the subsequent diagram of Figure \ref{ribbondiag6/fig} \(b), we simply
perform a 90$^\circ$ clockwise rotation and contract the $(1\ 2)$-labelled bar
(sliding consequently all the bands connecting it to the box).

\begin{Figure}[htb]{ribbondiag6/fig}{}{}
\centerline{\fig{}{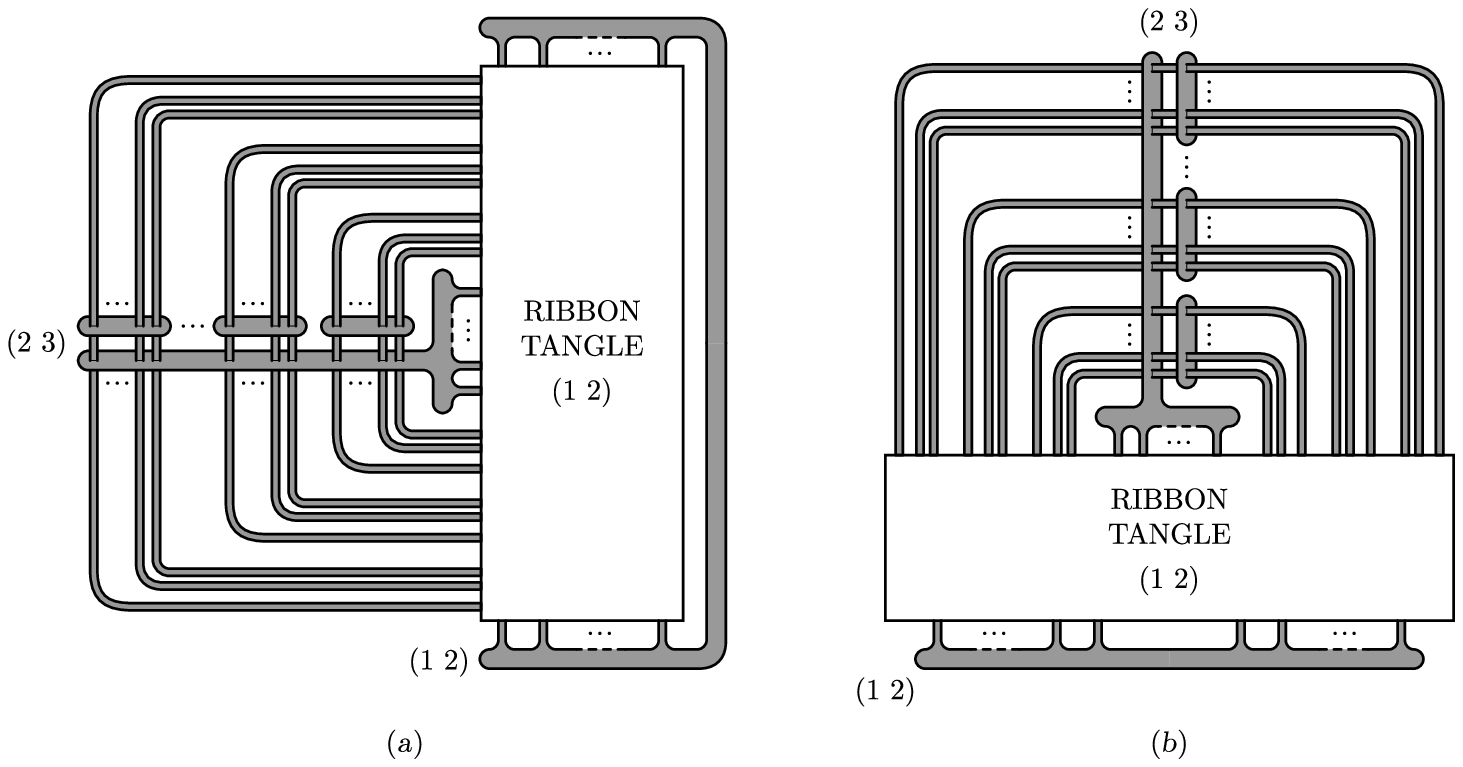}}
\end{Figure}

\begin{Figure}[htb]{ribbondiag7/fig}{}{}
\centerline{\fig{}{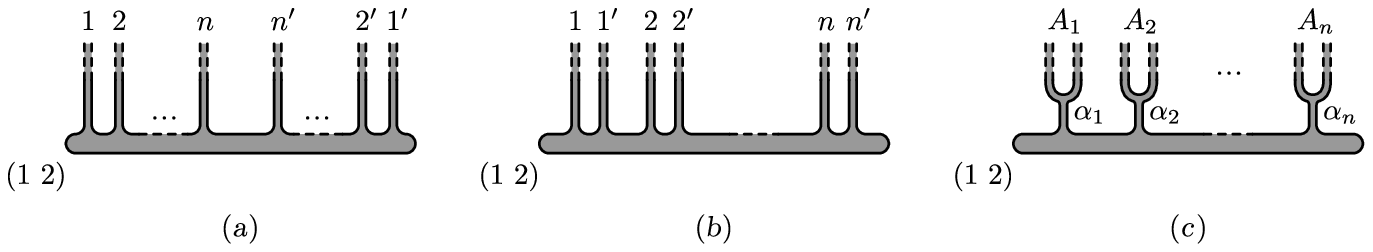}}
\end{Figure}

Figure \ref{ribbondiag7/fig} \(a) indicates the order in which the bands deriving
from the strings of the ribbon braid of Figure \ref{ribbondiag4/fig} \(b) are
attached to the bottom bar in Figure \ref{ribbondiag6/fig} \(b). Here, we numbered
by $i$ and $i'$ the two ends of the band corresponding to the $i$-th string. The same
modification depicted in Figure \ref{ribbondiag3/fig} we have already used before,
enables us to change this order, by pairing the two ends of the same band as in \(b).
After that, we can think of the $i$-th band as a (possibly non-orientable) closed
ribbon $A_i$, connected to the bottom bar by a band $\alpha_i$, as suggested in \(c).

Up to sliding the $(1\ 2)$-labelled bands entering the tangle box from the top edge
to the lateral ones, the resulting diagram looks like the one of Figure
\ref{cover2/fig}, except for the fact that the link $L' = L'_1 \cup \dots \cup L'_n$
formed by the cores of the $A_i$'s may not be vertically trivial. Indeed, the
triviality of the braid in Figure \ref{ribbondiag4/fig} \(b) implies that $L'$
is trivial, but not necessarily vertically trivial. On the other hand, each single
$L_i$ is vertically trivial, since its only self-crossings are the ones introduced
in Figure \ref{ribbondiag6/fig} and they satisfy the property pointed out when
introducing the figure. Therefore, we only need to worry about vertically separating
different $L'_i$'s.

To this end, let us observe that the triviality of $L'$ is enough to construct
disjoint disks $D_1, \dots, D_n$, with the same properties as in the proof
of Proposition \ref{diag/ribbon/thm}. We can use these disks in turn to vertically
separate the $L'_i$'s, just like we did there for proving the independence of
$F_K$ from the vertical order of the components of $L'$ (cf. Figure
\ref{indepAlpha/fig} and  Figure \ref{indepOrder/fig}).

To conclude the proof, it remains to verify that the $\alpha_i$'s can be put in the
right position as prescribed by the definition of $F_K$. We leave this trivial task
to the reader.
\end{proof}

\begin{remark}\label{trivialstate/rem}
As a consequence of Proposition \ref{tospecial/thm}, the link $L'$ in step \(c) of
the construction of $F_K$ at page \pageref{stepc} can be chosen to be any trivial 
state (not necessarily a vertically trivial one) of the link $L$ consisting of the
framed components of $K$, without losing the well-definedness of $F_K$ up to
4-stabilization, ribbon moves and 1-isotopy. The reason is that, by Proposition
\ref{diag/ribbon/diag/thm} and the observation immediately following its proof, the
resulting labelled ribbon surface $F_K$ does represent the same 2-deformation class of
4-dimentional 2-handlebodies as if we had chosen $L'$ to be a vertically trivial state
of $L$. Then, Propositions \ref{tospecial/thm} and \ref{2def/moves/thm} allow us to
conclude in a straigthforward way.
\end{remark}

As an example, in Figure \ref{gompf/fig} we present the Kirby diagram $K$ and the
corresponding ribbon surface $F_K$ for the Akbulut-Kirby 4-sphere $\Sigma_n$ with 
$n = 3$. The Kirby diagram is the same as the one drawn in figure 4 of \cite{Go91},
where it is shown that $\Sigma _n$ is diffeomorphic to $B^4$ for any $n$ and it is
also conjectured that it is not 2-equivalent to $B^4$ for $n \geq 3$. In the light
of Remark \ref{trivialstate/rem}, since in $K$ the link $L$ is already trivial, to
obtain $F_K$ we only need to thicken the undotted link components to $(1\
2)$-labelled ribbons with the right framings, then replace each dotted component
with a pair of parallel disks labelled by $(2\ 3)$, and finally connect by disjoint
bands the two $(1\ 2)$-labelled components and two of the $(2\ 3)$-labelled
disks, one for each pair.

\begin{Figure}[htb]{gompf/fig}{}{}
\centerline{\fig{}{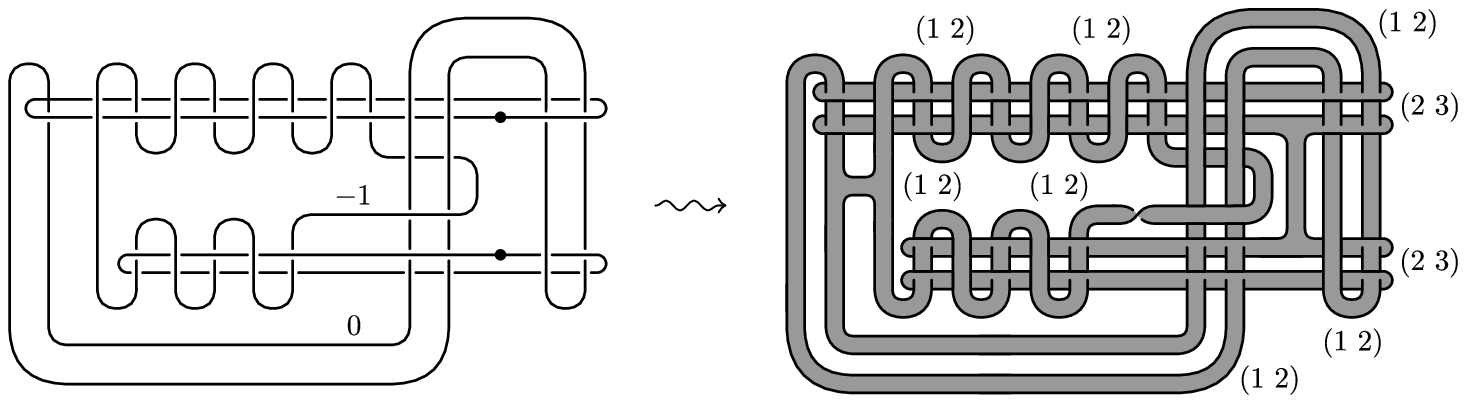}}
\end{Figure}

At this point, we are ready to prove our first equivalence theorem. Actually, this
is exclusively a matter of collecting the results we have already got.

\begin{proof}[Theorem \ref{equiv4/thm}]
The ``if'' part of the theorem is a special case of Propositions
\ref{movesR/2def/thm}. The ``only if'' part immediately follows from Propositions
\ref{to3fold/thm}, \ref{tospecial/thm}, \ref{diag/ribbon/diag/thm} and
\ref{2def/moves/thm}. Namely, given two labelled ribbon surfaces $F$ and $F'$
representing 2-equivalent 4-dimensional 2-handlebodies as branched covering of $B^4$
of the same degree $d \geq 4$, we can apply Propositions \ref{to3fold/thm} and
\ref{tospecial/thm} to transform them into $d$-stabilizations of certain $F_K$ and
$F_{K'}$, through moves $R_1$ and $R_2$. By the ``if'' part of the theorem and
Proposition \ref{diag/ribbon/diag/thm}, the two Kirby diagrams $K$ and $K'$ are
2-equivalent. Hence, by Proposition \ref{2def/moves/thm} $F_K$ and $F_{K'}$ are
related by moves $R_1$ and $R_2$.
\end{proof}

Before of going on to prove the other equivalence theorems, let us consider the
following proposition, which summarizes all we have said until now about branched
covering representation of (possibly disconnected) 4-dimensional 2-handlebodies.

\begin{proposition} \label{bijectivity/thm}
The map $F \mapsto K_F$ induces a bijective correspondence between 4-dimensional
2-handlebodies up to 2-deformation and labelled (orientable) ribbon surfaces,
representing them as simple branched coverings of $B^4$, up to labelling conjugation,
labelled 1-isotopy, stabilization and ribbon moves $R_1$ and $R_2$.
In particular, for handlebodies with $c$ connected components, the coverings can be
assumed to have degree $\leq 3c$ and two such coverings representations of
2-equivalent handlebodies can be related involving only coverings of degree $\leq
3c+1$.
\end{proposition}

\begin{proof}
By Proposition \ref{surjectivity/thm}, we already know that the correspondence in
the statement is surjective, being any 4-dimensional 2-handlebody with $c$ connected
components a $3c$-fold simple covering of $B^4$ branched over a ribbon surface (that
can be made orientable by Remark \ref{orient/rem}). 

To prove the injectivity, let us consider two coverings representing 2-equivalent
2-handlebodies. Since 2-deformation preserves connectedness, there is a bijective
correspondence between the components of the two handlebodies such that
corresponding components are 2-equivalent. Up to labelling conjugation, we can
assume that the sheets of the two coverings forming the corresponding components are
equally numbered. Moreover, by Proposition \ref{to3fold/thm} and destabilization, we
can reduce to 3 the maximum number of sheets for each component. Then, we can apply
Theorem \ref{equiv4/thm} to each pair of corresponding components in turn, leaving
unchanged the other ones. In this way, if the original coverings have degree $\leq
3c$, then all the intermediate coverings involved in relating them have degree $\leq
3c+1$.
\end{proof}

Having established our main result about branched covering representation of
4-dimensional 2-handlebodies, we pass to prove theorem \ref{equiv4b/thm}
concerning the case when they have diffeomorphic boundaries.

\begin{proof}[Theorem \ref{equiv4b/thm}]
As we observed in the Introduction, moves $P_\pm^{\pm1}$ and $T^{\pm1}$ do not change
the labelled boundary link up to labelled isotopy, so that they also preserve the
boundary of the covering manifold up to diffeomorphism.
Thus, taking into account Theorem \ref{equiv4/thm} and Proposition \ref{stmoves/thm}
\(b), we only need to show that such moves can be used to interpret blowing
up/down and 1/2-handle trading (cf. Figure \ref{diag4/fig}) for an ordinary Kirby
diagram $K$ in terms of the labelled ribbon surface $F_K$.\break Without loss of
generality, we can assume $K$ to be in standard form.

By definition of $F_K$, it is clear that moves $P_\pm$ obviously correspond positive
and negative blowups. Figure \ref{trading/fig} describes the sequence of moves
needed in order to replace the disk $C_i$ corresponding to the $i$-th 1-handle of
$K$ with the ribbon $A_{n+1}$ representing the new 2-handle deriving from the
trading.

\begin{Figure}[htb]{trading/fig}{}{}
\centerline{\fig{}{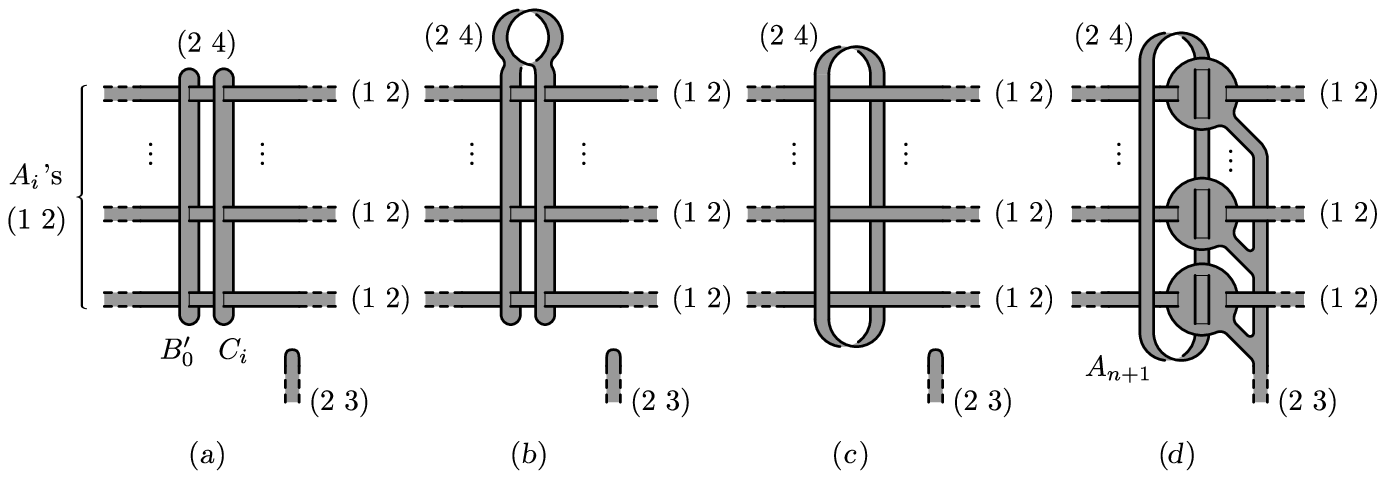}}
\end{Figure}

Diagram \(a) consists of the disk $C_i$, after we have operated on it as in Figure
\ref{indepDiag4/fig}, together with the parallel disk $B'_0$ resulting from that
modification. We perform a move $T$ and labelled 1-isotopy respectively to obtain
\(b) and \(c). This gives us an annulus representing the trivially framed attaching
loop of the new 2-handle. Then, we arrange the ribbon surface like in \(d),
according to a vertically trivial status of the framed link, by inserting some small
$(2\ 3)$-labelled disks as in \ref{indepOrder/fig}. Finally, we join the resulting 
annulus $A_{n+1}$ to $A_0$, by creating a new band $\alpha_{n+1}$ through a move
$R_3$, and we restore the stabilizing disk, as we did in the proof of Proposition
\ref{diag/ribbon/thm} by reversing the process of Figure \ref{indepAlpha/fig}.
\end{proof}

Our next goal is to derive Theorem \ref{equiv3/thm} from Theorem \ref{equiv4b/thm}.
The crucial point here is that any simply labelled link in $S^3$ can be transformed
through Montesinos moves into the boundary of a simply labelled ribbon surface in
$B^4$ (see Proposition \ref{link/ribbon/thm}). This follows quite directly from
Theorem B of \cite{MP01} about liftable braids, which we state here as Lemma
\ref{braids/thm} after having recalled a couple of definitions. 

\medskip

A simply labelled braid is called a {\sl liftable braid} when the two
labellings at its ends coincide. By an {\sl interval} we mean any braid that is
conjugate to a standard generator in the braid group. 
Actually, to make both the terms ``liftable'' and ``interval'' meaningful, one
should think of braids as self-homeomorphisms of the disk in the usual way (see
\cite{BW85} or \cite{MP01}), but this is not relevant in the present context. 

Of course, a labelled interval, as well as a standard generator, may or may not be
liftable depending on the labelling. We say that a labelled interval $x$ is of {\sl
type $i$} if $x^i$ is the first positive power of $x$ which is liftable. It is not
difficult to realize that conjugation preserves interval types and that each
interval is of type 1, 2 or 3 (cf. Lemma 2.4 of \cite{BW85} or Lemma 2.3 of
\cite{MP01}). 

The labelled intervals $x$, $y$ and $z$, whose first liftable positive powers are
depicted in Figure \ref{intervals/fig}, are the standard models for the three types
above. Namely, any labelled interval of type 1, 2 or 3 is respectively a conjugate of
$x^{\pm1}$, $y^{\pm1}$ or $z^{\pm1}$. Evidently, in the figure only the two
non-trivial strings of each labelled braid are drawn, the other ones being just
horizontal arcs with arbitrary labels. Moreover, in the labelling of each single
braid, we assume that $i$, $j$, $k$ and $l$ are all different.

\begin{Figure}[htb]{intervals/fig}{}{}
\centerline{\fig{}{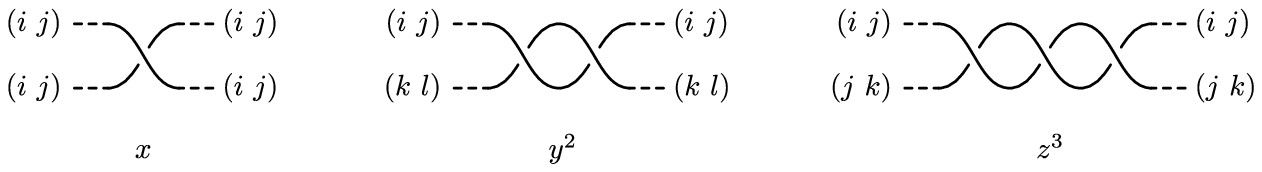}}
\end{Figure}

The main result of \cite{MP01} is the lemma below, which essentially says that any
liftable braid is a product of conjugates of labelled braids like the ones in Figure
\ref{intervals/fig}.

\begin{lemma} \label{braids/thm}
Any liftable braid is a product liftable powers of intervals.
\end{lemma}

We emphasize that the lemma holds without restrictions on the degree $d$ of the
labelling. However, it is worth observing that the case of $d = 2$ is trivial (every
braid is liftable in this case), while the case of $d = 3$ differs from the general
one for the absence of intervals of type 2. This special case was previously proved
in \cite{BW85} (cf. also \cite{BW94}), but the proof of Lemma \ref{braids/thm} given
in \cite{MP01} does not depend on \cite{BW85}.

\medskip

The relevant consequence of Lemma \ref{braids/thm} in the present context is the
following branched covering counterpart of the vanishing of the oriented cobordism
group $\Omega_3$.

\begin{proposition} \label{link/ribbon/thm}
Any labelled link $L \subset S^3$ representing a (possibly disconnected) $d$-fold
simple branched covering of $S^3$ is equivalent, up to labelled isotopy and moves
$M_1$ and $M_2$, to the boundary of labelled ribbon surface $F \subset B^4$
representing a $d$-fold simple branched covering of $B^4$.
\end{proposition}

\begin{proof}
Up to labelled isotopy, we can assume that the link $L$ is the closure $\widehat B$
of simply labelled braid $B$ (for example, we can use the labelled version of the well
known Alexander's braiding procedure). Of course, $B$ has to be a liftable braid.
Then, Lemma \ref{braids/thm} tells us that, up to labelled isotopy, we can think of
$B$ a product of conjugates of braids like $x^{\pm1}$, $y^{\pm2}$ or $z^{\pm3}$ (see
Figure \ref{intervals/fig}). Since braids $y^{\pm2}$ and $z^{\pm3}$ can be obviously
trivialized respectively by moves $M_2$ and $M_1^{\mp1}$, we can reduce ourselves to
the case when $B$ is a product of liftable intervals. 

In this case, a simply labelled ribbon surface $F \subset B^4$ bounded by $L$ can be
easily constructed from the band presentation of $B$ (see \cite{Ru83,Ru85})
determined by its factorization into liftable intervals. 
Namely, we start with a disjoint union of labelled trivial disks in $B^4$, spanned by
the labelled trivial braid obtained from $B$ by trivializing all the terms $x^{\pm1}$
appearing in the factorization above. Then, we attach to these disks a labelled
twisted band for each such term (see Figure \ref{covdiag7/fig} for a simple example,
where all the liftable intervals are standard generators).

Notice that the 3-dimensional diagram of the resulting surface may or may not form
ribbon intersection, depending on the conjugating braids of the liftable intervals
in the factorization of $B$ (cf. \cite{Ru83,Ru85}). In any case, the labelling
consistency when attaching the bands is ensured by the liftability of the intervals.
\end{proof}

\begin{proof}[Theorem \ref{equiv3/thm}]
As we said in the Introduction, it has been known for a long time, since the early
work of Montesinos, that moves $M_1$ and $M_2$ are covering moves. That is they, as
well as labelled isotopy and stabilization, do not change the covering manifold up
to diffeomorphism (see Section \ref{prelim/sec} for a proof of this fact).
Therefore, nothing more has to be added about the ``if'' part of the theorem.

The ``only if'' part follows from Proposition \ref{link/ribbon/thm} and Theorem
\ref{equiv4b/thm}, taking into account that the restriction of moves $R_1$ and $R_2$
to the boundary can be realized by moves $M_1$ and $M_2$ (see observation before of
Theorem \ref{equiv3/thm} in the Introduction), while moves $P_{\pm}$ and $T$
preserve the boundary up to labelled isotopy.
\end{proof}

Let us conclude this section with the proof of our last equivalence theorem. This
is Theorem \ref{equiv3g/thm}, which extends the previous Theorem \ref{equiv3/thm}
to possibly non-simple coverings of $S^3$ branched over an embedded graph.

\begin{proof}[Theorem \ref{equiv3g/thm}]
We have already observed in Section \ref{prelim/sec} that moves $S_1$ and $S_2$ are
covering moves, as they are applications of the coherent monodromies merging
principle. Hence, we have only to show that they allow us to transform any labelled
graph into a simply labelled link. We proceed in two subsequent steps: 1) we make the
labelling simple, by performing moves $S_1$ on the edges; 2) we make the graph into a
link, by performing moves $S_2$ on the vertices.

Let $G \subset R^3$ be a labelled embedded graph, endowed with a given graph
structure without loops (that is every edge has distinct endpoints).
We make the labelling simple, by operating on the edges of $G$ one by one. Each
time, we assume, up to labelled isotopy, that the edge $e$ under consideration is
not involved in any crossing. Denoting by $\sigma \in \Sigma_d$ the label of $e$, we
consider a coherent factorizations $\sigma = \tau_1 \dots \tau_k$ into
transpositions (any minimal factorization of $\sigma$ is coherent). Then, we split
$e$ into $k$ edges $e_1, \dots, e_k$ with the same endpoints, such that $e_i$ is
labelled by $\tau_i$, for each $i = 1, \dots, k$. To do that, we perform $k-1$
moves $S_1$, which progressively isolate the traspositions $\tau_i$ as labels of new
edges. Once all egdes of $G$ have been managed in this way, we are left with a simply
labelled graph which we still denote by $G$.

Now, we operate on the vertices of $G$ one by one, in order to make $G$ into a link.
Let $v$ be a vertex of $G$ and $e_1, \dots, e_h$ be the edges of $G$ having $v$ as
an endpoint, numbered according to the counterclockwise order in which they appear
around $v$ in the planar diagram of $G$. Since the total monodromy $\tau_1 \dots
\tau_h$ around $v$ must be trivial, $h$ must be even and the edges around $v$, can
be reodered, up to labelled isotopy, in such a way that $\tau_i = \tau_{h-i+1}$, for
every $i = 1, \dots, h/2$. This immediately follows from the well known
classification of the branched coverings of $S^2$, if one looks at a small 2-sphere
around $v$ transversal to $G$ (cf. \cite{BE79} or \cite{MP01}).
Then, by $h/2 - 1$ applications of move $S_2$, we replace the vertex $v$ by $h/2$
non-singular vertices $v_1, \dots, v_{h/2}$, such that $v_i$ is a common endpoint of
$e_i$ and $e_{h-i+1}$, for each $i = 1, \dots, h/2$.\break We leave to the reader to
verify that the sequence $\tau_1, \dots, \tau_{h/2}$ is coherent and that this
suffices for the needed moves $S_2$ to be performable. Obviously, after all the
singular vertices of $G$ have been replaced by non-singular ones, we are done.
\end{proof}

\medskip
\section{Final remarks\label{remarks/sec}}

First of all, we emphasize that the maps $F \mapsto K_F$ and $K \mapsto F_K$,
introduced respectively in Sections \ref{todiagram/sec} and \ref{toribbon/sec}
(see also Remark \ref{trivialstate/rem}), give an effective way to represent
4-dimensional 2-handlebodies up to 2-deformations as simple coverings of $B^4$
branched over ribbon surfaces, through generalized Kirby diagrams and Kirby calculus.

Effectiveness is preserved when passing to 3-manifolds too. In particular, being
the proof of Lemma \ref{braids/thm} in \cite{MP01} constructive, Proposition
\ref{link/ribbon/thm} and Theorem \ref{equiv3g/thm} (together with the map $F \mapsto
K_F$) enable us to define a procedure for obtaining a surgery description of a
closed orientable  3-manifold from any presention of it as a branched covering of
$S^3$ (cf. \cite{Hr01} and \cite{Hr03} for the 3-fold case).

Thus, it seems reasonable to expect recognition algorithms and effectively
computable invariants for closed orientable 3-manifolds (cf. \cite{Ma03}), based on
branched covering representation of them.

\smallskip

Secondly, we point out that our results, other than a different approach to covering
moves independent on \cite{Pi91}, \cite{Pi95} and \cite{Ap03}, also provide the
following new line of proof for the Hirsch-Hilden-Montesinos representation theorem:
start with the Alexander theorem \cite{Al20} to represent any closed oriented
3-manifold by a covering of $S^3$ branched over the 1-skeleton of a 3-simplex; make
such covering simple and its branching set into link, as in the proof of Theorem
\ref{equiv3g/thm}; apply Propositions \ref{link/ribbon/thm} and \ref{to3fold/thm}
in the order, to lower the degree of the covering.

Hopefully, the same ideas could be useful to make some progress in the branched
covering representation of smooth closed 4-manifolds. These are known to be 5-fold
simple coverings of $S^4$ branched over non-singular surfaces (see \cite{Pi95} and
\cite{IP02}), but it is an open problem whether the degree can be lowered from 5 to
4. Moreover, any result on covering moves relating diffeomorphic coverings of $S^4$
is still missing. Theorem \ref{equiv4/thm} together with the results of \cite{Mo77}
could give raise to a likely approach to this problem.

\medskip

Finally, we conclude with some remarks about the relation between  the present work
and some open problems in the topology of 4-dimensional 2-handlebodies.

\smallskip

A fundamental problem in 4-dimensional topology is the distinction between
homeomorphism and diffeomorphism classes. The only known invariants  able to 
detect such difference are the Donaldson and Seiberg-Witten invariants of smooth
closed manifolds. These have also been used (cf. \cite{Akb91} and Theorem 8.3.18 in
\cite{GS99}) to distinguish between homeomorphic but non diffeomeorphic 
4-dimensional 2-handlebodies in the cases when a standard way of closing them is
available. Invariants  defined directly on handlebodies and hopefully in purely
topological terms are missing and auspicable.

But there is even more delicate question which naturally arises in the topology of
4-dimensional 2-handlebodies and we have already mentioned in Section
\ref{prelim/sec}: is there difference between 2-equivalence classes and 
diffeomorphism classes? In \cite{Go91} (cf. the example after Remark
\ref{trivialstate/rem}) Gompf conjectures that the answer is yes, and offers a list
of possible counterexamples. Of course, detecting this phenomena can not rely any
more on invariants of smooth manifolds as the Seiberg-Witten ones.

Using the Hennings framework, in \cite{BM03} have been constructed invariants of
4-dimensional 2-handlebodies under 2-deformations, and these invariants depend on
the choice of an unimodular ribbon Hopf algebra. 

In the second part of this work \cite{BP05} we substantially improve this construction,
by showing that the map $F \rightarrow K_F$ between equivalence classes of ribbon
surfaces and Kirby diagrams factors through a bijective map onto the closed morphisms
of a universal category ${\cal H}^r$. The objects of ${\cal H}^r$ form a free
$(\otimes,1)$-algebra on a single object $H$, and $H$ is required to be a braided
ribbon Hopf algebra in ${\cal H}^r$. There is a standard procedure of ``braiding'' a
unimodular ribbon Hopf algebra $A$ associating to it a category ${\cal H}_A$ and a
functor ${\cal H}^r \rightarrow {\cal H}_A$. Therefore the invariants from \cite{BM03}
can be considered as particular examples of the new construction. But the result in
\cite{BP05} is much stronger: it actually gives a complete algebraic description of
4-dimensional 2-handlebodies and we hope is that this  would offer new approaches to
the open problems in the 4-dimensional topology mentioned above. Moreover, in
\cite{BP05} the result above is used to obtain an analogous algebraic description of
the boundaries of 4-dimensional 2-handlebodies, i.e. 3-dimensional manifolds, which
resolves for closed manifolds the problem posed by Kerler in \cite{Ke02} (cf.
\cite[Problem 8-16 (1)]{O02}).

\smallskip

Actually, the present work offers yet another possible approach towards studing the
difference between 2-deformations and diffeomorphisms: that is by relating it to the
difference between 1-isotopy and isotopy of ribbon surfaces.

We recall that 1-isotopy of ribbon surfaces in $B^4$ was derived from embedded
1-deformation of embedded 2-dimensional 1-handlebodies in $B^4$, by forgetting the
handlebody structure. On the other hand, once one has suitably defined embedded
2-deformation of embedded 2-dimensional 2-handlebodies, isotopy of arbitrary surfaces
in $B^4$ could be derived from it in a similar way.
Then, isotopy between ribbon surfaces differs from 1-isotopy just for allowing also
addition/deletion of embedded cancelling pairs of 1/2-handles and 2-handle isotopy.
This isotopy may involve non-ribbon intersections, such as double loops and triple
points, in the diagram.

In different words, we can say that two ribbon surfaces are 1-isotopic if and only if
they are isotopic through ribbon surfaces (of course, except for a finite number of
intermediate stages whose diagram is not self-tranversal). As we said in Section
\ref{prelim/sec}, we do not know whether isotopy relation between ribbon surfaces
coincides with 1-isotopy relation or not.

Anologously, since diffeomorphism of 4-dimensional 2-handlebodies is the same as
3-equivalence, the difference between 2-deformations and diffeomorphisms is in 
the addition/deletion of cancelling pairs of 2/3-handles and 3-handle isotopy. 
Moreover, the connection established in the previous sections, between labelled
1-isotopy of ribbon surfaces in $B^4$ and 2-deformation of 4-dimensional
2-handlebodies, through branched coverings and covering moves, can be at least
partially extended. More precisely, attaching a labelled 2-handle to the branching
surface $F \subset B^4$ corresponds to attaching a 3-handle to the covering
4-dimensional handlebody $H$, in such a way that any cancelling pair of 2/3-handles
of $H$ can be represented by a cancelling pair of labelled 1/2-handles of $F$.

A good staring point for studing this problem could be the example of the
Akbulut-Kirby  sphere $\Sigma_n$ (see figure \ref{gompf/fig} for the case of $n
= 3$). The proof given in \cite{Go91} of the fact that $\Sigma_n$ is diffeomorphic to
$B^4$, is based exactly on the intelligent introduction of a cancelling pair of
2/3-handles, changing the attaching maps by isotopy and eventually cancelling them
againts other handles. It would be interesting to see if, at least in this case,
these moves correspond to changing the branching surface by isotopy.

\end{document}